\documentclass[oneside]{amsart}

\usepackage[utf8]{inputenc}
\usepackage[T1]{fontenc}

\usepackage[a4paper]{geometry}
\usepackage{ragged2e}
\usepackage[pagebackref, hidelinks]{hyperref}
\hypersetup{
colorlinks= true, 
citecolor = blue, 
linkcolor = blue, 
urlcolor = blue, 
anchorcolor = blue
}

\usepackage{amsmath, amsthm, amsfonts}
\usepackage{amscd, amssymb, latexsym}
\usepackage[all]{xy}
\usepackage{mathtools} % for psmallmatrix

\usepackage{tabularx, array, diagbox}
\usepackage{tikz} %{pgf, pgfplots}
\usetikzlibrary{matrix, arrows.meta, positioning, quotes}

\usepackage{comment, color}
\usepackage{enumitem} 
\usepackage{soul}
\usepackage{caption}

\newtheorem{Definition}{Definition}[section]
\newtheorem{Theorem}[Definition]{Theorem}
\newtheorem{Lemma}[Definition]{Lemma}
\newtheorem{Proposition}[Definition]{Proposition}
\newtheorem{Corollary}[Definition]{Corollary}
\newtheorem{Example}[Definition]{Example}
\newtheorem{Remark}[Definition]{Remark}
\newtheorem{Question}[Definition]{Question}
\newtheorem{Conjecture}[Definition]{Conjecture}

\newcommand{\Ecf}[1]{\ensuremath{\lfloor\,{#1}\,\rfloor}}
\newcommand{\Card}[1]{\ensuremath{\left|{#1}\right|}}

\newcommand{\N}{\mathbb{N}}
\newcommand{\Z}{\mathbb{Z}}
\newcommand{\Q}{\mathbb{Q}}
\newcommand{\R}{\mathbb{R}}

\newcommand{\K}{\mathbb{K}}

\newcommand{\M}{\mathbb{M}}

\newcommand{\Primes}{\mathbb{P}}

\newcommand{\Sph}{\mathbf{S}}
\newcommand{\Proj}{\mathbf{P}}
\newcommand{\Torus}{\mathbb{T}}

\DeclareMathOperator{\Homeo}{Homeo}
\DeclareMathOperator{\Mod}{Mod}
\DeclareMathOperator{\Mat}{M}
\DeclareMathOperator{\GL}{GL}
\DeclareMathOperator{\SL}{SL}
\DeclareMathOperator{\PSL}{PSL}
\DeclareMathOperator{\PGL}{PGL}

\DeclareMathOperator{\Id}{\mathbf{1}}
\DeclareMathOperator{\val}{val}
\DeclareMathOperator{\im}{im}

\newcommand{\Al}{\mathcal{A}}
\DeclareMathOperator{\PP}{\mathcal{P}}
\DeclareMathOperator{\Shift}{\sigma}
\DeclareMathOperator{\CoInv}{\mathcal{G}}
\DeclareMathOperator{\Inf}{Inf}
\DeclareMathOperator{\Index}{\mathcal{I}}
\DeclareMathOperator{\Denjoy}{\mathcal{D}}

\DeclareMathOperator{\Gal}{\Gal}
\DeclareMathOperator{\covol}{\covol}

\title{Isogenies of minimal Cantor systems:\\ from Sturmian to Denjoy and interval exchanges}

\author{Scott Schmieding and Christopher-Lloyd Simon}
\date{}

\begin{document}

%\maketitle

\thanks{SS was partially supported by NSF grant DMS-2247553.}

\begin{abstract}
    This work is motivated by the study of continued fraction expansions of real numbers: we describe in dynamical terms their orbits under the multiplicative action by $\Z$ and the projective action of $\PGL_2(\Z)$, hence of $\PGL_2(\Q)$.
    A real number gives rise to a Sturmian system encoding a rotation of the circle, which is a particular case of an interval exchange system.
    It is well known that $\PGL_2(\Z)$-equivalence of real numbers, characterized by the tails of their continued fraction expansions, amounts to flow equivalence of Sturmian systems.
    We show that the multiplicative action of $m\in \Z$ on a real number corresponds to taking the $m$th-power followed by what we call an infinitesimal 2-asymptotic factor of its Sturmian system. 
    
    This leads us to introduce the notion of isogeny between zero-dimensional systems: it combines flow equivalences, virtual conjugacies and infinitesimal asymptotic equivalences.
    We develop tools for classifying systems up to isogeny involving cohomological invariants and states.
    We then use this to give a complete description of $\PGL_2(\Q)$-equivalence of real numbers in terms of Sturmian systems. Considering powers of Sturmian systems leads naturally to the realm of Denjoy systems, and we classify Denjoy systems up to isogenies preserving this class via the action of $\PGL_{2}(\Q)$ on their invariants.
    We also investigate eventual flow equivalence and show that it implies topological conjugacy for Sturmian and Denjoy systems.
    For this, we prove the following novel arithmetic result which we believe is of independent interest: Consider $\alpha,\beta\in \R\setminus\Q$; if for all but finitely many prime integers $n\in \Primes$ the numbers $n\alpha,n\beta$ are $\PGL_2(\Z)$-equivalent, then $\alpha \equiv \pm\beta \bmod{\Z}$.
    
    In another direction, we consider interval exchanges satisfying Keane's condition.  
    We characterize flow equivalence in terms of interval-induced subsystems (or the tails of their paths in the bilateral Rauzy induction diagram). 
    We describe their 2-AI factors in terms of minimal models. 
    Finally we find rational invariants for isogeny involving the length modules and SAF invariants of the associated ergodic measures.
    This leads to a conjecture for their classification up to isogeny, which we prove in the totally ergodic case.
    %
    % We believe that the isogeny classification of self-induced interval exchanges is related to the commensurator of its Veech group in $\PSL_2(\R)$.
    % \textbf{Keywords :} Sturmian systems, Denjoy systems, Interval exchang transformations, Flow equivalence, 
\end{abstract}

\maketitle

\renewcommand{\contentsname}{Plan of the paper}
\setcounter{tocdepth}{1}
\tableofcontents

%\newpage

\section{Introduction}

\subsection{Motivation: arithmetic of continued fraction expansions}

Every positive real number $\gamma$ admits a Euclidean continued fraction expansion:
\begin{equation*}
%\textstyle
\Ecf{c_0,c_1,\dots}=c_0+\frac{1}{c_1+\dots}
%\Ecf{a_0,a_1,a_2,\dots}=a_0+\frac{1}{a_1+\frac{1}{a_2+\dots}}
\quad \mathrm{with} 
\quad c_j \in \N
\quad \mathrm{and} 
\quad \forall j>0,\; c_j>0.
\end{equation*}
Such an expansion is infinite if and only if $\gamma$ is irrational, in which case it is unique. 
Moreover, this expansion is periodic if and only if $\gamma$ is quadratic.
% A longstanding conjecture \cite[Section 9]{Shallit_Reals-bounded-partial-quotients:survey_1992} states that if $c_i$ are bounded, then $\gamma$ must be either rational, quadratic or transcendental.

The group $\PGL_2(\Z)$ acting faithfully on $\R\Proj^1$ is generated by the inversion $\gamma\mapsto 1/\gamma$ and the translation $\gamma \mapsto \gamma+1$.
The positive powers of these two transformations generate the submonoid $\PGL_2(\N)$ preserving and acting freely on $(0,\infty)$.
Two numbers $\alpha, \beta\in (0,\infty)$ belong to the same orbit when some tails of their continued fraction expansions $\alpha_i= \Ecf{a_i,\dots}$ and $\beta_j= \Ecf{b_j,\dots}$ coincide.
In particular the set of rationals forms a single orbit, and the quadratic irrationals are grouped according to the periodic parts of their expansion. 
Hence the action of $\PGL_2(\N)\subset \PGL_2(\Z)$ on $(0,\infty)\subset \R\Proj^1$ is well understood, as its orbits correspond to those of the action of the monoid $\N$ on $(0,1)$ generated by the Gauss map $\gamma \mapsto \gamma_1$ which shifts the expansion $(c_i)\mapsto (c_{i+1})$.

The multiplicative action of $\N_{>0}\subset \Z^*$ on $(0,\infty)\subset \R$, which is very natural from the arithmetic viewpoint, is much more cryptic from the viewpoint of continued fraction expansions. Various works have studied the evolution of continued fractions under integer multiplication.

Indeed, it has recently been shown in \cite{Aka_Continued-fractions-quadratics_2020} (building on work in \cite{Aka-Shapira_Continued-fractions-quadratics_2018}) that for a fixed quadratic irrational $\gamma$, there is an infinite set of primes $\mathbb{P}_\gamma$ such that the continued fraction expansions along $\{q\gamma \colon q \in \mathbb{P}_\gamma\}$ approach the Gauss-Kuzmin statistics describing those of a random number.
Under the generalised Riemann hypothesis, a similar result is shown where $\mathbb{P}_\gamma$ is a full density set of primes, and even a full density set of natural numbers.
Thus, while the continued fraction expansions of $\gamma$ and $q\gamma$ determine one another, they become de-correlated in the limit.
This reveals that the actions of $\N_{>0}$ by multiplication and by the Gauss map are somewhat asymptotically decorrelated.

The multiplicative action of $\Q^*$ and the action of $\PGL_2(\Z)$ on $\R\Proj^1$ generate the action of $\PGL_2(\Q)$. Note that all these actions preserve (the projective line over) each real field $\K \subset \R$, in particular the rationals, consisting of a single orbit, and each real quadratic field which partitions into exactly two orbits.

\begin{comment}
% \begin{Remark}[Low complexity: symbolic vs arithmetic]
The rational and quadratic numbers are the only ones which are simple both from the symbolic viewpoint of continued fraction expansions and the arithmetic viewpoint of algebraic number theory.
%
On the one hand, a longstanding conjecture \cite[Section 4]{Shallit_Reals-bounded-partial-quotients:survey_1992} states that if $c_j$ is bounded then $\gamma$ must be either rational, quadratic or transcendental. 

On the other hand, McMullen \cite{McMullen_diophantine-quadratic-field_2009} shows that for every square-free $\delta\in \N_{>1}$, there is a bound $M_{\delta}\in \N_{>1}$ such that $\Q(\sqrt{\delta})$ contains infinitely many quadratic irrational numbers whose continued fractions expansions have entries $\le M_{\delta}$. (It is also conjectured that one may always choose $M_{\delta}=2$.)
% \end{Remark}
\end{comment}

\subsection{From continued fractions to Sturmian systems}

The diophantine approximation properties of $\gamma \in (0,1)$ are encoded by its sequences of continuants $c_j\in \N_{\ge 1}$ for $j\in\N_{>0}$ and of fractional parts $q\gamma-\lfloor q\gamma \rfloor \in [0,1)$ for $q\in \N_{>0}$.

This leads to the dynamical setting: consider the torus $\Torus = \R\bmod{\Z}$ together with the rotation $R_\gamma \colon \Torus \to \Torus$ defined by $t \mapsto t+\gamma\bmod{1}$ which under the identification $\Torus=[0,1)$ implements a piecewise translation on the two intervals $I_0= [0, 1-\gamma)$ and $I_1 = [1-\gamma, 1)$. 
The orbit $R_\gamma^q(0)$ of $0$ recovers the fractional parts of $q\gamma$.
Recording the sequence of intervals containing the points $R_{\gamma}^{q}(0)$ yields the Sturmian sequence $x_{\gamma} \in \{0,1\}^{\Z}$ associated to $\gamma$.

Geometrically, the line \(\R\cdot \begin{psmallmatrix} \gamma \\ 1 \end{psmallmatrix}\) in the Cartesian plane $\R^2$ intersects the integral translates of vertical and horizontal coordinate lines according to the bi-infinite sequence $x_\gamma \in \{0,1\}^{\Z}$.
The generators of the Euclidean monoid $J\colon \gamma\mapsto 1/\gamma$ and $R \colon \gamma\mapsto \gamma+1$ acting on real numbers $\gamma \in (0,\infty)$ lift to an action on words over $\{0,1\}$ via the substitutions defined by $J(0)=1, J(1)=0$ and $R(0)=01, R(1)=1$, and the sequence $x_{\gamma} \in \{0,1\}^\Z$ can then be derived from the continued fraction expansion. Conversely, one can recover the continued fraction expansion of $\gamma$ using the S-adic expansion of its Sturmian sequence $x$.

The compact space $\{0,1\}^\Z$ is endowed with the shift homeomorphism $\Shift(u)_{n} = u_{n+1}$, and the closure of the orbit of our Sturmian sequence $x_{\gamma}$ gives a subshift $X_\gamma \subset \{0,1\}^\Z$ called the Sturmian subshift associated to the parameter $\gamma$.
Hence the one-dimensional dynamical system given by the rotation $R_\gamma\colon \Torus\to \Torus$ is encoded by the zero-dimensional system given by the Sturmian subshift $\Shift_\gamma \colon X_{\gamma} \to X_{\gamma}$.

One connection between the arithmetic of real irrationals and the dynamics of Sturmian systems occurs via the flow equivalence relation. Given a system $(X,S)$, its suspension is the space $\Sigma_{S}X = (X\times \R)\bmod{\Z}$ where the diagonal action of $n\in \Z$ is by $(x,s)\to (S^n(x), s-n)$. The suspension comes with the continuous flow $\Phi_{s'} \colon (x,s)\mapsto (x,s+s')$, and we say two systems $(X,S)$ and $(Y,T)$ are \emph{flow equivalent} when there is a homeomorphism between their suspension spaces $\Sigma_{S} X$ and $\Sigma_{T} Y$ which preserves the leaves of the flow and their orientation induced by the flow. 

% Fokkink~\cite{Fokkink_structure-trajectories_1991} established the following correspondence between arithmetic and dynamics: two real irrationals $\alpha, \beta\in \R\setminus \Q$ are $\PGL_{2}(\Z)$-equivalent if and only if the associated Sturmian systems $(X_\alpha,\Shift_\alpha)$ and $(X_\beta,\Shift_\beta)$ have homeomorphic suspensions. (Note that Fokkink uses different language, referring to the suspension of Sturmian subshifts as Denjoy continua.) It follows immediately from this (see Theorem~\ref{thm:flow-equivalence-Sturmian-systems} for details) that Sturmian subshifts with parameters $\alpha,\beta$ are flow equivalent if and only if $\alpha$ and $\beta$ are $\PGL_{2}(\Z)$-equivalent.

Our work is motivated by the following correspondences between the arithmetic of real irrationals $\alpha, \beta\in \R\setminus \Q$ and the dynamics of their associated Sturmian systems $(X_\alpha,\Shift_\alpha), (X_\beta,\Shift_\beta)$.
First, it is easy to see that $(X_\alpha,\Shift_\alpha)$ and $(X_\beta,\Shift_\beta)$ are topologically conjugate if and only if $\alpha= \pm \beta \bmod{\Z}$ (equal up to sign and $\Z$-addition).
Second, Fokkink~\cite{Fokkink_structure-trajectories_1991} showed that $(X_\alpha,\Shift_\alpha)$ and $(X_\beta,\Shift_\beta)$ are flow equivalent if and only if $\alpha, \beta$ are $\PGL_{2}(\Z)$-equivalent (we reprove this in Theorem~\ref{thm:flow-equivalence-Sturmian-systems} using our modern methods involving states and coinvariants).
% This characterisation up to flow equivalence goes back to Fokkink~\cite{Fokkink_structure-trajectories_1991} (who used different language, referring to the suspension of Sturmian subshifts as Denjoy continua), and we will explain it in Theorem~\ref{thm:flow-equivalence-Sturmian-systems} using our modern methods involving states and coinvariants.

In light of this result and the discussions in the previous paragraphs, one of our main motivations is to describe the projective action of $\PGL_2(\Q)$ on $\R\setminus \Q$ in terms of Sturmian systems.

The action of $\Z$ by multiplication on irrationals leads us to consider powers of systems and introduce the notion of \emph{infinitesimal 2-asymptotic factors} of systems.
Indeed, we will see that the Sturmian system $(X_{m\alpha}, \Shift_{m\alpha})$ is an infinitesimal 2-asymptotic factor of the $m$-th power of the Sturmian system $(X_\alpha, \Shift_\alpha^m)$.
However, the $m$-th power $(X_\alpha, \Shift_\alpha^m)$ is no longer Sturmian.
We are thus naturally led to extend the class of Sturmian systems (particular encodings of rotations on the circle) to the class of Denjoy systems (Cantor extensions of Denjoy homeomorphisms of the circle).

With this motivation in mind, we now introduce a few definitions in more detail to give an account of our main constructions and results, including some applications to the Sturmian, Denjoy and interval exchange systems.

Section~\ref{sec:systems} introduces the main definitions and running examples of dynamical systems, which will be discussed in the applications ending each of the following sections.

\subsection{Isogenies between dynamical systems}

By a (dynamical) \emph{system} $(X,S)$ we mean a non-empty compact topological space $X$ together with a homeomorphism $S \colon X\to X$.
% A factor map between two systems $(X,S)$ and $(Y,T)$ is a surjective function $\pi \colon (X,S) \to (Y,T)$ such that $\pi \circ S = T\circ \pi$. 
% Two systems $(X,S)$ and $(Y,T)$ are \emph{topologically conjugate} when there exists $h \in \Homeo(X,Y)$ such that $T = h\circ S \circ h^{-1}$

For a system $(X,S)$ and $m\in \N_{>0}$, its $m$-th power is the system $(X, S^m)$, where as usual $S^{m}$ means $S$ composed with itself $m$ times.
Given an equivalence relation $\mathcal{R}$ on systems, we say that two systems $(X,S)$ and $(Y,T)$ are \emph{virtually $\mathcal{R}$} when there exist $m,n\in \N_{>0}$ such that the systems  $(X,S^m)$ and $(Y,T^n)$ are $\mathcal{R}$.
% We say that they are \emph{eventually $\mathcal{R}$} when there exists $m\in \N_{>0}$ such that for all $n>m$ the systems $(X,S^n)$ and $(Y,T^n)$ are $\mathcal{R}$.

The suspension of a system $(X,S)$ consists of its mapping torus $\Sigma_{S}X = (X\times \R)\bmod{\Z}$ where the diagonal action of $n\in \Z$ is by $(x,s)\to (S^n(x), s-n)$, together with the continuous flow $\Phi_{s'} \colon (x,s)\mapsto (x,s+s')$.
Two systems $(X,S)$ and $(Y,T)$ are called \emph{flow equivalent} when their suspension spaces $\Sigma_{S} X$ and $\Sigma_{T} Y$ are homeomorphic by a homeomorphism respecting the orientation induced by the flow.

We must finally introduce the notion of \emph{infinitesimal 2-asymptotic (2-AI) factor} and the corresponding \emph{infinitesimal 2-asymptotic (2-AI) equivalence} relation. This requires a lengthy and technical discussion involving the group of coinvariants of a system (Section~\ref{sec:coinvariants}) and the index map of a 2-asymptotic factor (Section~\ref{sec:2-AI-equivalence}), so we take it as a black box for this introduction.
Let us just mention one useful characterization that we prove in Corollary~\ref{cor:2-asym-factors-ue-minimal-cantor}. To a zero-dimensional system $(X,S)$ is associated a unital ordered group of coinvariants $\CoInv_S$, and each invariant Borel probability measure $\mu$ on the system defines a state on that group $\tau_{\mu}\colon \CoInv_S \to \R_{\ge 0}$ given by integration against that measure.
For uniquely ergodic minimal Cantor systems, the existence of a 2-AI factor map between the systems implies that the images in $\R$ of their respective states coincide after tensoring with the rationals. Moreover, for such systems, if the rational images of the states coincide then every 2-asymptotic factor map between them must be infinitesimal.

Thus 2-AI equivalence is imposed on us after considering powers of systems (as in virtual topological conjugacy). 
Indeed, we prove in Corollary~\ref{cor:2-AI-factor-Denjoy-power-m} that the Sturmian system $(X_{m\alpha}, \Shift_{m\alpha})$ and the Denjoy system $(X_\alpha, \Shift_\alpha^m)$ are 2-AI equivalent. 

\begin{comment}
\begin{Proposition}
    The Sturmian system $(X_{m\alpha}, \Shift_{m\alpha})$ and the Denjoy system $(X_\alpha, \Shift_\alpha^m)$ are 2-AI equivalent.
\end{Proposition}
\end{comment}
\begin{comment}
Let $(X,S)$ be a system. 
Two points $x,x^{\prime} \in X$ are called asymptotic, denoted $x\Bumpeq x'$, when $d(S^{n}(x),S^{n}(x^{\prime})) \to 0$ as $n \to \pm \infty$.
% This is an equivalence relation.
% Note that for a factor map $\pi \colon (X,S)\to (Y,T)$ we have $x\Bumpeq x' \implies \pi(x)\Bumpeq \pi(x')$.
A factor map $\pi \colon (X,S) \to (Y,T)$ is called \emph{2-asymptotic} (2-A) when there exists a subset of \emph{blown-up points} $B \subset Y$ such that
\begin{enumerate}[noitemsep]
\item For $y\in Y$, we have $\Card{\pi^{-1}(y)}> 1$ if and only if $y \in B$.
\item For every $y\in B$, we have $\Card{\pi^{-1}}(y) = 2$.
\item For every $y \in B$, the pair of points in $\pi^{-1}(y)=\{x,x'\}$ are asymptotic.
\end{enumerate}

The definition of \emph{infinitesimal} 2-A factors requires the group of co-invariants.
\end{comment}
We are now ready to formulate one of the main definitions in this work (Definition~\ref{def:isogeny}).
\begin{Definition}[isogeny]
    We say that two systems are \emph{isogenous} when they are related by a sequence of virtual flow equivalences and 2-AI equivalences.
\end{Definition}

One of our first main results (Theorem~\ref{thm:isogenous-Sturmian}) characterizes the action of $\PGL_{2}(\Q)$ on $\R \setminus \Q$ in terms of isogeny of Sturmian systems.

\begin{Theorem}[isogenies of Sturmian systems]
    The following are equivalent:
\begin{itemize}[noitemsep]
    \item The points $\alpha,\beta\in \R\Proj^1\setminus \Q\Proj^1$ belong to the same $\PGL_{2}(\Q)$-orbit.
    \item The Sturmian subshifts $(X_{\alpha},\Shift_{\alpha})$, $(X_{\beta},\Shift_{\beta})$ are isogenous.
\end{itemize}
\end{Theorem}

Our methods also enable us to characterise isogenies within the class of Denjoy systems, which generalize Sturmian systems according to the following construction. 
A homeomorphism of the circle $F\colon \Torus \to \Torus$ with irrational rotation number $\alpha$ which is not conjugate to the rotation $R_\alpha$ is uniquely characterized up to topological conjugacy by its action on its unique minimal set $\Denjoy_F$ which is a Cantor set.
The restriction of $F$ to $\mathcal{D}_{F}$ yields a Denjoy system $(\Denjoy_{F}, F)$ which is characterised up to conjugacy by its invariants: the rotation number $\alpha$ and a countable set of points $A\subset \Torus$ (see Subsection~\ref{subsec:denjoy-systems}).

We now state the characterisation of Denjoy systems up to isogeny achieved in Theorem~\ref{thm:isogenous-Denjoy}.

\begin{Theorem}[isogenies of Denjoy systems]
For $k \in \{0,1\}$, consider Denjoy systems $(\Denjoy_{F_k}, F_k)$ with invariants $(\alpha_k, A_k)$.
\begin{enumerate}
\item
If there exists 
$M= (m_{ij}) \in \PGL_2(\Q)$ such that:
\begin{equation*}
    \alpha_1 = \frac{m_{11}\cdot\alpha_0+m_{12}}{m_{21}\cdot\alpha_0+m_{22}} 
    \qquad \mathrm{and} \qquad 
   \Q A_1 \equiv \frac{1}{m_{21}\cdot\alpha_0+m_{22}} \cdot \Q A_0
\end{equation*}
then the Denjoy systems $(\Denjoy_{F_k}, F_k)$ are isogenous.
\item
If the Denjoy systems are isogenous through a chain of Denjoy systems, then there exists 
$M= (m_{ij}) \in \PGL_2(\Q)$ such that:
\begin{equation*}
    \alpha_1 = \frac{m_{11}\cdot\alpha_0+m_{12}}{m_{21}\cdot\alpha_0+m_{22}} 
    \qquad \mathrm{and} \qquad 
    \Q A_1 \equiv \frac{1}{m_{21}\cdot\alpha_0+m_{22}} \cdot \Q A_0.
\end{equation*}
\end{enumerate}
\end{Theorem}

\subsection{Generalizations to interval exchange systems}

We also discuss interval exchange systems (symbolic encodings of interval exchange transformations, abbreviated IES), which also generalise Sturmian systems (encoding irrational rotations).
More precisely, an IES on $d$ symbols is a subshift obtained from the canonical encoding of an IET on $d$ intervals. 
Given an IES $(X,\Shift)$, its ergodic probability measures correspond (via length of segments) to the IETs $T\colon (0,1]\to (0,1]$ with that encoding.

We focus on those which generalize irrational rotations, namely the IETs satisfying Keane's infinite distinct orbits (ido) condition. This is equivalent to saying that the associated IES has complexity $n\mapsto (d-1)n+1$, and also to the fact that it is infinitely many times renormalizable via Rauzy-Veech induction (the analogue of continued fraction expansions).

\begin{comment}
The fact that a real number is irrational if and only if its continued fraction expansion is infinite generalizes to the fact that an IES has ido if and only if its Rauzy induction is infinite.

The fact that an irrational is badly approximable by rationals if and only if its continued fraction expansion has bounded entries generalizes to the fact that an IES with ido satisfies Boshernitzan's criterion if and only if its Rauzy induction yields a primitive proper S-adic representation.

The fact that real numbers with infinite periodic continued fraction expansions correspond to slopes of foliations of the torus which are fixed by elements of the mapping class group generalizes to the fact that IES with ido and infinite periodic Rauzy induction correspond to the attractive foliations of pseudo-Anosov mapping class of translation surfaces.
\end{comment}

We prove an analog to Fokkink's theorem characterising flow equivalence of IES in terms of two-sided Rauzy induction~\ref{thm:flow-equivalence-IES}.
We describe 2-A factors of IES in subsection~\ref{subsec:IES-2A-minimal-model}, by providing a preferred representative (related to the minimal models in \cite[Section 2.1]{Dahmani_Groups-IET_2019}).

For an IES $(X,\Shift)$, its invariant measures correspond to the IETs $T\colon (\alpha_0,\alpha_d]\to (\alpha_0,\alpha_d]$ with that symbolic encoding, which form a cone in $\R_+^d$.
Let us recall that the IETs form a group under composition.
Arnoux-Fathi introduced a surjective homomorphism $\operatorname{SAF}\colon \operatorname{IET}\to \R\wedge_\Q \R$, which is a kind of scissor-congruence invariant: they showed that it is invariant by inducing IETs on smaller intervals (see~\cite[Proposition 2.3]{Veech_metric-theory-IET-2-ApproxByPrimitive_1984}).
Then Sah showed that $\operatorname{SAF}\colon \operatorname{IET}\to \R\wedge_\Q \R$ corresponds to the abelianization (see~\cite[Theorem 1.3]{Veech_metric-theory-IET-3-SAF_1984}).

\begin{Definition}[rational invariants]
Consider an IES $(X,\Shift)$ with ido and let $\mu_1,\dots,\mu_e$ be its $e\in \N_{>0}$ ergodic probability measures. Each $\mu_i$ corresponds to an IET $T_{\mu_i}\colon (0,1]\mapsto (0,1]$ with that encoding, to which we associate the $\Q$-vector space $\Q\otimes \im(\tau_{\mu_i})\subset \R$ of dimension at most $d-1$ generated by the $\mu_i$-measures of cylinders in $X$, and the vector $\operatorname{SAF}(T_{\mu_i})\in \R\wedge_\Q \R$.
    
We define the \emph{rational invariants} of $(X,\Shift)$ 
as the product of subspaces $\prod_i \Q\otimes \im(\tau_{\mu_i})\subset \R^e$ up to multiplication by $\lambda \in \R_{>0}$ on $\R^e$ together with the product of vectors $\prod_i \operatorname{SAF}(T_{\mu_i})\in (\R\wedge_\Q \R)^e$ up to multiplication by $m/n\in \Q_{>0}$.
\end{Definition}

We propose the following Conjecture~\ref{conj:isogenies-IES=Rational-Invariants}.

\begin{Conjecture}[rational invariants of isogenies]
    If two IES with ido are isogenous then they have the same rational invariants.
\end{Conjecture}

We prove this conjecture in the totally ergodic case for isogenies staying within the class of IESs, in Proposition~\ref{prop:isogenies-total-ergodic-IES}.

\begin{Proposition}[isogenous totally ergodic IES have same rational invariants]
    If totally ergodic IES with ido are related by isogenies through IES, then they have the same rational invariants.
\end{Proposition}

We will see in Corollary~\ref{cor:isogenies-IES-ido-LR} that this proposition applies to IES with ido which are linearly recurrent, in particular (Corollary~\ref{cor:isogenies-self-induced-IES}) to those which are self-induced, and hence  (Corollary~\ref{cor:quadratic-IES}) to the minimal IES with lengths in a quadratic field (characterized by \cite{Boshernitzan-Carroll_quadratic-IET_1997} as satisfying the generalisation of Lagrange's theorem on the periodicity of continued fractions).

Self-induced minimal IES are precisely those arising from the stable laminations of pseudo-Anosov automorphisms of translation surfaces: the isogeny relation thus yields an equivalence relation between such pseudo-Anosov automorphisms which we believe is worth studying.
% Those yield a Teichm\"uller ray in the space of translation surfaces, 
In particular, we propose Conjecture~\ref{conj:veech-groups} relating isogenies of minimal self-induced IES and commensurators of their Veech groups (the automorphism groups of their underlying translation surfaces in $\PSL_2(\R)$), supported by the fact that isogenies preserve their rational invariants.
\begin{table}[h]
    \centering
    \begin{tabular}{||c||c|c|c||}
         \hline 
         Equivalence 
         & Sturmian: $\alpha\in \Torus$ 
         & Denjoy $\alpha \in A\subset \Torus$ 
         & IES with ido 
         %& SFT: $M\in \Mat_n(\N)$ 
         \\ \hline \hline
         Conjugacy 
         & $\pm \alpha \bmod{\Z}$ 
         & $\pm (\alpha, A) \bmod{\Z}$ 
         & ??? 
         \\
         & Lemma~\ref{lem:equal-rank2-lattices} 
         & Markley : §\ref{subsec:denjoy-systems} 
         & ???
         \\ \hline
         Flow 
         & $\PGL_2(\Z)\cdot \alpha$ 
         & $\PGL_2(\Z)\cdot (\alpha, A)$ 
         & Rauzy tail
         \\
         & Theorem~\ref{thm:flow-equivalence-Sturmian-systems} 
         & Theorem~\ref{thm:flow-equivalence-Denjoy-systems} 
         & Theorem~\ref{thm:flow-equivalence-IES}
         \\ \hline
         Virtual 2-AI 
         & $\Q \cdot \alpha \bmod{\Z}$ 
         & $(\Q \alpha,\Q A) \bmod{\Z}$ 
         & ??? 
         \\
         & Cor~\ref{cor:2-AI-equiv-Denjoy} \&~\ref{cor:2-AI-factor-Denjoy-power-m} 
         & Cor~\ref{cor:2-AI-equiv-Denjoy} \&~\ref{cor:2-AI-factor-Denjoy-power-m} 
         & ??? 
         \\ \hline 
         Isogeny & $\PGL_2(\Q)\cdot \alpha$ & $\PGL_2(\Q)\cdot (\alpha, A)$ & $\Q$-invariants
         \\
         (preserving type) 
         & Theorem~\ref{thm:isogenous-Sturmian} 
         & Theorem~\ref{thm:isogenous-Denjoy} 
         & Conj~\ref{conj:isogenies-IES=Rational-Invariants} \& Prop~\ref{prop:isogenies-total-ergodic-IES}
        \\ \hline
    \end{tabular}
    \caption{Summary of the classification results towards isogeny}
    \label{tab:isogenies}
\end{table}

Let us conclude the discussion on isogenies with Table~\ref{tab:isogenies} summarizing the classification results achieved in this work.
The notion of isogeny arose from our questions about Sturmian systems and $\PGL_2(\Q)$ equivalence, for which we achieved the desired classification results.
We extended this classification to Denjoy systems and certain IES using rational invariants, in restriction to isogenies that remain in this class. 
To understand the structure of isogenies in general (not restricted to a certain class of systems), we were led to ask Questions
\ref{quest:factoring-asymptotic-orbits} and
\ref{quest:terminal-2-A-factors}.

\subsection{Eventual flow equivalence}
Given an equivalence relation $\mathcal{R}$ on systems, we say that two systems $(X,S)$ and $(Y,T)$ are \emph{eventually $\mathcal{R}$} if for all but finitely many $n \in \N$, the systems $(X,S^{n})$ and $(Y,T^{n})$ are $\mathcal{R}$-equivalent. 
%%there exists $m\in \N_{>0}$ such that for all $n>m$ 
%the systems $(X,S^n)$ and $(Y,T^n)$ are $\mathcal{R}$.
In general, we have the following implications:
% (Since topological conjugacy implies flow equivalence, it follows that that eventual topological conjugacy implies eventual flow equivalence and virtual conjugacy implies virtual flow equivalence.)
%
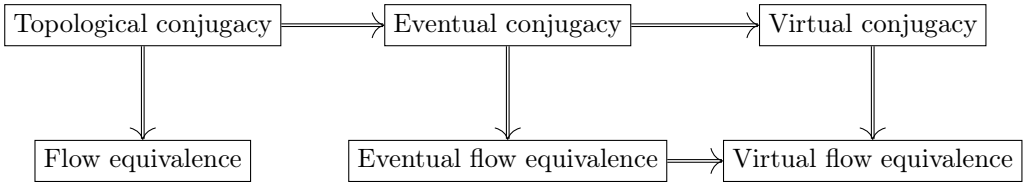
\begin{center}
\begin{tikzpicture}[node distance=4.3cm and 1cm,
box/.style = {draw, align=center}]

    % Define nodes
    \node (A) [box] at (0, 0) {Topological conjugacy};
    \node (B) [box, right of=A, xshift=.5cm] {Eventual conjugacy};
    \node (C) [box, right of=B, xshift=.5cm] {Virtual conjugacy};
    \node (D) [box, below of=A,yshift=2.5cm] {Flow equivalence};
    \node (E) [box, below of=B,yshift=2.5cm] {Eventual flow equivalence};
    \node (F) [box, below of=C,yshift=2.5cm] {Virtual flow equivalence};

    % Draw arrows
    \draw[->, double] (A) -- (B);
    \draw[->, double] (B) -- (C);
    \draw[->, double] (A) -- (D);
    \draw[->, double] (B) -- (E);
    \draw[->, double] (C) -- (F);
    % \draw[->, double] (D) -- (E);
    \draw[->, double] (E) -- (F);
    % \draw[->, dots] (E) -- (A);
\end{tikzpicture}
\end{center}

The notion of eventual topological conjugacy and its connection with shift equivalence of nonnegative square matrices lies at the heart of Williams' Problem and the classification of shifts of finite type (as explained in the survey~\cite{Boyle-Schmieding_Stable-Algebra_2024}). 
Furthermore, for shifts of finite type, the relation between eventual topological conjugacy and flow equivalence had been studied in~\cite{Boyle_algebraic-aspects-symbolic-dynamics_2000}, and has recently been reexamined in~\cite{Boyle_shift-equiv-implies-flow-equiv_2024}, where it was shown that eventual conjugacy implies flow equivalence.

We focus on eventual flow equivalence, and show that it implies topological conjugacy for Sturmian systems (Theorem~\ref{thm:eventual-flow-Sturmian}) and Denjoy systems (Theorem~\ref{thm:eventual-flow-Denjoy}).

\begin{Theorem}[eventual flow equivalence for Sturmian systems]
\label{Intro_thm:eventual-flow-Sturmian}
Two Sturmian systems 
% $(X_{\alpha},\Shift_{\alpha})$ and $(X_{\beta},\Shift_{\beta})$ 
are eventually flow equivalent if and only if they are topologically conjugate.
% (namely $\alpha=\pm\beta\bmod{\Z}$).
\end{Theorem}

\begin{Theorem}[eventual flow equivalence for Denjoy systems]
\label{Intro_thm:eventual-flow-Denjoy}
Two Denjoy systems are eventually flow equivalent if and only if they are topologically conjugate.
\end{Theorem}

It follows from the previous theorem that even for Sturmian subshifts, flow equivalence does not imply eventual flow equivalence. This provides another motivation to study our notion of isogeny. 
We also give some examples (Examples~\ref{ex:virtual-flow-not-flow},~\ref{ex:virtual-flow-not-flow-quadratic}) of Sturmian subshifts which are not flow equivalent, but have powers which are flow equivalent.

The following diagrams summarize the various implications between virtual and eventual equivalencies (topological conjugacy and flow) that we obtain for Sturmian systems (Subsection \ref{subsec:eventual-flow-Sturmian}) and Denjoy systems (Subsection \ref{subsec:eventual-flow-Denjoy}).

\begin{center}
\begin{tikzpicture}[node distance=4.5cm and 1cm,
box/.style = {draw, minimum height=12mm, inner xsep=1mm, align=center},
%sy+/.style = {yshift= 2mm},
%sy-/.style = {yshift=-2mm},
every edge quotes/.style = {align=center}]
    % Define nodes
    \node (A) [box] at (0, 0) {Topological conjugacy \\ $\pm \alpha \bmod{\Z}$ \\ Lemma~\ref{lem:equal-rank2-lattices}};
    \node (B) [box, right of=A, xshift=.5cm] {Eventual conjugacy \\ $\pm \alpha \bmod{\Z}$ \\ Proposition~\ref{prop:eventual-conj-to-conj-Sturmian}};
    \node (C) [box, right of=B, xshift=.5cm] {Virtual conjugacy \\ $\pm \alpha \bmod{\Q}$ \\ Proposition~\ref{prop:virtual-conjugacy-Sturmian}};
    \node (D) [box, below of=A, yshift=2.5cm] {Flow equivalence \\ $\alpha \bmod{\PGL_2(\Z)}$ \\ Theorem~\ref{thm:flow-equivalence-Sturmian-systems}};
    \node (E) [box, below of=B,yshift=2.5cm] {Eventual flow equivalence \\ $\pm \alpha \bmod{\Z}$ \\ Theorem~\ref{thm:eventual-flow-Sturmian}};
    \node (F) [box, below of=C, yshift=2.5cm] {Virtual flow equivalence \\ $n\alpha\bmod{\PGL_2(\Z)}$ \\ Remark \ref{rem:virtual-flow-equivalence-Sturmian}};

    % Draw arrows
    \draw[->, double] (A) -- (B);
    \draw[->, double] (B) -- (C);
    \draw[->, double] (A) -- (D);
    \draw[->, double] (B) -- (E);
    \draw[->, double] (C) -- (F);
    % \draw[->, double] (D) -- (E);
    \draw[->, double] (E) -- (F);
    \draw[->, double] (B) -- (A);
    \draw[->, double] (E) -- (B);
    \draw[->, double] (E) -- (A);
    \draw[->, double] (A) -- (E);
\end{tikzpicture}
\end{center}

\begin{center}
\begin{tikzpicture}[node distance=4.5cm and 1cm,
box/.style = {draw, minimum height=12mm, inner xsep=1mm, align=center},
%sy+/.style = {yshift= 2mm},
%sy-/.style = {yshift=-2mm},
every edge quotes/.style = {align=center}]
    % Define nodes
    \node (A) [box] at (0, 0) {Topological conjugacy \\ $\pm (\alpha, A) \bmod{\Z}$ \\ Lemma~\ref{lem:equal-rank2-lattices}};
    \node (B) [box, right of=A, xshift=.5cm] {Eventual conjugacy \\ $\pm (\alpha,A) \bmod{\Z}$ \\ Theorem \ref{thm:eventual-flow-Denjoy} };
    \node (C) [box, right of=B, xshift=.5cm] {Virtual conjugacy \\ $(\pm\alpha \bmod{\Q}, A)$ \\ Proposition \ref{prop:virtual-conjugacy-Denjoy}} ;
    \node (D) [box, below of=A, yshift=2.5cm] {Flow equivalence \\ $(\alpha,A) \bmod{\PGL_2(\Z)}$ \\ Theorem \ref{thm:flow-equivalence-Denjoy-systems} };
    \node (E) [box, below of=B,yshift=2.5cm] {Eventual flow equivalence \\ $\pm (\alpha, A) \bmod{\Z}$ \\ Theorem \ref{thm:eventual-flow-Denjoy}};
    \node (F) [box, below of=C, yshift=2.5cm] {Virtual flow equivalence \\ $(n\alpha, A)\bmod{\PGL_2(\Z)}$ \\ Remark~\ref{rem:virtual-flow-equivalence-Denjoy} };

    % Draw arrows
    \draw[->, double] (A) -- (B);
    \draw[->, double] (B) -- (C);
    \draw[->, double] (A) -- (D);
    \draw[->, double] (B) -- (E);
    \draw[->, double] (C) -- (F);
    % \draw[->, double] (D) -- (E);
    \draw[->, double] (E) -- (F);
    \draw[->, double] (B) -- (A);
    \draw[->, double] (E) -- (B);
    \draw[->, double] (E) -- (A);
    \draw[->, double] (A) -- (E);
\end{tikzpicture}
\end{center}

In contrast, Mike Boyle communicated to us that among shifts of finite type there exist pairs that are eventually flow equivalent but not eventually conjugate (hence not topologically conjugate), and pairs that are flow equivalent but not eventually flow equivalent (see Theorem~\ref{thm:boyleexamples}).

\begin{Question}[from eventual flow equivalence to flow equivalence or topological conjugacy]
    In which classes of systems does eventual flow equivalence imply topological conjugacy? In which classes of systems does eventual flow equivalence imply flow equivalence?
\end{Question}

It may be that eventual flow equivalence implies topological conjugacy among large classes of minimal systems, but the proof would require a deeper understanding of topological invariants under flow equivalence and powers, which would have to subsume the difficulties of the self-induced Sturmian case (see proof of Theorem~\ref{thm:eventual-PGL2Z} as well as Example~\ref{ex:p-eventually-PSL2Z-equivalent} and Question~\ref{quest:weak-eventual-PGLZ}).
Indeed, our characterisations of eventual flow equivalence for Sturmian and Denjoy systems rely on our Theorem~\ref{thm:eventual-PGL2Z}, a novel arithmetic result which we believe is of independent interest.
\begin{Theorem}[eventual $\PGL_2(\Z)$-equivalence]
    Consider $\alpha,\beta\in \R\setminus\Q$. If for all but finitely many prime integers $n\in \Primes$ the numbers $n\alpha,n\beta$ are $\PGL_2(\Z)$-equivalent, then $\alpha \equiv \pm\beta \bmod{\Z}$.
    % Consider irrationals $\alpha,\beta\in \R\setminus \Q$.
    % Consider a subset of rational primes $\Primes'\subset \Primes$ having full Dirichlet density (say all but finitely many).
    % If for all $n\in \Primes'$, the numbers $n\alpha,n\beta$ are $\PGL_2(\Z)$-equivalent, then $\alpha=\pm \beta \bmod{\Z}$.
\end{Theorem}

Let us mention the following Question~\ref{quest:weak-eventual-PGLZ}, hinting at the arithmetic subtleties of eventual $\PGL_2(\Z)$ equivalence.

\begin{Question}[$\N'$-virtual $\PGL_2(\Z)$-equivalence]
    For a subset of $\N' \subset \N$, say that $\alpha,\beta$ are $\N'$-virtually $\PGL_2(\Z)$-equivalent when for all $n\in \N'$ the numbers $n\alpha,n\beta$ are $\PGL_2(\Z)$-equivalent.

    For which subsets $\N'$ does $\N'$-virtual $\PGL_2(\Z)$-equivalence imply $\PGL_2(\Z)$-equivalence, or further equality up to sign and addition by an integer?
    
    Theorem~\ref{thm:eventual-PGL2Z} replies both positively when $\N'$ is a subset of primes of full Dirichlet density while Example~\ref{ex:p-eventually-PSL2Z-equivalent} shows that both are false when $\N'$ consists of the powers of a single prime.
    
    We believe that it is false if $\N'$ is a finite subset of primes (Example~\ref{ex:p-eventually-PSL2Z-equivalent} suggests a constructive way to prove this), but what if $\N'$ is a subset of primes of positive Dirichlet density?
\end{Question}

\subsection*{Acknowledgements}
We sincerely thank Pierre Arnoux, Mike Boyle, John Chaika, Vincent Delecroix, Pascal Hubert, Chris Leininger for their continuous interest in reading and effort in responding to our questions.
% informative help and discussions regarding this material.
We are also grateful to Menny Akka, Daniel Martin and François Maucourant for discussions concerning eventual $\PGL_2(\Z)$-equivalence.

% This material is based upon work supported by the National Science Foundation under Grant No. DMS-2247553.

\subsection*{Disclosure}
This work did not benefit from the use of 
artificial intelligence except for typographical corrections.

\section{Background}
\label{sec:systems}
A \emph{system} $(X,S)$ is a non-empty compact metric space $X$ with a homeomorphism $S \colon X\to X$.
A \emph{morphism} of systems is a continuous map $\pi \colon (X,S) \to (Y,T)$ such that $\pi S = T \pi$, and it is called a \emph{factor} when it is surjective. We say $(X,S)$ and $(Y,T)$ are semi-conjugate if there exists a factor $\pi \colon (X,S) \to (Y,T)$.

% For every system $(X,S)$, the space of Borel probability measures on $X$ that are $T$-invariant form a non-empty convex set.

A system $(X,S)$ is \emph{minimal} when its only closed $S$-invariant subsets are $\emptyset$ and $X$, or equivalently when every $S$-orbit is dense. A system $(X,S)$ is called a \emph{Cantor system} when $X$ is a Cantor space, and is called a \emph{zero-dimensional} system when $X$ is zero-dimensional.

For a system $(X,S)$ we let $\mathcal{M}(X,S)$ denote the space of finite Borel $S$-invariant measures on $X$, and $\mathcal{M}_{p}(X,S) \subset \mathcal{M}(X,S)$ the subspace of probability measures. 

A system $(X,S)$ is \emph{uniquely ergodic} when there is a unique Borel probability measure on $X$ that is $S$-invariant.

For a property $\mathtt{P}$ of systems, a system $(X,T)$ is \emph{totally $\mathtt{P}$} when the system $(X,T^{n})$ has property $\mathtt{P}$ for every $n \in \N$.

\subsection{Subshifts}
\label{subsec:subshifts}

Fix a finite set $\Al$ of cardinality $\Card{\Al} \in \N_{>1}$ called the alphabet and let $\Al^{\star}=\bigcup_{n=1}^{\infty} \Al^n$ be the set of finite words over $\mathcal{A}$.
The set of bi-infinite sequences $\Al^\Z$ with the product topology is a compact space.
A basis of neighbourhoods for the topology is given by the clopen sets $\operatorname{Cyl}_k(w)= \{x\in \Al^\Z \colon x_{[k,k+n-1]}=w\}$ for $k \in \Z$ and $w\in \Al^n$.
It is metrized by $$d(x,x') = 2^{-\inf\{|k|\in \N \colon x_{k} \ne x'_{k}\}}.$$
The shift $\Shift \colon \Al^{\Z} \to \Al^{\Z}$ defined by $\Shift(x)_{n} = x_{n+1}$ is a homeomorphism.

A subsystem $(X,\Shift)$ of $(\Al^{\Z}, \Shift)$ is called a \emph{subshift} on $\Card{\Al}$ symbols.
A word $w\in \Al^{\star}$ is a subword of $x\in X$ when there exists $k\in \Z$ such that $x\in \operatorname{Cyl}_k(w)$.
The \emph{language} of $(X,\Shift)$ is the set $\Al^{\star}_X \subset \Al^{\star}$ of all subwords of all elements.
There is an equivalence between subshifts of $(\Al, \Shift)$ and languages in $\Al^\star$ which are factorial (closed by taking subwords) and extendable (every word is a proper subword of some other word).
The factor complexity of $(X,\Shift)$ is the function $C_X\colon \N\to \N$ defined by $C_X(n)=\Card{\Al^{\star}_X \cap \Al^n}$.

The subshift $(X,\Shift)$ is called \emph{uniformly recurrent} if for every $x\in X$, every word in $\Al^\star_X$ occurs in the sequence $\Shift^k(x)$ infinitely often and with bounded gaps.
This is equivalent to the minimality of $(X,\Shift)$.

The subshift $(X,S)$ has a \emph{recurrence function} $R_X\colon \N\to \N$ which to $n\in \N$ associates the minimum $R(n)$ of the integers $r\in \N$ such that every word in $\Al^r_X$ contains every word in $\Al^n_X$.
One may show (see~\cite[Theorem 2.3]{Berthe-Delecroix_beyond-substitutive-S-adic_2014}) that $C_X\le R_X$.

A subshift $(X,S)$ is called \emph{linearly recurrent} when there exists $C\in \R_{>0}$ such that for all $n\in \N$ we have $R_X(n)\le Cn$.
Linear recurrence implies uniform recurrence (equivalently minimality), sublinear factor complexity (as $C_X\le R_X$), and unique ergodicity (this is \cite[Theorem 2.5]{Berthe-Delecroix_beyond-substitutive-S-adic_2014}, which refers to \cite{Durand-Host-Skau_Substitutional-systems-Bratteli-dimension-groups_1999}).

Let us recall a criterion due to Boshernitzan for linear recurrence of minimal subshifts (see \cite[Theorem 21]{Delecroix_translation-surfaces-Salta_2015} for a proof).
For a subshift $(X, \Shift)$ with an invariant probability measure $\mu\in \mathcal{M}_{p}(X,\Shift)$, define the function $\epsilon_\mu\colon \N\to \R$ by $\epsilon_\mu(n)=\min\{\mu(\operatorname{Cyl}_0(w)\cap X) \colon w \in \Al^n_X \}$ (the minimal measure of cylinders of width $n$).
The minimal subshift $(X,\Shift)$ is linearly recurrent if and only if any $\mu\in \mathcal{M}_{p}(X,\Shift)$ satisfies $\inf n\epsilon_\mu(n)>0$.
It follows that if a subshift is totally minimal and linearly recurrent, then it is totally linearly recurrent, hence totally uniquely ergodic.

Let us finally mention that the notions of uniform recurrence and linear recurrence can be characterised in terms of S-adic representations (see~\cite[Section 5]{Berthe-Delecroix_beyond-substitutive-S-adic_2014}).
In particular, it is shown in \cite{Durand_characterization-substitutive-sequences-return-words_1998} that a subshift is linearly recurrent if and only if it is a strongly primitive and proper S-adic subshift.

\subsection{Rotations and Sturmian systems}

On the circle $\Torus=\R\bmod{\Z}$, let $R_\alpha  \colon \Torus \to \Torus$ denote the rotation by $\alpha \in \Torus$.
If $\alpha$ is rational then $R_\alpha$ is periodic. 
If $\alpha$ is irrational then $(\Torus, R_\alpha)$ is minimal and uniquely ergodic, with Lebesgue measure on $\Torus$ the unique invariant probability measure.

For an irrational $\alpha\in \Torus$ we define the Sturmian system $(X_\alpha, \Shift_\alpha)$ as the simplest symbolic encoding of $R_\alpha$, that is a subshift on two symbols which factors onto $(\Torus, R_\alpha)$. 
For this, we construct a point $x \in \{0,1\}^{\Z}$ by
\begin{equation*}
x_{n} = 
\begin{cases}
        0 & \text{if } n \alpha \in [0,1-\alpha)\\
        1 & \text{if } n \alpha \in [1 - \alpha,0)
\end{cases}
\end{equation*}
and let $X_{\alpha}$ be the closure of the $\Shift$-orbit of $x$ in $\{0,1\}^{\Z}$. 

The language of $X_{\alpha}$ has complexity $C_{X_\alpha}(n) = n+1$ (which is the smallest possible among aperiodic sequences, by Morse-Hedlund~\cite{HedlundMorse1940}).
Conversely, every recurrent subshift of complexity $C_{X}(n) = n+1$ is obtained by the previous construction for a unique $\alpha\in \Torus$. 

It is well known \cite{Fogg_Subsitutions-in-dyn-arit-comb_2002} that for irrational $\alpha\in \Torus$ the system $(X_{\alpha},\Shift_{\alpha})$ is minimal and uniquely ergodic.

The irrational $\alpha\in \R\setminus \Q$ has a unique Euclidean continued fraction expansion:
\begin{equation*}
%\textstyle
\Ecf{a_0,a_1,\dots}=a_0+\frac{1}{a_1+\dots}
%\Ecf{a_0,a_1,a_2,\dots}=a_0+\frac{1}{a_1+\frac{1}{a_2+\dots}}
\quad \mathrm{with} 
\quad a_j \in \N
\quad \mathrm{and} 
\quad \forall j>0,\; a_j>0.
\end{equation*}
Thus it is the fixed point of the infinite product of matrices $R^{a_0}L^{a_1}\dots$ where 
\begin{equation*}
    L =
    \begin{psmallmatrix}
    1 & 0\\
    1 & 1
    \end{psmallmatrix},
    \qquad
    R =
    \begin{psmallmatrix}
    1 & 1\\
    0 & 1
    \end{psmallmatrix}
\end{equation*}
act by linear fractional transformations on $\R\Proj^1$.

Such a sequence can be generated from any seed by applying its S-adic representation, which is determined by the continued fraction expansion of $\alpha$ (see~\cite{Schmieding-Simon_hexp_2024}).

\subsection{Homeomorphisms of the circle and Denjoy systems}
\label{subsec:denjoy-systems}

Let us briefly recall, following \cite{PSS_C-star-Denjoy_1986}, the main results of Poincaré~\cite{Poincare_courbes-equa-diff_1885} concerning homeomorphisms of the circle and their rotation numbers, and of Markley~\cite{Markley_homeo-circle_1970} describing the invariants for conjugacy classes in $\Homeo^+(\Torus)$ of Denjoy homeomorphisms.

%A \emph{semi-conjugation} from the system $(X,S)$ to the system $(Y,T)$ is a surjective continuous function $H\colon X\to Y$ such that $H\circ S = T\circ H$.

Let $F\in \Homeo^+(\Torus)$ be an orientation preserving homeomorphism of the circle.
Choose a lift to the universal cover $\tilde{F}\colon \R\to \R$: it is determined by $\lfloor \Tilde{F}(0)\rfloor \in \Z$ as any two lifts differ by an integer.
There exists $\tau(\tilde{F}) \in \R$ such that for all $x\in \R$ we have $\tfrac{1}{n}(\title{F}^n(x)-x)\to \tau(\tilde{F})$, and the value $\tau(\tilde{F}) \in \R\bmod{\Z}$ does not depend on the lift: it is called the \emph{rotation number} of $F$ and denoted $\rho(F)\in \Torus$.

The rotation number $\rho \colon \Homeo^+(\Torus)\to \Torus$ is invariant by conjugacy in $\Homeo^+(\Torus)$.
More generally, if $F_1, F_2\in \Homeo^+(\Torus)$ and there exists a continuous surjective map $H\colon \Torus\to \Torus$ satisfying $HF_{1} = F_{2}H$, then $\rho(F_1)=\pm \rho(F_2)$ where the sign depends on the orientation preserving or reversing behaviour of $H$.
Moreover, for all $n\in \Z$ we have $\rho(F^n)=n\rho(F)$.

If $\rho(F)$ is rational then $F$ has a periodic point, and $F$ is conjugate to $R_{\rho(F)}$.
If $\rho(F)$ is irrational, then every point $x\in \Torus$ has an orbit $\{F^n(x)\colon n\in \Z\}$ with the same cyclic order structure as $\{n \cdot \rho(F)\colon n\in \Z\}$, hence $F$ is semi-conjugate to $R_{\rho(F)}$ by a map $H\colon \Torus \to \Torus$ which is unique up to post-composing by a rotation.

% Consequently, when $\rho(F)$ is irrational, $f$ is conjugate to $R_{\rho(F)}$ if and only if the semi-conjugacy $h$ is a homeomorphism.
We say that $F$ is a \emph{Denjoy homeomorphism} when it is semi-conjugate but not conjugate to $R_{\rho(F)}$, that is when $\rho(F)$ is irrational and the semi-conjugacy $H\colon \Torus \to \Torus$ is not a homeomorphism (in which case it is necessarily orientation preserving).

Fix a Denjoy homeomorphism $F\in \Homeo^+(\Torus)$ and let $H\colon \Torus \to \Torus$ be a semi-conjugacy from $F$ to $R_\alpha$.
The homeomorphism $F$ is uniquely ergodic with invariant probability measure $dH$. 
% (assigning to the interval $[x_1,x_2]\subset \Torus$ the measure $H(x_1)-H(x_2)$).
The support of $dH$, denoted $\Denjoy_F\subset \Torus$, is homeomorphic to a Cantor space, and the \emph{Denjoy system} $(\Denjoy_F, F)$ is the unique minimal $F$-invariant subsystem of $(\Torus, F)$.
The restriction of $H$ to $ \Denjoy_{F}$ yields a factor $(\Denjoy_{F}, F) \to (\Torus, R_{\rho(F)})$.

The open set $\Torus\setminus \Denjoy_F$ is a disjoint union of open intervals $I_n\subset \Torus$ indexed by $\N$, and $H$ collapses the closure of each interval $I_n$ to a point.
The map $F$ uniquely determines that set of points $\{H(I_n)\colon n\in \N\}$ up to post-composition by a rotation: denote by $Q(F)$ its rotation class.
The set $\{I_n\colon n\in \N\}$ with its cyclic order and $F$-orbit structures is isomorphic to the set $Q(F)$ with its cyclic order and $R_{\rho(F)}$-orbit structures.

The Denjoy homeomorphism $F$ is characterised up to topological conjugacy in $\Homeo^+(\Torus)$ by $\left(\rho(F), Q(F)\right)$ and in $\Homeo(\Torus)$ by $\pm \left(\rho(F), Q(F)\right)$~\cite{Markley_homeo-circle_1970}. Similarly, the Denjoy system $(\Denjoy_F, F)$ is characterised up to cyclic-order preserving topological conjugacy by $\left( \rho(F), Q(F) \right)$ and up to topological conjugacy by $\pm \left(\rho(F), Q(F)\right)$ (see~\cite[Remark 3]{PSS_C-star-Denjoy_1986}).
From now on, for a pair consisting of an irrational $\alpha\in \Torus$ and a countable $R_\alpha$-invariant subset $A \subset \Torus$ containing the origin, we denote such a Denjoy homeomorphism by $(\Torus, F_{\alpha, A})$ and its minimal subsystem by $(\Denjoy_{\alpha, A}, F_{\alpha, A})$.

The following is a key fact for us. The first part is folklore (for a proof see for instance~\cite{Masui_Denjoy_2009}).
\begin{Proposition}[from Sturmian to Denjoy]
Let $\alpha\in \Torus$ be irrational and $A=\{R_\alpha^n(0)\colon n\in \Z\}$.
\begin{enumerate}[noitemsep]
\item
The Denjoy system $(\Denjoy_{\alpha, A}, F_{\alpha, A})$ is topologically conjugate to the Sturmian system $(X_\alpha,\Shift_\alpha)$.
\item
For all $m\in \Z \setminus \{0\}$, the Denjoy systems $(\Denjoy_{\alpha, A}, F_{\alpha, A}^m)$ and $(\Denjoy_{m\alpha, A}, F_{m\alpha, A})$ coincide, and are topologically conjugate to the system $(X_\alpha,\Shift^m_\alpha)$, which is not Sturmian when $\lvert m \rvert >1$.
\end{enumerate}
\end{Proposition}

Note that a Denjoy system factors onto a unique Sturmian system, given by its rotation number $\alpha \bmod{\Z}$.

\subsection{Interval exchange transformations and their subshifts}
\label{subsec:interval-exchange-systems}

\begin{comment}
Let $I\subset \R$ be a bounded open interval. An interval exchange transformation (IET) on $I$ is a bijective map $T\colon D_{T} \to D_{T^{-1}}$ such that $D_{T}, D_{T^{-1}}$ are subsets of $I$ whose complements are finite sets (with the same cardinality) and the restriction of $T$ to each component of $D_T$ is a
translation onto some component of $D_{T^{-1}}$.

The connected components of $D_T$ are intervals denoted $I_k$ which we label by $k\in \Al=\{1,\dots, d\}$ in increasing order.
The \emph{discontinuities} of $T$ are the points $\alpha_k\in \R$ such that $I_k=(\alpha_{k-1}, \alpha_k)$ and the interval lengths $\lambda_k=\alpha_k-\alpha_{k-1}$ define the \emph{length vector} $\lambda \in \R_{>0}^{\Al}$ of $T$. 

An interval exchange transformation (IET) $T$ on $d$ intervals is an orientation preserving right-semi-continuous piecewize-isometry of an interval $[\alpha_0,\alpha_d) \subset \R$ with $d-1$ of discontinuities. Hence there exists a partition of $[\alpha_0,\alpha_d)$ into $d$ intervals $I_k= [\alpha_{k-1}, \alpha_{k})$ indexed by $\Al$ which are permuted by $T$ according to a permutation $\sigma \colon \Al \to \Al$ and in restriction to each of them it $T$ is a translation, given by the formula:
\begin{equation*}
    T(x\in I_k) = (x-\alpha_{k-1}) +  \sum_{l>k, \sigma(l)<\sigma(k)} \lambda_l \: - \sum_{l<k, \sigma(l)>\sigma(k)} \lambda_l.
\end{equation*}
\end{comment}

Fix our alphabet $\Al=\{1,\dots, d\}$ for $d\in \N_{\ge 1}$.
Consider a permutation $\sigma \colon \Al \to \Al$, and length parameters $\lambda\in (\R_{>0})^\Al$. 
We fix $\alpha_0=0$, define the \emph{discontinuities} $\alpha_k=\sum_{1}^k \lambda_j$ for $k\in \Al$, and consider the intervals $I_k=(\alpha_{k-1},\alpha_{k}]$ of lengths $\lambda_k=\alpha_k-\alpha_{k-1}$.

The interval exchange transformation (IET) with permutation $\sigma$ and lengths $\lambda$ is the function $T\colon (\alpha_0,\alpha_d] \to (\alpha_0,\alpha_d]$ permuting the intervals $I_k= (\alpha_{k-1}, \alpha_{k}]$ according to $\sigma$. It is given by the formula
\begin{equation*}
    T(x\in I_k) = (x-\alpha_{k-1}) +  \sum_{l>k, \sigma(l)<\sigma(k)} \lambda_l \: - \sum_{l<k, \sigma(l)>\sigma(k)} \lambda_l.
\end{equation*}
Two IETs are called \emph{equivalent} when they have the same permutation and their length parameters yield proportional vectors in $\R_{>0}^\Al$.
Every equivalence class contains a unique representative with $\sum \lambda_k=1$, defined on the interval $(0,1]$.
% , determined by its permutation $\pi \colon \Al \to \Al$ and length vector $\lambda \in \R_{>0}^\Al$ with $\sum \lambda_k = 1$.

The \emph{symbolic encoding} of $T$ is the subshift $(X,\Shift) \subset (\Al^\Z, \Shift)$ defined as follows.
We define $X_o \subset (\alpha_0,\alpha_d]$ to be the complement of the orbits of all discontinuities: a point $x\in X_o$ yields a sequence of labels $l_j \in \Al$ recording the intervals $I_{l_j}$ containing the $T^j(x)$.
This defines a map $X_o \to \Al^\Z$ and the closure of its image is $X$, a subshift. 
Such a subshift $(X,\Shift)$ will be called an \emph{interval exchange system} (IES). 

For $m\in \Z$, the $m$-th power $T^m$ is an interval exchange with discontinuities given by $\alpha^{(m)}_k\in \bigcup_{n=1}^{m} T^{1-n}\left(\{\alpha_0,\dots,\alpha_d\}\right)$, of which there are at most $(d-1)\lvert m\rvert+1$, with equality unless some discontinuities belong to the same orbit.
The intervals $(\alpha^{(m)}_{k}, \alpha^{(m)}_{k+1}] \subset (\alpha_0,\alpha_d]$ whose endpoints are consecutive discontinuities for $T^m$ correspond to cylinders $\operatorname{Cyl}_0(w)\subset X$ of words $w\in \Al_X^m$.

We say that $T$ satisfies \emph{Keane's infinite distinct orbit condition} (or has ido) when the past-orbits of the discontinuities $\{T^{-n}(\alpha_k)\colon n\in \N\}$ are all infinite and disjoint.
This implies that the map $X_o \to X$ is injective and its inverse, which is partially defined on its image, extends to a unique map $X\to [\alpha_0,\alpha_d]$.
%Note that the image of the cylinder set $\operatorname{Cyl}_0(k)\subset X$ is $I_k\subset (\alpha_0,\alpha_d]$.

If $T$ has ido, then $(X,\Shift)$ has complexity function $C_X(n)=(d-1)n+1$ and it is minimal (Keane showed \cite{Keane_IET_1975} that the ido condition implies that every $T$-orbit is dense).
Conversely, if a minimal IES $(X,\Shift)$ has $C_X(n)=(d-1)n+1$, then every IET with that encoding has ido (so it makes sense to say whether a minimal IES has ido).

If a minimal subshift on $d$ unordered letters of complexity $n\mapsto (d-1)n+1$ comes from an IES, then the order on its alphabet $\Al$ and the permutation $\pi\colon \Al \to \Al$ are uniquely determined up to the action of either the reflection group $\Z/2$ generated by the order reversing involution $k\mapsto d+1-k$ or else the dihedral group $\Z/d\rtimes \Z/2$ obtained as its semi-direct product by the cyclic group $\Z/d$ generated by the cyclic order preserving permutation $k\mapsto k+1\bmod{d}$.
% (These are the only two subgroups of $\mathfrak{S}(\Al)$ preserving the contiguity relation defined by $\lvert i-j\rvert=1\bmod{d}$ which contain the reflection group $\Z/2$.)

Given a minimal subshift $(X,\Shift)\subset(\Al,\Shift)$ of complexity $C_X(n)=(d-1)n+1$ on the ordered alphabet $\Al$, one may decide from its language whether there exist a permutation $\pi\colon \Al\to \Al$ and lengths $\lambda\colon \Al \to \R_{>0}$ yielding an IET with ido having that encoding (\cite{Ferenczi-Zamboni_Languages-k-IET_2008}).
When so, the permutation $\pi\colon \Al \to \Al$ is uniquely determined, and it follows from \cite[Section 8.1]{Yoccoz_IET-translation-surfaces_2010} that the possible length parameters $\lambda$ are in bijective correspondence with the invariant measures of $(X, \Shift)$.

The permutation $\pi$ is called \emph{irreducible} when $\pi^{-1}(\{1,\dots, k\})=\{1,\dots, k\} \iff k\in \{1,d\}$; this means that $T$ does not restrict to an IET on a smaller interval.
If $T$ has ido then $\pi$ is irreducible.
If $\pi$ is irreducible and the rational relations between the lengths $\lambda_k$ are generated by $\sum \lambda_k=1$ then $T$ has ido.
Moreover, Masur \cite{Masur_IET-MF_1982} and Veech \cite{Veech_Gauss-measure-induction-IET-ae-UE_1982} showed that for every irreducible $\pi$, almost all probability vectors $\lambda$ yield an IES which is uniquely ergodic.

Note that a Denjoy system $(\Denjoy_{\alpha, A}, F_{\alpha, A})$ with $A\bmod{R_\alpha}$ of cardinal $d\in \N_{\ge 2}$ is topologically conjugate to an IES on $d$ intervals and its permutation $\pi$ is a $d$-cycle, which preserves the cyclic order on $\Al\bmod{d}$ in the sense that if $(i,j,k)$ are cyclically ordered then $(\pi(i),\pi(j),\pi(k))$ is also cyclically ordered.
Conversely, all such cyclic-order-preserving $d$-cycles arise from Denjoy systems.
For instance, the Sturmian system $(X_\alpha, \Shift_\alpha)$ is the IES associated to the transposition $\tau$ of $d=2$ intervals with lengths $\lambda_1 = \alpha$ and $\lambda_2 = 1-\alpha$.

\section{Coinvariants, states and the suspension}
\label{sec:coinvariants}
\subsection{Coinvariants, states and cohomology of the suspension}

Let $(X,S)$ be a system and denote by $C(X,\Z)$ the abelian group of continuous integer valued functions on $X$.
It is endowed with the coboundary homomorphism $f\mapsto f - f\circ S$, whose image defines the subgroup of coboundaries $\textrm{cbd}_{S}(X,\Z)$ and whose cokernel defines the group of \emph{coinvariants}
\[\CoInv_{S} = C(X,\Z) / \textrm{cbd}_{S}(X,\Z).\]

A homomorphism $\tau \colon \CoInv_{S} \to \R$ is called \emph{positive} if $\tau([f])$ is nonnegative whenever $[f]=[g]$ for some nonnegative $g \in C(X,\Z_{\ge 0})$, and a positive homomorphism $\tau$ is called a \emph{state} when $\tau([1]) = 1$. We let $\mathcal{P}(X,S)$ denote the space of states for a system $(X,S)$.

% The following is shown for the case of minimal cantor systems in % Herman, R. H., Putnam, I. F., and Skau, C. F., Ordered Bratteli diagrams, dimension groups and topological dynamics, Internat. J. Math. 3. 
The following is very well known; a proof can be found in~\cite[Section 1.6]{Boyle-Handelman_orbit-flow-oredered-cohomology_1996}.

\begin{Proposition}[states as measures]
\label{prop:states-as-measures}
Let $(X,S)$ be a zero-dimensional system. The space of states is in bijection with the space $\mathcal{M}_{p}(X,S)$ of $S$-invariant Borel probability measures on $X$ by the map sending ${\mu\in \mathcal{M}_{p}(X,S)}$ to the state $\tau_\mu \in \mathcal{P}(X,S)$ defined by 
$$\tau_{\mu}(f) = \int_{X}f \, d\mu.$$

More precisely, the space of states and the space of $S$-invariant Borel probability measures on $X$ both have the structure of a non-empty metrizable Choquet simplex. 
% in a Fréchet space (a locally convex topological vector space that is metrizable and complete).
The map $\mu \mapsto \tau_\mu$ is an isomorphism, in particular it preserves the convex structure.    
\end{Proposition}

\begin{comment}
\begin{proof}
% One may think of $C(X,\Z)$ as a lattice in the Banach space $C(X,\R)$, and the space of continuous real linear forms on $C(X,\Z)$ is the topological dual of the Banach space $C(X,\R)$.
The space of continuous real linear forms on $C(X,\Z)$ is the topological dual of the Banach space $C(X,\R)$.
% A continuous linear form $C(X,\R) \to \R$ is called positive when it takes non-negative values on non-negative functions $f \in C(X,\R_{\ge 0})$.
By the Riesz–Markov–Kakutani representation theorem, since $X$ is compact, the map sending a Borel probability measure $\mu\in \mathcal{P}(X)$ to the continuous linear form $f \mapsto \int_X f \, d\mu$ on $C(X,\R)$ extends to an isomorphism between the space of signed measures on $X$ and the topological dual of $C(X,\R)$.

The space of continuous linear forms on the quotient $\CoInv_{S}$ correspond to the measures that vanish on the subspace of coboundaries $\mathrm{cbd}(X,S) \subset C(X,\Z)$, that is the space of signed $S$-invariant measures on $X$.
The positive measures and the probability measures correspond to the positive linear forms and the states.
% More precisely, the space of states and the space of invariant probability measures on $(X,S)$ both have the structure of a standard simplex in a Fréchet space (a locally convex topological vector space that is metrizable and complete).
% The map $\mu \mapsto \tau_\mu$ is an affine isomorphism, in particular it preserves the convex structure.
\end{proof}
\end{comment}

A morphism of systems $\pi\colon (X,S)\to (Y,T)$ yields an induced homomorphism of groups $\pi^*\colon \CoInv_{T} \to \CoInv_{S}$ given by $\pi^{*}([f]) = f \circ \pi$ and an associated affine morphism of simplices $\pi_*\colon \PP(X,S) \to \PP(Y,T)$ given by $\pi_{*}(\tau)([f]) = \tau(\pi^{*}([f]))$.
For minimal Cantor systems, the induced map $\pi^*$ on coinvariants is injective, and hence the map $\pi_*$ on states is surjective (see~\cite[Prop. 3.1]{Glasner-Weiss-1995}).

We define the \emph{subgroup of infinitesimals} $\Inf_{S}\subset \CoInv_{S}$ as the set of classes $[f] \in \CoInv_{S}$ such that $\tau([f]) = 0$ for every state $\tau$ on $\CoInv_{S}$, or equivalently such that $\int_{X} f\, d\mu = 0$ for every $S$-invariant Borel probability measure $\mu$ on $X$; in formula, 
$$\Inf_{S} = \bigcap_{\tau \in \mathcal{P}(X,S)} \ker\left(\tau \colon \CoInv_{S} \to \R\right).$$

Let $(X,S)$ be a system.
The product space $X \times \R$ is endowed with the diagonal action of the lattice $\Z\subset \R$ by $n\cdot(x,s)=(S^{n}(x),s-n)$.
The quotient $\Sigma_S X := X\times \R \bmod{\Z}$ is endowed with the $\R$-action by left translation of the right component, namely the flow $\Phi_{s^{\prime}}(x,s)=(x,s+s^{\prime})$.
We call $\Sigma_S X$ the \emph{suspension} of the system $(X,S)$.

Note that if $(M,S)$ is a minimal subsystem of $(X,S)$ then $\Sigma_S M$ is a minimal subsystem of $\Sigma_S X$. Conversely, if $N$ is a minimal subsystem of $\Sigma_S X$ then the set $M=\{x\in X\colon (x,0)\in N\}$ yields a minimal subsystem $(M,S)$ of $(X,S)$.
Moreover, the projection $\Sigma_S X \to X$ yields an isomorphism $\mathcal{P}(\Sigma_S X,\Phi) \to \mathcal{P}(X,S)$ between the respective spaces of invariant Borel probability measures (see~\cite[Lemma 3.8]{ItzaOrtiz_Denjoy-flow_2015}).

For a topological space $Y$, let $\pi^{1}(Y)$ denote the group of homotopy classes of maps from $Y$ to the unit circle $\Sph^{1}$, for which the group operation is pointwise multiplication in the target $\Sph^1$. 
If $Y$ is compact Hausdorff then $\pi^{1}(Y)$ is isomorphic to the first integral \v{C}ech cohomology group $\check{H}^{1}(Y,\Z)$.

\sloppy Now suppose that $(X,S)$ is a zero-dimensional system. 
Given a map ${\varphi \colon \Sigma_{S}X \to \Sph^{1}}$, there exists for every $z \in \Sigma_{S}X$ a continuous function $\phi_{z} \colon \R \to \R$ satisfying $\phi_z(0)=0$ and $\varphi(z+s) = \varphi(z)e^{2 \pi i \phi_{z}(s)}$ for all $s \in \R$. 
Denote by $C_{+}(\Sigma_{S}X,\Sph^{1})$ the set of maps $\varphi \colon \Sigma_{S}X \to \Sph^{1}$ such that we may choose $\phi_{z}$ to be non-decreasing for all $z \in \Sigma_{S}X$. 
The abelian group $\pi^{1}(\Sigma_{S}X)$ is endowed with the pre-order defined by the positive cone
\begin{equation*}
\pi^{1}_{+}(\Sigma_{S}X) = \{ [\varphi] \colon \varphi \in C_{+}(\Sigma_{S}X,\Sph^{1})\} \subset \pi^{1}(\Sigma_{S}X).
\end{equation*}
For $f \in C(X,\Z)$, we define $\varphi_{f} \in C(\Sigma_{S}X,\Sph^{1})$ by
\begin{equation*}
\varphi_{f}(x,s) = e^{2 \pi i s f(x)}, \quad 0 \le s < 1.
\end{equation*}
The following result is well known, and a proof may be found in  \cite[Proposition 4.5]{Boyle-Handelman_orbit-flow-oredered-cohomology_1996}.

\begin{Proposition}[from coinvariants to homotopy classes]
\label{prop:exponentialiso}
Let $(X,S)$ be a zero-dimensional system. 
The exponential map yields an isomorphism of groups:
\begin{equation}
\begin{gathered}
\CoInv_{S} \to \pi^{1}(\Sigma_{S}X)\\
[f] \longmapsto [\varphi_{f}]
\end{gathered}
\end{equation}
Moreover if $(X,S)$ is minimal, then the isomorphism respects the order structures, and hence gives an isomorphism of ordered abelian groups
\begin{equation*}
(\CoInv_{S},\CoInv_{S}^{+}) \cong (\pi^{1}(\Sigma_{S}X),\pi^{1}_{+}(\Sigma_{S}X)).
\end{equation*}
\end{Proposition}

\subsection{Coinvariants of Denjoy and interval exchange systems}

Recall from Subsection~\ref{subsec:denjoy-systems} the definition of Denjoy homeomorphisms as well as the invariants classifying them up to (orientation preserving) conjugacy in $\Homeo(\Torus)$.

Fix an irrational $\alpha\in \Torus$ and an $R_\alpha$-invariant countable subset $A\subset \Torus$ containing the origin.
Let $N\in \N\sqcup \{\infty\}$ denote the number of $R_{\alpha}$-orbits in $A$ and choose a point $a_i\in A$ in each of these orbits in such a way that $a_1=\alpha$.
Such pairs $(\alpha, A)$ parametrize Denjoy homeomorphisms $(\Torus, F_{\alpha, A})$ and Denjoy systems $(\Denjoy_{\alpha, A}, F_{\alpha, A})$ up to cyclic order preserving conjugacy.

The surjective map $H\colon \Torus \to \Torus$ yielding the factor $(\Torus, F_{\alpha, A}) \to (\Torus, R_{\alpha})$ restricts to a factor $\pi_{\alpha, A} \colon (\Denjoy_{\alpha, A}, F_{\alpha, A}) \to (\Torus, R_{\alpha})$.
Both systems $(\Torus, F_{\alpha, A})$ and $(\Denjoy_{\alpha, A}, F_{\alpha, A})$ are uniquely ergodic, with invariant probability measures given by $dH$ on $\Torus$ and its restriction to its support $\Denjoy_{\alpha, A}$.

The following is \cite[Theorem 5.3]{PSS_C-star-Denjoy_1986}.

\begin{Theorem}[Denjoy systems]
\label{thm:Denjoy-coinvariants-state}
Consider a Denjoy system $(\Denjoy_{\alpha, A}, F_{\alpha, A})$, and choose as above some representatives $a_i\in A$ of $A\bmod{R_\alpha}$ with $a_1=\alpha$. Then $\CoInv_{F_{\alpha, A}} \cong \bigoplus_0^N \Z$ and its unique state $\tau \colon \CoInv_{F_{\alpha, A}} \to \R$ has image $\Z + \sum_{1}^{N} \Z \cdot a_i$. The group of infinitesimals $\Inf_{F_{\alpha, A}}$ corresponds to the $\Z$-module of relations among elements in $A\subset \Torus$. In particular, for a Sturmian system $(X_{\alpha},\Shift_{\alpha})$, we have $\CoInv_{\Shift_{\alpha}} \cong \Z ^{2}$ and the image under the unique state is $\Z + \alpha \Z$.
\end{Theorem}

Let us briefly outline some of the ideas in the proof, as it will be useful later.

\begin{proof}[Ideas of the proof]
We identify points in $\Torus = \R\bmod{\Z}$ by their representatives in $[0,1)\subset \R$. Under $\pi_{\alpha, A}$, each point $a \in A \subset \Torus$ has two preimages $a^\pm \in \Denjoy_{\alpha, A}$ which are the endpoints of an interval $(a^-, a^+) \subset \Torus\setminus \Denjoy_{\alpha, A}$ collapsed to $a$ by $H\colon \Torus \to \Torus$.
The characteristic functions of the closed intervals $[0,a], [a,1]\subset \Torus$ are discontinuous, but their lifts $\chi_{a}^-, \chi_a^+ \colon \Denjoy_{\alpha, A}\to \{0,1\}$ belong to $C(\Denjoy_{\alpha, A}, \Z)$, and they sum to $\chi_a^-+\chi_a^+=1$.
Since the measure $dH$ is supported on $\Denjoy_{\alpha, A}$ and trivial on all intervals of the form $(a^-, a^+)\subset \Torus \setminus \Denjoy_{\alpha, A}$, the trace satisfies $\tau(\chi_{a}^-)=a\in [0, 1)$ and $\tau(\chi_{a}^+)=1-a\in [0, 1)$. The core of the proof then consists in showing that the quotient of $\CoInv_{\alpha, A}$ by coboundaries is generated by the (classes of) the constant function $1$ and of the $\chi_{a_i}^-=1-\chi_{a_i}^+$ for our representatives $a_i$ of $A\bmod{R_\alpha}$, and showing that these classes are all linearly independent.

Observe that the group of infinitesimals $\Inf_{F_{\alpha, A}}$ corresponds to the $\Z$-module of relations of the $\Z$-submodule of $\R$ generated by $\tau(1)=1$ and the $\tau(\chi_{a_i}^-) = a_i \in (0, 1)$, which amounts to the stated description.
\end{proof}

We now turn to the symbolic codings of IETs satisfying Keane's ido condition.
Their coinvariants and trace follow from \cite{Putnam_Cstar-algebras-minimal-homeo-Cantor_1989} or \cite{Putnam_Cstar-algebras-IET_1992}.

\begin{Theorem}[interval exchange systems]
    Consider an IET associated to a permutation of $\Al=\{1,\dots ,d\}$ and discontinuities $0=\alpha_0<\dots < \alpha_d = 1$, which satisfies Keane's ido condition.
    Let $(X,\Shift)\subset (\Al^\Z, \Shift)$ be the associated IES, which is minimal. Denote by $\chi_k \colon X\to \{0,1\}$ the characteristic function of the cylinder $\operatorname{Cyl}_0(k)\subset X$, corresponding to that of the interval $(\alpha_{k-1}, \alpha_{k}]\subset (0,1]$, of length $\lambda_k$. The group $\CoInv(X,\Shift)$ is isomorphic to a free abelian group $\Z^d$, and the classes of the characteristic functions $\chi_k$ form a basis.
    The trace map $\tau \colon \CoInv(X,\Shift) \to \R$ coming from the Lebesgue measure on $[0,1)$ sends $\chi_k$ to $\lambda_k$.
\end{Theorem}

\section{Flow equivalence}

\subsection{Flow equivalence and states}

A \emph{flow equivalence} between two systems $(X,S)$ and $(Y,T)$ is a homeomorphism between their suspensions $\Sigma_{S}X \to \Sigma_{T}Y$ which maps flow orbits bijectively to flow orbits, and preserves the orientation on each flow orbit. 
Two systems are called \emph{flow equivalent} when there exists a flow equivalence between their suspensions.

Note that a topological conjugacy $\psi \colon (X,S) \to (Y,T)$ induces a flow equivalence ${\Psi \colon \Sigma_{S}X \to \Sigma_{T}Y}$ defined by $\Psi \colon (x,t) \mapsto (\psi(x),t)$.

Now consider zero-dimensional systems $(X,S)$ and $(Y,T)$ with a flow equivalence ${\Psi \colon \Sigma_{S}X \to \Sigma_{T}Y}$. 
The map $\Psi$ induces an isomorphism of groups
\begin{equation*}
\begin{gathered}
\Psi^{*} \colon \pi^{1}(\Sigma_{T}Y) \to \pi^{1}(\Sigma_{S}X)\\
\Psi^{*}([\varphi]) = [\varphi \circ \Psi]
\end{gathered}
\end{equation*}
which preserves the positive cones, since $\Psi$ preserves the orientation induced by the flow. 
Hence it follows from Proposition~\ref{prop:exponentialiso} that the flow equivalence $\Psi \colon \Sigma_{S}X \to \Sigma_{T}Y$ also induces an isomorphism $\Psi^{*} \colon \CoInv_{T} \to \CoInv_{S}$. When $(X,S)$ is minimal, the standard isomorphism given in Proposition~\ref{prop:exponentialiso} respects the order (see~\cite[Prop. 4.5]{Boyle-Handelman_orbit-flow-oredered-cohomology_1996}), thus $\Psi^{*} \colon \CoInv_{T} \to \CoInv_{S}$ is an isomorphism of ordered groups. 
An explicit method to compute $\Psi^{*}$ without appealing to Proposition~\ref{prop:exponentialiso} can be found in~\cite[Sec. 4]{Schmieding-Yang_Map-subshift_2021}. 
Note that in general, $\Psi^{*}([1])$ need not be equal to $[1]$.

When $(X,S)$ and $(Y,T)$ are minimal, the induced isomorphism $\Psi^{*} \colon \CoInv_{T} \to \CoInv_{S}$ which is order preserving thus induces a linear isomorphism between the spaces of positive homomorphisms
\begin{equation}
\begin{gathered}
\Psi_{*} \colon \mathcal{P}(\CoInv_{S}) \to \mathcal{P}(\CoInv_{T})\\
\Psi_{*}(\tau)(g) = \tau(\Psi^{*}(g)).
\end{gathered}
\end{equation}

We summarize this in the following proposition.

\begin{Proposition}[flow equivalence yields isomorphism of coinvariants]
\label{prop:flow-coinvariants}
Suppose that $(X,S)$ and $(Y,T)$ are minimal Cantor systems and ${\Psi \colon \Sigma_{S}X \to \Sigma_{T}Y}$ is a flow equivalence. 
The map $\Psi$ induces a linear isomorphism between the spaces of invariant measures $\Psi_{*} \colon \mathcal{M}(X,S) \to \mathcal{M}(Y,T)$. In particular, for every $\mu_X \in \mathcal{M}_p(X,S)$ there exist unique $\mu_{Y} \in \mathcal{M}_{p}(Y,T)$ and $\lambda\in \R_{>0}$ such that $\tau_{\mu_{Y}}=\lambda \cdot \Psi_*(\tau_{\mu_{X}})$, and hence $\im(\tau_{\mu_{X}}) = \lambda \cdot \im(\tau_{\mu_{Y}})$.
\end{Proposition}

%The first statement follows from the fact that flow equivalence is equivalent to topological conjugacy of induced systems, together with the fact that if $(C,S_{C})$ is an induced system of $(X,S)$, then the map $\mu \mapsto \nu$ where $\nu(A) = \mu(A)/\mu(X)$ is a linear bijection between $\mathcal{M}(S)$ and $\mathcal{M}(S_{C})$; see~\cite[Prop. 2]{Durand-Ormes-Petite_self-induced-systems_2018} for details.

\begin{Remark}[average stretch factor]
Suppose that $(X,S)$ and $(Y,T)$ are minimal and uniquely ergodic with unique invariant probability measures $\mu_{X}$ and $\mu_{Y}$ respectively. Then both of $\PP(X,S)$ and $\PP(Y,T)$ are rays, given by scalar multiples of $\tau_{\mu_{X}}$ and $\tau_{\mu_{Y}}$. Hence if $\Psi \colon \Sigma_{S}X \to \Sigma_{T}Y$ is a flow equivalence then there exists a unique $\lambda\in\R_{>0}$ such that $\Psi_{*}(\tau_{\mu_{X}}) = \lambda \cdot \tau_{\mu_{Y}}$, so for all $[f] \in \CoInv_{T}$ we have $\tau_{\mu_{X}}(\Psi^{*}(f)) = \lambda \cdot \tau_{\mu_{Y}}(f)$.

The number $\lambda$ can be interpreted as the average expansion of the flow equivalence on leaves, in the following sense. Given the flow equivalence $\Psi \colon \Sigma_{S}X \to \Sigma_{T}Y$,
the function
$\alpha_{\Psi} \colon \Sigma_{f}X \times \R \to \R$
defined implicitly by
$\Psi(z+t) = \Psi(z) + \alpha_{\Psi}(z,t)$
is well-defined, continuous, and satisfies the cocycle condition
%\begin{equation*}
$\alpha_{\Psi}(z,t_{1} + t_{2}) = \alpha_{\Psi}(z,t_{1}) + \alpha_{\Psi}(z+t_{1},t_{2})$.
%\end{equation*}
% The cocycle $\alpha_{f}$ is always continuous (a proof can be found in Proposition 6.2 of \url{https://arxiv.org/pdf/1506.02694.pdf}).
Since $(X,S)$ and $(Y,T)$ are uniquely ergodic, so are the corresponding flows $(\Sigma_{S}X,\R)$ and $(\Sigma_{T}Y,\R)$.
Using ergodicity of the flow and the ergodic theorem, there exists $\lambda_{\Psi} \in \R$ such that for all $z \in \Sigma_{S}X$, we have
$\lim_{t\to \infty}\frac{1}{t}\alpha_{\Psi}(z,t)=\lambda_{\Psi}$.
This $\lambda_{\Psi}$ is precisely the $\lambda$ given by Proposition~\ref{prop:flow-coinvariants} 
(see~\cite{Schmieding-Yang_Map-subshift_2021} for details).
\end{Remark}

\subsection{Flow equivalence and discrete cross sections}

Let $(X,S)$ be a zero-dimensional system. 
A \emph{discrete cross section} of $(X,S)$ is a closed subset $C \subset X$ such that ${X = \{S^{i}(x) \colon x \in C, i \in \N_{>0}\}}$ and the function $r_{C} \colon C \to \N_{>0}$ defined by ${r_{C}(x) = \min \{i \in \N_{>0} \colon S^{i}(x) \in C\}}$ is continuous. Note that the continuity of $r_C$ implies that $C$ will be clopen.
Associated to a discrete cross section $C \subset X$ is a return system $(C,S_{C})$ defined by
\begin{equation*}
\begin{gathered}
S_{C} \colon C \to C\\
S_{C} \colon x \mapsto S^{r_{C}(x)}x.
\end{gathered}
\end{equation*}
The systems $(X,S)$ and $(C,S_{C})$ are always flow equivalent, and there is a natural flow equivalence $\Psi \colon \Sigma_{S}X \to \Sigma_{S_C} C$ defined as follows.
For each $x \in C$, $\Psi$ maps the segment in $\Sigma_{S}X$ with endpoints $(x,0)$ and $(x,r_C(x)) = (S_{C}(x),0)$ linearly via ${t \mapsto \frac{1}{r_{C}(x)}t}$ onto the segment in $\Sigma_{S_{C}}C$ with endpoints $(x,0)$ and $(x,1) = (S_{C}(x),0) \in \Sigma_{S_{C}}X$.
That $\Psi$ is well-defined and is a flow equivalence then follows from the fact that $C$ is a discrete cross section of $(X,S)$.

Note that if $(X,S)$ is minimal, then any clopen subset is a discrete cross section, and any return system $(C,S_{C})$ is also minimal.

It follows from the above discussion that if $(C,S_{C})$ and $(D,T_{D})$ are return systems for $(X,S)$ and $(Y,T)$, then any topological conjugacy $\psi \colon (C,S_{C}) \to (D,T_{D})$ induces a flow equivalence between their suspensions, and hence a flow equivalence between the suspensions of the initial systems
\begin{equation*}
\begin{gathered}
\Psi \colon \Sigma_{S}X \to \Sigma_{T}Y \\
\Psi \colon (x,t) \mapsto (\psi(x),t).
\end{gathered}
\end{equation*}

A \emph{return equivalence} between systems $(X,S)$ and $(Y,T)$ is a triple $(\psi,C,D)$ where ${C \subset X, D \subset Y}$ are discrete cross sections and $\psi \colon (C,S_{C}) \to (D,T_{D})$ is a topological conjugacy. 

The following result is originally due to Parry and Sullivan in \cite{Parry-Sullivan_topo-invariant-flows_1975}, and is fundamental in the study of flow equivalence; we give a version which is worded differently, more in line with the proof given in \cite{Boyle-Carlsen-Eiler_Flow-isotopy_2017,
      Boyle-Carlsen-Eilers_flow-equivalence-subshifts_2017}.

\begin{Theorem}[Parry-Sullivan]
\label{thm:parrysullivanequivalence}
Between zero-dimensional systems, 
% $\Sigma_{S}X$ and $\Sigma_{T}Y$
every flow equivalence is isotopic to one which is induced by a return equivalence.
% given by return systems $(C,S_{C})$ and $(D,T_{D})$.
\end{Theorem}

More background on flow equivalence can be found in \cite{Boyle-Carlsen-Eiler_Flow-isotopy_2017,
      Boyle-Carlsen-Eilers_flow-equivalence-subshifts_2017}.

\begin{Remark}[flow codes]
% The term flow code arose from the case where $(X,S)$ and $(Y,T)$ are subshifts: the map $\psi$ is induced by a word block code; see \cite[Appendix A]{Boyle-Chuysurichay_MCG-subshift-finite-type_2018}.
When $(X,T)$ and $(Y,S)$ are subshifts, return equivalences are induced by \emph{flow codes}, a notion developed in~\cite{Boyle-Carlsen-Eiler_Flow-isotopy_2017,
      Boyle-Carlsen-Eilers_flow-equivalence-subshifts_2017} (see also the appendix in~\cite{Boyle-Chuysurichay_MCG-subshift-finite-type_2018}).
For subshifts, flow codes encode flow maps between their suspensions in an analogous way to how block codes encode continuous shift commuting maps.
\end{Remark}

\subsection{Flow equivalence of Sturmian, Denjoy and interval exchange systems}

\subsubsection*{Flow equivalence of Sturmian systems}

Fix $\alpha, \beta \in \R\setminus \Q$ and consider the associated Sturmian subshifts $(X_{\alpha},\Shift_{\alpha}), (X_{\beta},\Shift_{\beta})$ with unique invariant Borel probability measures $\mu_{\alpha}, \mu_{\beta}$. 

The classification of Sturmian systems up to flow equivalence was first proved by Fokkink in~\cite{Fokkink_structure-trajectories_1991} (although using different language); an alternative proof was given in~\cite{Barge-Williams_Denjoy-classification_2000}. 
We give an alternative short proof here.

\begin{Theorem}[flow equivalence of Sturmian systems]
\label{thm:flow-equivalence-Sturmian-systems}
For $\alpha,\beta\in \R \setminus \Q$, the following are equivalent:
\begin{enumerate}[noitemsep]
    \item The points $\alpha,\beta\in \R\Proj^1\setminus \Q\Proj^1$ belong to the same $\PGL_{2}(\Z)$-orbit.
    \item The Sturmian subshifts $(X_{\alpha},\Shift_{\alpha})$ and $(X_{\beta},\Shift_{\beta})$ are flow equivalent.
\end{enumerate}
\end{Theorem}

\begin{proof}
Recall that a flow equivalence $\Psi \colon \Sigma_{\Shift_{\alpha}}X_{\alpha} \to \Sigma_{\Shift_{\beta}}X_{\beta}$ induces by Proposition~\ref{prop:flow-coinvariants} an isomorphism $\Psi^{*} \colon \CoInv_{\Shift_{\beta}} \to \CoInv_{\Shift_{\alpha}}$ and an induced map on the spaces of positive homomorphisms $\Psi_{*} \colon \mathcal{P}(\CoInv_{\Shift_{\alpha}}) \to \mathcal{P}(\CoInv_{\Shift_{\beta}})$. 
Letting $\tau_{\alpha}$ and $\tau_{\beta}$ denote the unique states on $\mathcal{G}_{\Shift_{\alpha}}$ and $\mathcal{G}_{\Shift_{\beta}}$ respectively, there exists $\lambda > 0$ such that $\Psi_{*}(\tau_{\alpha}) = \lambda \tau_{\beta}$. 
It follows from Theorem~\ref{thm:Denjoy-coinvariants-state} that $\Z+\Z\alpha = \lambda(\Z + \Z\beta)$, and we deduce (see Lemma~\ref{lem:proportional-rank2-lattices}) that $\alpha$ and $\beta$ are $\PGL_{2}(\Z)$-equivalent.

For the converse, suppose $A\alpha = \beta$ for some $A \in \GL_{2}(\Z)$. The factor maps $H_{\alpha} \colon X_{\alpha} \to \mathbb{T}$ and $H_{\beta} \colon X_{\beta} \to \mathbb{T}$ extend to factor maps $\pi_{\alpha} \colon \Sigma_{\Shift_{\alpha}}X_{\alpha} \to \mathbb{T}^{2}, \pi_{\beta} \colon \Sigma_{\Shift_{\beta}}X_{\beta} \to \mathbb{T}^{2}$ which intertwine the flows on the suspensions with the flows on $\mathbb{T}^{2}$ in the direction of $(\alpha,1)$ and $(\beta,1)$ respectively. The automorphism $T_{A} \colon \mathbb{T}^{2} \to \mathbb{T}^{2}$ induced by $A$ lifts to a homeomorphism $\tilde{T}_{A} \colon \Sigma_{\Shift_{\alpha}}X_{\alpha} \to \Sigma_{\Shift_{\beta}}X_{\beta}$. If $\tilde{T}_{A}$ does not respect the orientation on flow lines, then we may compose this with the self-homeomorphism $\Sigma_{\Shift_{\beta}}X_{\beta} \to \Sigma_{\Shift_{\beta}}X_{\beta}$ obtained by lifting the automorphism $-\Id$ on $\mathbb{T}^{2}$ to obtain a flow equivalence.
\end{proof}

\begin{Remark}[flow equivalence preserving cyclic order on the fiber]
\label{rem:flow-equivalence-Denjoy-preserving-cyclic-order}
Irrationals $\alpha, \beta\in \R\setminus \Q$ belong to the same $\PSL_{2}(\Z)$-orbit if and only if the Sturmian subshifts $(X_{\alpha},\Shift_{\alpha})$ and $(X_{\beta},\Shift_{\beta})$ are related by a flow equivalence which preserves the cyclic order of the fibers. The orientation reversing homeomorphism of $\Torus=\R/\Z$ given by $x \mapsto 1-x$ conjugates $R_{\alpha}$ to $R_{-\alpha}$ and induces a conjugacy (hence a flow equivalence) between the Sturmian systems $(X_{\alpha},\Shift_{\alpha})$ and $(X_{-\alpha},\Shift_{-\alpha})$ which reverses the cyclic order on the fibers.
\end{Remark}

\subsubsection*{Flow equivalence of Denjoy systems}

The following generalisation of Theorem~\ref{thm:flow-equivalence-Sturmian-systems} was proved by Putnam in the context of $C^{*}$-algebras. The proof in \cite[Theorem 6]{Putnam_Morita-equivalence-Denjoy_1988} is stated just for homeomorphisms of the suspensions, but it can be seen from the proof that one can find a homeomorphism respecting orientation of the flow (see also \cite[6.1.5]{AransonBelitskyZhuzhoma1996}).

% Recall that for $\lambda \in \R$ and $A\subset \Torus$, the subset $\lambda A \subset \Torus$ is defined by $\lambda A = \{\lambda (a+k) \colon a\in \tilde{A}, k\in \Z\}$.

\begin{Theorem}[flow equivalence of Denjoy systems]
\label{thm:flow-equivalence-Denjoy-systems}
For $k\in \{0,1\}$, the Denjoy systems  $(\Denjoy_{F_k}, F_k)$ with invariants $(\alpha_k, A_k)$ (satisfying $0\in A_k$) are flow equivalent if and only if there exists 
$M= (m_{ij}) \in \GL_2(\Z)$ such that:
\begin{equation}
    \label{eq:Denjoy-flow-equivalence}
    \alpha_1 = \frac{m_{11}\cdot\alpha_0+m_{12}}{m_{21}\cdot\alpha_0+m_{22}} 
    \qquad \mathrm{and} \qquad 
    A_1 \equiv \frac{1}{m_{21}\cdot\alpha_0+m_{22}} \cdot A_0
\end{equation}
where $\equiv$ denotes equality of subsets of $\Torus$ up to rotation.

% where $\tilde{A}_k\subset \R$ is the lift of $A_k \subset \R/\Z$, and $\equiv$ denotes equality of subsets of $\R$ up to translation. 
\end{Theorem}

\begin{Remark}[flow equivalence preserving the cyclic order on fibers]
Similar to the case of Sturmian systems, for $k\in \{0,1\}$, the Denjoy systems  $(\Denjoy_{F_k}, F_k)$ with invariants $(\alpha_k, A_k)$ (satisfying $0\in A_k$) are related by a flow equivalence preserving the cyclic order on the fibers if and only if there exists $M= (m_{ij}) \in \GL_2(\Z)$ satisfying the relations in the previous Theorem~\ref{thm:flow-equivalence-Denjoy-systems} with $\det(M)>0$ and $\frac{1}{m_{21}\cdot\alpha_0+m_{22}}>0$. 

The orientation reversing homeomorphism of $\Torus=\R/\Z$ given by $x \mapsto -x$ induces a conjugacy (hence a flow equivalence) between the Denjoy systems $(\Denjoy_{\alpha, A},F_{\alpha, A})$ and $(\Denjoy_{-\alpha, -A}, F_{-\alpha, -A})$ which reverses the cyclic order on the fibers.
\end{Remark}

\subsubsection*{Flow equivalence of interval exchange systems}

Our final goal for this section is to generalize Fokkink's theorem on flow equivalence to interval exchanges, of which Sturmian systems correspond to the case of alphabets on $d=2$ letters.
The continued fraction expansion is generalized by the Rauzy-Veech induction.

Consider $(X,\Shift)\subset (\Al, \Shift)$ a minimal IES with ido on an unordered alphabet.
We saw that up to global inversion, the order on $\Al$ and the permutation $\pi$ are uniquely determined (and that $\pi$ must be irreducible, although it may have factors). %, so we may fix the order on $\Al$ and the permutation $\pi$.

Let us recall the description of all IETs with that symbolic encoding.
It follows from \cite[Section 8.1]{Yoccoz_IET-translation-surfaces_2010} that the length vectors $\lambda$ of IETs $T\colon (0,1]\to (0,1]$ having that symbolic encoding correspond to the invariant probability measures of $(X, \Shift)$.

The proof relies on Rauzy induction, which to any IET associates a path in the Rauzy diagram, and this path is infinite if and only if it satisfies ido.
Two IETs with ido yield the same path in the Rauzy diagram if and only if their symbolic encodings are subshifts on the same number of letters which are topologically conjugate, in which case their length parameters correspond to invariant measures for that system.

In short, for minimal IES on the same ordered alphabet, the topological conjugacy is equivalent (up to the order reversing involution) to having the same path in the Rauzy diagram 
(a generalisation of the fact that two real numbers are equal up to inversion if and only if they have the same continued fraction expansion).

The following Theorem~\ref{thm:flow-equivalence-IES} can thus be viewed as a generalization of Fokkink's theorem about flow equivalence of Sturmian systems to IES.

For $T\colon (\alpha_0,\alpha_d] \to (\alpha_0,\alpha_d]$ an IET with ido, the first return on a non-empty subinterval $I\subset (\alpha_0,\alpha_d]$ yields an interval exchange with ido on $d'\in \{d, d+1, d+2\}$ intervals where $d'-d\in \{0,1,2\}$ is the number of endpoints of $I$ which do not belong to the past orbits of the discontinuities $\{T^{-n}(\alpha_k)\colon n\in \N, 0\le k \le d\}$.
In particular, the cylinder-intervals defined by consecutive discontinuities of $T^{-m}$ for some $m\in \N_{>0}$ have $d'=d$. 
They correspond to the cylinders in the associated IES.

For a subshift $(X,\Shift)$ and a cylinder 
$C\subset X$, the return system $(C,\Shift_C)$ will be called a cylinder-induced system.

The following sharpens the Parry-Sullivan Theorem~\ref{thm:parrysullivanequivalence} in the case of interval exchange systems. 

\begin{Theorem}[flow equivalence of IES]
\label{thm:flow-equivalence-IES}
% Consider minimal IES with ido, which may be on different alphabets. They are flow equivalent if and only if they they admit cylinder-induced systems which are topologically conjugate.
Two minimal IES with ido are flow equivalent if and only if they admit cylinder-induced systems which are topologically conjugate.
\end{Theorem}

\begin{proof}
    Let $(X,\Shift)$ and $(Y,\Shift)$ be minimal IES with ido.
    
    Suppose first that they admit cylinder-induced IETs which are topologically conjugate. 
    By minimality these cylinder-induced systems are return systems, so their topological conjugacy yields a flow equivalence between the initial systems by Theorem~\ref{thm:parrysullivanequivalence}.
    
    Now suppose that they are flow equivalent.
    Choose corresponding IETs, and apply the zippered-rectangle construction \cite[Section 4]{Yoccoz_IET-translation-surfaces_2010} for any choice of valid height parameters, to obtain two marked translation surfaces (the marking is a horizontal segment satisfying \cite[Corollary 5.5]{Yoccoz_IET-translation-surfaces_2010}).
    By assumption, the translation surfaces are isomorphic, but the isomorphism may not preserve the marked transversal.

    Hence we may identify the translation surfaces, so that the two IETs and height parameters correspond to two transverse segments $I=(\alpha_0,\alpha_1), J=(\beta_0,\beta_1)$ (satisfying the conditions in \cite[Corollary 5.5]{Yoccoz_IET-translation-surfaces_2010}).
    Note that every non-empty open horizontal segment yields an IET by following the vertical flow, since the return time is finite by Keane's condition (see \cite[subsection 3.1]{Yoccoz_IET-translation-surfaces_2010}).
    
    Now consider any point $x\in I$ in the first transversal which does not belong to the orbit of a discontinuity, and follow the vertical flow until it first arrives to a point $y\in J$ in the second transversal, say in time $t\in \R_{\ge 0}$.
    There exists an open interval neighbourhood of $x$ denoted $(\alpha_0^{\prime}, \alpha_1^{\prime})\subset I$ which flows vertically to an interval neighbourhood of $y$ denoted $(\beta_0^{\prime}, \beta_1^{\prime})\subset J$ in the same constant time $t$.
    The vertical flow yields a conjugacy between the IETs on $(\alpha_0^{\prime}, \alpha_1^{\prime}]$ and $(\beta_0^{\prime}, \beta_1^{\prime}]$.

    Moreover, we may assume that $\alpha_0^{\prime}, \alpha_1^{\prime}$ (hence $\beta_0^{\prime}, \beta_1^{\prime}$) belong to the orbits of distinct discontinuities.
    Hence the IETs induced on $(\alpha_0^{\prime}, \alpha_1^{\prime}]$ and $(\beta_0^{\prime}, \beta_1^{\prime}]$ have the same number of intervals given by $d-1=2g-2+s$ where $g$ is the genus of the punctured surface and $s$ is the number of conical singularities.
    
    In fact we may even assume that $(\alpha_0^{\prime}, \alpha_1^{\prime}]$ corresponds to a cylinder interval for the IET on $I=(\alpha_0, \alpha_1]$, and hence $(\beta_0^{\prime}, \beta_1^{\prime}]$ will also be a cylinder interval for the IET on $J=(\beta_0, \beta_1]$.
    The associated symbolic subshifts yield cylinder-induced systems of $(X,S)$ and $(Y,T)$ that are topologically conjugate.
\end{proof}

\begin{Remark}[generalizing the class of cylinder-induced systems]
    The previous theorem also holds replacing the class of cylinder-induced systems by any class of induced systems (defined in the first order theory of IES) whose associated cross-sections satisfy the following properties:% are unions of consecutive cyclinders form a basis for the topology of $X$ and that are closed under induction.
    \begin{itemize}[noitemsep]
        \item[-] they are unions of consecutive cylinders
        \item[-] they form a basis for the topology of $X$
        \item[-] they are closed under induction
    \end{itemize}
    Special cases of such collections are given by the ``child intervals'' considered in \cite{Boshernitzan-Carroll_quadratic-IET_1997}, and the ``admissible intervals'' in \cite{Dolce-Francesco_IET-admissibility-induction_2017}.
    
    % We say that a set of admissible intervals for $T$ forms a \emph{basis of neighbourhoods} when every point in the complement $X_o\subset (\alpha_0,\alpha_d]$ of the orbits of discontinuities has a basis of neighbourhoods consisting of intervals in our family.
\end{Remark}

\begin{Remark}[cones of invariant measures are rationally equivalent]
    By Proposition~\ref{prop:flow-coinvariants}, if two minimal IES with ido on alphabets of the same cardinality are flow equivalent, then their cones of invariant measures are mapped to one another by a matrix in $\GL(d,\N)$. 
\end{Remark}

\section{2-asymptotic equivalences}
\label{sec:2-AI-equivalence}
We now introduce a certain equivalence relation on systems that will serve as a key ingredient for our definition of isogeny in the following section.

\begin{Definition}[transitive closure of factor equivalences]
\label{def:equiv-factors}
Consider a property $\mathtt{P}$ which to a factor map $\pi \colon (X,S) \to (Y,T)$ associates a truth value $\mathtt{P}(\pi) \in \{\top, \bot\}$ and such that $\mathtt{P}(\operatorname{id})=\top$.

Two systems $(Y_{1}, T_{1})$ and $(Y_2, T_{2})$ are called \emph{elementary-$\mathtt{P}$-equivalent} when there exists a system $(X,S)$ and factor maps $\pi_i \colon (X,S)\to (Y_i, T_i)$ satisfying $\mathtt{P}$.
% $\pi_1, \pi_2$ satisfying $P$ :
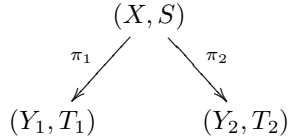
\begin{equation*}
\begin{aligned}
\xymatrix@C-2pc{
 & (X,S) \ar_{\pi_{1}}[dl] \ar^{\pi_{2}}[dr]& \\
(Y_1,T_1) & & (Y_2,T_2) \\
}
\end{aligned}
\end{equation*}
Thus elementary-$\mathtt{P}$-equivalence is a  reflexive and symmetric relation on systems.
Its transitive closure defines an equivalence relation on systems, called \emph{$\mathtt{P}$-equivalence}.
\end{Definition}

\subsection{2-asymptotic factor maps}

Let $(X,S)$ be a system. We say two points $x,x^{\prime} \in X$ are \emph{asymptotic}, denoted $x\Bumpeq x'$, when $d(S^{n}(x),S^{n}(x^{\prime})) \to 0$ as $n \to -\infty$ and $n\to +\infty$.
This is an equivalence relation.

Note that for a factor map $\pi \colon (X,S)\to (Y,T)$ we have $x\Bumpeq x' \implies \pi(x)\Bumpeq \pi(x')$.

\begin{Definition}[2-asymptotic]
\label{def:2-asymp}
A factor map $\pi \colon (X,S) \to (Y,T)$ is called \emph{2-asymptotic} (2-A) when there exists a countable subset of \emph{blown-up points} $B \subset Y$ such that
\begin{enumerate}[noitemsep]
\item For $y\in Y$, we have $\Card{\pi^{-1}(y)}> 1$ if and only if $y \in B$.
\item For every $y\in B$, we have $\Card{\pi^{-1}(y)} = 2$.
\item For every $y \in B$, the pair of points in $\pi^{-1}(y)=\{x,x'\}$ are asymptotic.
\end{enumerate}
% We denote by $A=\pi^{-1}(B)$ the set of points belonging to an asymptotic pair which gets contracted by $\pi$.
\end{Definition}

Note that a topological conjugacy is a 2-A factor, since one may choose the set $B$ in the definition to be empty.

\begin{comment}
\begin{Remark}[finite number of asymptotic pairs for subshifts]
When $(X,S)$ is a subshift, a 2-A factor map can only identify finitely many pairs of asymptotic $S$-orbits (by compactness argument).    
\end{Remark}
\end{comment}

\begin{comment}
% We may restrict this equivalence relation between minimal Cantor systems, where we require the $(Y_i, S_i)$ to be minimal, but not necessarily the $(X_j, T_j)$.
% Note however that in a chain of elementary 2-A equivalencies between minimal Cantor systems $(Y_i, S_i)$, we do not require that the systems $(X_j, T_j)$ on top to be minimal.
\end{comment}

\begin{Remark}[2-A factors between minimal Cantor systems]
For factors $\pi\colon (X,S) \to (Y,T)$ between minimal Cantor systems, the definition of 2-A amounts to that found in~\cite{GPS_asymptotic-index_2001}: $\pi$ is at most $2$-to-$1$, and $ \forall \epsilon > 0$ the set $\{(x_{1},x_{2}) \in X \times X \colon \pi(x_{1}) = \pi(x_{2}), \; d(x_{1},x_{2}) \ge \epsilon\}$ is finite.
\end{Remark}

\begin{Remark}[2-A factors induce isomorphism of invariant measures]
\label{remark:2afactormeasureiso}
If $\pi \colon (X,S) \to (Y,T)$ is a 2-A factor map between minimal Cantor systems, then the pushforward map on measures $\pi_{*} \colon \mathcal{M}(S) \to \mathcal{M}(T)$ is bijective; in particular it sends ergodic measures to ergodic measures.
Indeed, a factor map always induces a surjective pushforward map on the space of measures (see for instance~\cite[3.11]{DenkerGrillenbergerSigmundBook}), and it is injective because $\pi$ is injective off a countable set and no $S$-invariant measure can have atoms by minimality.
\end{Remark}

Consider a 2-A factor $\pi \colon (X,S) \to (Y,T)$. 
Choose a representative $y\in B\bmod{T}$ in every blown-up orbit, and an order on its preimage $\pi^{-1}(y)=\{x, x'\}$.
For such an asymptotic pair $(x, x')$, there is a group homomorphism $C(X,\Z) \to \Z$ defined by
\begin{equation*}
    % \Index_{y} \colon C(X,\Z) \to \Z
    % \qquad \mathrm{by} \qquad
    f \mapsto \sum_{n \in \Z}f(S^{n}(x))-f(S^{n}(x'))
\end{equation*}
which only depends on the orbit $y\bmod{T}$ and the order of its preimage $(x, x')\bmod{S}$.
This homomorphism vanishes on coboundaries, so it induces a homomorphism $\Index_{y} \colon \CoInv_{S} \to \Z$.
For every $f \in C(X,\Z)$, there are only a finite number of $y\in B\bmod{T}$ such that $\Index_y(f)\ne 0$, hence following \cite{GPS_asymptotic-index_2001}, we may define the \emph{index homomorphism}
\begin{equation*}
    \Index \colon \CoInv_{S} \to \bigoplus_{B\bmod{T}} \Z
    \qquad \mathrm{by} \qquad
    \Index = \bigoplus_{y\in B\bmod{T}} \Index_{y}.
\end{equation*}

Recall that for any factor $\pi\colon (X,S)\to (Y,T)$ between minimal Cantor systems, the map $\pi^*\colon \CoInv_{T}\to \CoInv_{S}$ is injective. When $\pi$ is 2-A, the cokernel of $\pi^*$ is identified with the index map by \cite[Theorem 3.1]{GPS_asymptotic-index_2001} as follows.

\begin{Theorem}[index exact sequence]
\label{thm:GPS_im(pi*)=ker(index)}
Let $\pi \colon (X,S) \to (Y,T)$ be a 2-A factor map between minimal Cantor systems. The index homomorphism $\Index\colon \CoInv_{S} \to \bigoplus_{B\bmod{T}}\Z$ is surjective and its kernel is the image of the induced homomorphism $\pi^{*}\colon \CoInv_{T} \to \CoInv_{S}$. In other words, there is a short exact sequence of abelian groups:
\begin{equation*}
    0 \to \CoInv_{T} \xrightarrow{\pi^*} \CoInv_{S} \xrightarrow{\Index} \bigoplus_{B\bmod{T}}\Z \to 0.
\end{equation*}
% In formula $\forall [f] \in \CoInv_{S}, \quad \Index(f) = 0 \iff [f] \in \im(\pi^{*})$.
\end{Theorem}

\subsection{Infinitesimal 2-asymptotic equivalence}

\begin{Definition}[infinitesimal] 
We say that a 2-A factor $\pi \colon (X,S) \to (Y,T)$ with index homomorphism $\Index$ is \emph{infinitesimal} when $\Index(\CoInv_{S})/\Index(\Inf_{S})$ is torsion, namely $$\left(\Index(\CoInv_{S}) / \Index(\Inf_{S})\right) \otimes \Q = 0.$$
\end{Definition}

Note that if $\pi$ is a topological conjugacy, the associated index map is the zero map, and hence $\pi$ is an infinitesimal 2-A factor map. Thus we may define the notion of infinitesimal 2-A equivalence.

\begin{Remark}[full-rank implies infinitesimal]
    Consider a 2-A factor $\pi \colon (X,S) \to (Y,T)$ with index homomorphism $\Index$. If the restriction $\Index \colon \Inf_{S} \to \bigoplus_{B\bmod{T}}\Z$ has torsion cokernel then $\pi$ is infinitesimal, namely $\Index(\CoInv_{S}) / \Index(\Inf_{S})$ is torsion. When $(X,S)$ and $(Y,T)$ are minimal Cantor systems, these conditions are equivalent, since in this case $\Index \colon \CoInv_{S} \to \bigoplus_{B\bmod{T}}\Z$ is surjective by Theorem~\ref{thm:GPS_im(pi*)=ker(index)}.
\end{Remark}

The following propositions connect infinitesimal 2-A factors with the rational images of respective states.
\begin{Proposition}[2-AI factors and states]
\label{prop:2A_infinitesimal-implies-0-measured}
Consider an infinitesimal 2-asymptotic factor $\pi \colon (X,S) \to (Y,T)$ whose index homomorphism $\Index$ satisfies $\ker(\Index)\subset \pi^*(\CoInv_{T})$ (as it is the case when $(X,S)$ and $(Y,T)$ are minimal Cantor systems by Theorem~\ref{thm:GPS_im(pi*)=ker(index)}). Then for any state $\tau_X$ of $(X,S)$, its push-forward state $\tau_Y=\pi_*(\tau_X)$ of $(Y,T)$ satisfies
\begin{equation*}
    \im(\tau_{X}) \otimes \Q = \im(\tau_{Y}) \otimes \Q.
\end{equation*}

% When both $(X,S)$ and $(Y,T)$ are uniquely ergodic, this relation holds between the states $\tau_{\mu_X}$ and $\tau_{\mu_Y}$ associated to their unique invariant probability measures. 
\end{Proposition}

\begin{proof}
% By definition of $\tau_Y=\pi_*(\tau_X)$ we have $\im(\tau_Y)\subset \im(\tau_X)$.
Let $[f] \in \CoInv_{S}$. 
Since $\Index(\CoInv_{S})/\Index(\Inf_{S})$ is torsion, we may choose $r\in \N_{>0}$, and $g \in \Inf_{S}$ such that $r\Index(f) = \Index(g)$.
Thus $\Index(rf - g) = 0$ so the assumption $\ker \Index = \im(\pi^*)$ ensures that we may choose $h \in C(Y,\Z)$ such that $\pi^{*}(h) = rf - g$. 

Now consider a state $\tau_X$ of $(X,S)$ and let $\tau_Y=\pi_*(\tau_X)$ be the push-forward state of $(Y,T)$.
By definition of $\Inf_{S}$ we have $\tau_X(g)=0$ thus 
\begin{equation*}
    \tau_{Y}(h) = \pi_{*}(\tau_{X})(h) = \tau_{X}(\pi^{*}(h)) = \tau_{X}(rf - g)=r\tau_X(f)
\end{equation*}
proving that $\im(\tau_X)\otimes \Q = \im(\tau_Y)\otimes \Q$ as desired.
\end{proof}

% The following Proposition proves a converse under the stronger assumptions that $(X,S)$ is uniquely ergodic.

\begin{Proposition}[2-AI factors of UE systems]
\label{prop:2A-UE_0-measured-implies_infinitesimal}
Consider a 2-asymptotic factor $\pi\colon (X,S)\to (Y,T)$ whose index homomorphism $\Index$ satisfies $\ker(\Index)\supset \pi^*(\CoInv_{T})$ 
(as it is the case when $(X,S)$ and $(Y,T)$ are minimal Cantor systems by Theorem~\ref{thm:GPS_im(pi*)=ker(index)}). Assume that $(X,S)$ is uniquely ergodic with unique state $\tau_X$, and let $\tau_Y=\pi_*(\tau_X)$ be its push-forward state on $(Y,T)$. If $\im(\tau_X)\otimes \Q = \im(\tau_Y)\otimes \Q$ then the 2-A factor $\pi$ is infinitesimal.
\end{Proposition}

\begin{proof}
    To show that $\Index(\CoInv_{S})/\Index(\Inf_{S})$ is torsion, we fix $b\in \Index(\CoInv_{S})$ and will construct $f'\in \Inf_S$ and $d\in \N_{>0}$ such that $\Index(f')=db$.
    
    We may choose $f\in \CoInv_{S}$ with $\Index(f)=b$.
    By assumption $\im(\tau_X)\otimes \Q=\im(\tau_Y)\otimes \Q$, so there exists $d\in \N_{>0}$ and $g\in \CoInv_{T}$ such that $\tau_X(f)=\tfrac{1}{d}\tau_Y(g)$.
    Let $f'=df-\pi^*(g)\in \CoInv_{S}$.
    On the one hand $$\tau_X(f')=d\tau_X(f)-\tau_X(\pi^*(g))=d\tau_X(f)-\pi_*(\tau_X)(g)=d\tau_X(f)-\tau_Y(g)=0$$ so $f'\in \ker(\tau_X)$ hence $f'\in \Inf_{S}$ by unique ergodicity of $(X,S)$.
    On the other hand $\Index(f')=d\Index(f)$ since $\pi^*(g)\in \pi^*(\CoInv_{T})\subset \ker \Index$ by assumption.
\end{proof}

Let us sum up the previous results in a particular case of interest.

\begin{Corollary}[2-asymptotic factors of uniquely ergodic minimal Cantor systems]\label{cor:2-asym-factors-ue-minimal-cantor}
    % Consider a 2-asymptotic factor map $\pi \colon (X,S)\to (Y,T)$ between minimal Cantor systems which are uniquely ergodic with states $\tau_X, \tau_Y= \pi_*(\tau_X)$.
    % The factor $\pi$ is infinitesimal if and only if $\im(\tau_X)\otimes \Q=\im(\tau_Y)\otimes \Q$.

    Consider uniquely ergodic minimal Cantor systems $(X,S), (Y,T)$ with states $\tau_X, \tau_Y$. Both of the following hold:

    \begin{enumerate}
    \item
    If there exists a 2-AI factor $\pi \colon (X,S) \to (Y,T)$, then $\im(\tau_X)\otimes \Q=\im(\tau_Y)\otimes \Q$.
    \item
    If $\im(\tau_X)\otimes \Q=\im(\tau_Y)\otimes \Q$ then every 2-A factor $\pi \colon (X,S) \to (Y,T)$ is infinitesimal.
    \end{enumerate}
    % If there is a 2-asymptotic factor map $\pi \colon (X,S)\to (Y,T)$, then it is infinitesimal if and only if $\im(\tau_X)\otimes \Q=\im(\tau_Y)\otimes \Q$.
    Thus if $(X,S)$ and $(Y,T)$ are 2-AI equivalent, then $\im(\tau_{X}) \otimes \Q = \im(\tau_{Y}) \otimes \Q $. 
\end{Corollary}

2-AI extensions can also arise from maximal equicontinuous factors of certain substitution subshifts, as in the following example.
%These examples are of interest in the study of aperiodic tilings associated to one-dimensional substitutions, such as the following.

\begin{Example}[period doubling substitution]
Let $(X,\Shift_{X})$ be the minimal subshift associated to the period doubling substitution defined by $\tau(a) = ab$, $\tau(b) = aa$. The group of continuous eigenvalues of the system $(X,\Shift_{X})$ is isomorphic to $\Z[\tfrac{1}{2}]$ giving a factor onto the $2$-adic odometer $\pi \colon (X,\Shift_{X}) \to (\mathcal{O}_{2},S)$ (see~\cite[§1.6]{Fogg_Subsitutions-in-dyn-arit-comb_2002}). This is a 2-AI factor, which can be seen from the calculations given in~\cite[§5.3.2]{Barge-Kallendonk-Schmieding_factors-tiling_2012}.    
\end{Example}

\subsection{\texorpdfstring{$\M$-measured equivalence}{M-equivalence}}

In this subsection, we fix a subgroup $\M\subset \R$ and generalise the notion of infinitesimal 2-A equivalence to that of $\M$-measured equivalence.

% This is not the main focus of our work, but serves to elucidate the relation between the index subgroup of a 2-A factor and the group of infinitesimals.

\begin{Definition}[$\M$-measured]
\label{def:M-measured-state}
Let $\pi \colon (X,S) \to (Y,T)$ be a 2-asymptotic factor. 
We say that $\pi$ is \emph{$\M$-measured with respect to a subset of states} $\mathcal{T} \subset \PP(X,S)$ when there exists a subgroup $G \subset \CoInv_{S}$ such that $\Index(\CoInv_{S})/\Index(G)$ is torsion and every $\tau_{X}\in \mathcal{T}$ satisfies $\tau_{X}(G) \subset \M$. In particular when $\mathcal{T}$ consists of
\begin{enumerate}[noitemsep]
\item a single state $\tau_X$ we say that $\tau_X$ is $\M$-measured with respect to $\pi$.
\item the whole set $\PP(X,S)$ we say that $\pi$ is \emph{$\M$-measured}.
\item the set of ergodic measures, we say that $\pi$ is \emph{ergodically $\M$-measured}.
\end{enumerate}
\end{Definition}

Note that for any subgroup $\M \subset \R$ a topological conjugacy is always $\M$-measured, and we may consider the notion of $\M$-measured 2-A equivalence between uniquely ergodic systems.

\begin{Remark}[infinitesimal $\iff$ $0$-measured]
A 2-asymptotic factor is infinitesimal if and only if it is $\{0\}$-measured. One direction of this is obvious: set $G = \Inf(\CoInv_{S})$. For the other direction, suppose $\pi$ is $\{0\}$-measured and let $G \subset \CoInv_{S}$ be the subgroup given by the definition. Then $\tau(G) = 0$ for every state $\tau\in\PP(X,S)$, so $G \subset \Inf(\mathcal{G}_{S})$. Since $\mathcal{I}(\mathcal{G}_{S})/\mathcal{I}(G)$ is torsion and $G \subset \Inf(\mathcal{G}_{S})$ it follows that $\mathcal{I}(\mathcal{G}_{S}) / \mathcal{I}(\Inf(\mathcal{G}_{S}))$ is torsion as well.

% Note that if $\M\subset \M'$, we have $\M$-measured implies $\M'$-measured.
\end{Remark}

\begin{Proposition}[$\M$-measured and states]
\label{prop:M-measured-state}
Consider a 2-asymptotic factor $\pi \colon (X,S) \to (Y,T)$ with index homomorphism $\Index$.
Fix a state $\tau_X$ of $(X,S)$ and let $\tau_Y=\pi_*(\tau_X)$ be its push-forward state on $(Y,T)$.
Let $\M\subset \R$ be a subgroup. When $\ker(\Index)\subset \pi^*(\CoInv_{T})$, we have the implication:
\begin{equation*}
    \text{$\tau_X$ is $\M$-measured} \implies (\im(\tau_{X})+\M) \otimes \Q = (\im(\tau_{Y})+\M) \otimes \Q.
\end{equation*}
When $\ker(\Index)\supset \pi^*(\CoInv_{T})$, we have the implication:
\begin{equation*}
    \text{$\tau_X$ is $\M$-measured} \impliedby (\im(\tau_{X})+\M) \otimes \Q = (\im(\tau_{Y})+\M) \otimes \Q.
\end{equation*}
Consequently, when $\ker(\Index)=\im(\pi^*)$, as it is the case when $(X,S)$ and $(Y,T)$ are minimal Cantor systems by Theorem~\ref{thm:GPS_im(pi*)=ker(index)}, we have the equivalence:
\begin{equation*}
    \text{$\tau_X$ is $\M$-measured} \iff (\im(\tau_{X})+\M) \otimes \Q = (\im(\tau_{Y})+\M) \otimes \Q.
\end{equation*}
\end{Proposition}

\begin{proof}

Suppose that $\ker(\Index)\subset \pi^*(\CoInv_{T})$ and that $\pi$ is $\M$-measured with respect to $\tau_{X}$.
Let $G\subset \CoInv_{S}$ be a subgroup such that $\Index(\CoInv_{S})/\Index(G)$ is torsion and $\tau_X(G)\subset \M$.
For every $[f] \in \CoInv_{S}$, we may choose $r\in \N_{>0}$, and $g \in G$ such that $r\Index(f) = \Index(g)$. 
Thus $\Index(rf - g) = 0$ so we may choose $h \in C(Y,\Z)$ such that $\pi^{*}(h) = rf - g$. 
By assumption $\tau_X(g)\in \M$ so 
\[\tau_{Y}(h) = \pi_{*}(\tau_{X})(h) = \tau_{X}(\pi^{*}(h)) = \tau_{X}(rf - g)=r\tau_X(f)-\tau_X(g)\in r\tau_X(f)+\M\]
proving $(\im(\tau_{X})+\M) \otimes \Q = (\im(\tau_{Y})+\M) \otimes \Q$ as desired.

Suppose that $\ker(\Index)\supset \pi^*(\CoInv_{T})$ and that $(\im(\tau_X)+\M)\otimes\Q = (\im(\tau_Y)+\M)\otimes\Q$.
Define the subgroup $G=\tau_X^{-1}(\M)\subset \CoInv_{S}$ and let us show that $\Index(\CoInv_{S})/\Index(G)$ is torsion.
Let $b\in \Index(\CoInv_{S})$ and choose $f\in \CoInv_{S}$ with $\Index(f)=b$.
Since $(\im(\tau_X)+\M)\otimes \Q=(\im(\tau_Y)+\M)\otimes \Q$, there exists $d\in \N_{>0}$, $m\in \M$ and $g\in \CoInv_{T}$ such that $\tau_X(f)=\tfrac{1}{d}(\tau_Y(g)+m)$.
Let $f'=df-\pi^*(g)\in \CoInv_{S}$.  
On the one hand $$\tau_X(f')=d\tau_X(f)-\tau_X(\pi^*(g))=d\tau_X(f)-\pi_*(\tau_X)(g)=d\tau_X(f)-\tau_Y(g)=m$$ thus $f'\in \tau_X^{-1}(\M)=G$.
On the other hand $\Index(f')=d\Index(f)$ since $\pi^*(g)\in \pi^*(\CoInv_{T})\subset \ker \Index$ by assumption.
\end{proof}

\begin{comment}
\begin{Remark}[$0$-dimensional]
\label{rem:0-dimensional}
    The notion of asymptotic factor $\pi\colon (X,S)\to (Y,T)$, especially its attributes infinitesimal and $\M$-measured, mostly apply to $0$-dimensional systems for which the group of coinvariants equals the $K$-theory group in degree $0$.

    We believe that the appropriate definitions for the more general setting should involve the maps induced on $K$-theory, but we preferred to avoid such generalities which were not needed for the current applications.
\end{Remark}
\end{comment}

\begin{Definition}[$\M$-equivalence]
    We say that two systems are ergodically $\M$-equivalent when they are related by a sequence of 2-A factors that are ergodically $\M$-measured.
\end{Definition}

\begin{Remark}[applications of $\M$-equivalence to classify non-uniquely ergodic systems]
    We will not need to refer to the definition of ergodic $\M$-equivalence in this work, since most of our applications will concern totally uniquely ergodic systems. However, we believe that it may be a useful notion to classify systems up to isogenies, in particular interval exchanges with ido which may not be totally uniquely ergodic (see Conjecture~\ref{conj:isogenies-IES=Rational-Invariants}).  
\end{Remark}

\subsection{2-A factors of Denjoy and Sturmian systems}

In this subsection, we describe all 2-A factors of Denjoy systems (up to topological conjugacy), and characterize which ones are infinitesimal or $\M$-measured.

The subsets of $\R$ containing $0$ can be added and scaled by real numbers. 
Hence the subsets of $\Torus = \R/\Z$ containing $0$ can be added and multiplied by a real number, by first lifting them along $\R \to \R/\Z$, then adding them or scaling them, and projecting back down along $\R\to \R/\Z$.
For a subring $\K\subset \R$ and a subset $A \subset \Torus$ containing $0$, the $\K$-submodule of $\Torus$ generated by $A$ is thus:
\begin{align*}
    \K{A}
    &= \{\textstyle{\sum}_i \lambda_i ( \tilde{a_i}+k_i) \bmod{1} \colon \tilde{a_i}\in \tilde{A}\cap [0,1), k_i\in \Z, \lambda_i\in \K\} \\
    &=
    \{\lambda_0+\textstyle{\sum}_i \lambda_i \tilde{a_i} \bmod{1} \colon \tilde{a_i}\in \tilde{A}\cap [0,1), \lambda_i\in \K\} \subset \Torus.
\end{align*}
In particular $\K{A} = \K{(A+\K/\Z)}$.

Fix an irrational $\alpha \in \Torus$ and an $R_\alpha$-invariant countable subset $A\subset \Torus$ containing $0$.
The factor $\pi_{\alpha, A} \colon (\Denjoy_{\alpha, A}, F_{\alpha, A}) \to (\Torus, R_{\alpha})$ is 2-asymptotic.
Each point $a \in A \subset \Torus$ has two preimages $a^\pm \in \Denjoy_{\alpha, A}$ which form a naturally ordered asymptotic pair for $\pi_{\alpha, A}$, and all asymptotic pairs arise this way.

Now for every $R_\alpha$-invariant subset $A'\subset A$ containing $0$, there is a factor map $(\Denjoy_{\alpha, A}, F_{\alpha, A}) \to (\Denjoy_{\alpha,A}, F_{\alpha, A'})$ obtained by identifying pairs of points $a^-, a^+\in \Denjoy_{\alpha,A}$ such that $\pi_{\alpha,A}(a^-) = \pi_{\alpha, A}(a^+)\in A\setminus A'$, and these are all its asymptotic pairs.

\begin{Lemma}[2-A factor of Denjoy]
\label{lemma:2AfactorsofDenjoyrotationnumbers}
    There is a 2-A factor of Denjoy systems $\pi \colon (\Denjoy_{\alpha,A},F_{\alpha, A}) \to (\Denjoy_{\alpha',A^{\prime}},F_{\alpha',A^{\prime}})$ if and only if there is $\epsilon \in \{\pm 1\}$ such that $\alpha' = \epsilon \alpha \bmod{\Z}$ and $A'\subset \epsilon A \subset \Torus$.
\end{Lemma}
\begin{proof}
    Let us show the implication (only if) as the converse (if) is clear.
    Since a Denjoy system admits a 2-A factor onto a rotation of the torus which is unique up to conjugacy, we have that $\pi_{\alpha, A}$ is conjugate to $\pi_{\alpha', A'} \circ \pi$ hence $\alpha'=\pm \alpha\bmod{\Z}$.
    Hence, up to composing with the map $x \mapsto -x$ on $\mathbb{T}$, the image of nontrivial asymptotic points in $\Denjoy_{\alpha',A^{\prime}}$ under $\pi_{\alpha^{\prime},A^{\prime}}$ must be contained in the image of nontrivial asymptotic points in $\Denjoy_{\alpha,A}$ under $\pi_{\alpha,A}$, and hence $A^{\prime} \subset \pm A$. 
\end{proof}

% We now characterise which ones are infinitesimal or $\M$-measured for a subgroup $\M\subset \R$.

\begin{Remark}[factor from Denjoy to circle-rotation is not 2-AI but may be $\M$-measured]  
The 2-A factor $\pi_{\alpha, A} \colon (\Denjoy_{\alpha, A}, F_{\alpha, A}) \to (\Torus, R_{\alpha})$ of uniquely ergodic systems is not infinitesimal. For a subgroup $\M\subset \R$, it is $\M$-measured if and only if $\Q{A}\subset \Q \otimes (\Z+\M) \bmod{\Z}$.

Let us prove this remark for completeness, as it involves computations shedding light on the structure of Denjoy systems, akin to those in the proofs of \cite[Theorem~\ref{thm:Denjoy-coinvariants-state}]{PSS_C-star-Denjoy_1986} and \cite[Theorem~\ref{thm:GPS_im(pi*)=ker(index)}]{GPS_asymptotic-index_2001}.
%
% Note that the system $(\Torus, R_\alpha)$ is not zero-dimensional, and the analysis of the map $\pi_{\alpha, A}$ will explain part of Remark~\ref{rem:0-dimensional}.
%
We cannot invoke Theorem~\ref{thm:GPS_im(pi*)=ker(index)}, but only one of the implications from Proposition~\ref{prop:M-measured-state} as $\ker(\Index)\supsetneq \im(\pi_{\alpha, A}^*)$.
\end{Remark}

\begin{proof}
Choose representatives $a_i\in A\bmod{R_\alpha}$ with $a_1=\alpha$. Let us describe the index homomorphism $\Index \colon \CoInv_{\alpha, A} \to \bigoplus_{A\bmod{R_{\alpha}}} \Z$ of $\pi_{\alpha, A}$.

Recall the constructions described in the proof of Theorem~\ref{thm:Denjoy-coinvariants-state}.
The group $\CoInv_{\alpha, A}$ has a basis consisting of the classes of $1$ and $\chi_{a_i}^-$ for our representatives $a_i\in A\bmod{R_\alpha}$, on which the trace evaluates as $\tau_{\alpha, A}(1)=1$ and $\tau_{\alpha, A}(\chi_{a_i}^-)=a_i$, and the subgroup $\Inf_{\alpha, A}$ corresponds to $\Z$-linear relations among these values.
The sum defining the homomorphism $\Index_{a_i} \colon \CoInv_{\alpha, A} \to \Z$ yields $\Index_{a_i}(\chi_{a_j}^-)=\delta_{i,j}$. %(since $a_i = a_j \bmod{R_\alpha}$ if and only if $i=j$).
Hence $\Index\colon \CoInv_{\alpha, A}\to \bigoplus_{A\bmod {R_\alpha}} \Z$ restricts to an isomorphism on $\bigoplus \chi_{a_i}^-$, so $\Index(\CoInv_{\alpha, A})/\Index(\Inf_{\alpha, A})=\left(\bigoplus \Z\cdot a_i \right)/\Inf_{\alpha, A}$. This is an abelian group with positive rank since $A$ contains irrational numbers, thus it is not torsion, proving that $\pi_{\alpha, A}$ is not infinitesimal.

Now consider a subgroup $\M\subset \R$. 
Observe that every continuous function in $C(\Torus, \Z)$ and the subgroup of $R_\alpha$-coboundaries is trivial, whence $\CoInv_{R_\alpha}=\Z$ and the trace has image $\tau_{\alpha}=\Z$.
The injective homomorphism $\pi_{\alpha, A}^* \colon \Z \to \Z+\bigoplus_{i}\Z\chi_{a_i}^-$ restricts to a bijection $\Z\mapsto \Z$ between the subgroups of constant functions. 
In particular $\im(\pi_{\alpha, A}^*)\subset \ker(\Index)$.
Since $(\Denjoy_{\alpha, A}, F_{\alpha, A})$ is uniquely ergodic, the Proposition~\ref{prop:M-measured-state} implies that if $\pi_{\alpha, A}$ is $\M$-measured, then $(\M + \im(\pi_{\alpha, A}))\otimes \Q = (\M+\Z)\otimes \Q$.
Conversely, assume $(\M + \im(\pi_{\alpha, A}))\otimes \Q = (\M+\Z)\otimes \Q$. Let us show that the subgroup $G=\tau_{\alpha, A}^{-1}(\M)$ of $\CoInv_{\alpha, A}$ is such that $\Index(\CoInv_{\alpha, A})/\Index(G)$ is torsion. 
Note that $G$ contains the subgroup $\Z$ of constants if and only if $\M$ contains $\Z$.
Since $\Index$ is an isomorphism in restriction to the complement $\bigoplus \chi_{a_i}^-$, we must show that $\CoInv_{\alpha, A}/(\Z+G)$ is torsion.
 After dividing by $\ker(\Inf_{\alpha, A})\subset G$, it is enough to show that $\tau_{\alpha, A}(\CoInv_{\alpha, A})/\tau_{\alpha, A}(\Z+G)= \im(\tau_{\alpha, A})/(\Z+\M)$ is torsion, which follows from the assumption.
\end{proof}

Now let us reap from the previous two subsections the corollaries of main interest.

\begin{Corollary}[which 2-A factors of Denjoy are infinitesimal and $\M$-measured]
\label{cor:2-AI-factor-Denjoy}
The 2-A factor $\pi \colon (\Denjoy_{\alpha,A},F_{\alpha, A}) \to (\Denjoy_{\alpha,A^{\prime}},F_{\alpha,A^{\prime}})$ is infinitesimal if and only if the $\Q$-modules $\Q{A'} \subset \Q{A}$ are equal. Moreover, if $\M\subset \R$ is a subgroup, then the 2-A factor $\pi$ is $\M$-measured if and only if the subsets $A'\cup \M\subset A\cup \M$ generate the same $\Q$-submodules of $\Torus$.
\end{Corollary}
\begin{proof}
The Denjoy systems are uniquely ergodic: let $\tau$ and $\tau'$ be their unique states so that $\pi_*(\tau)=\tau'$.
They are minimal Cantor so by Theorem~\ref{thm:GPS_im(pi*)=ker(index)} the index homomorphism $\Index$ of $\pi$ has kernel $\ker(\Index)=\im(\pi^*)$.

Hence the hypothesis of Propositions~\ref{prop:2A_infinitesimal-implies-0-measured},~\ref{prop:2A-UE_0-measured-implies_infinitesimal} and Proposition~\ref{prop:M-measured-state} are satisfied.
Thus $\pi$ is infinitesimal if and only if $\im(\tau)\otimes \Q=\im(\tau')\otimes \Q$ and it is $\M$-measured if and only if $(\im(\tau)+\M)\otimes \Q = (\im(\tau')+\M)\otimes \Q$.

It follows from Theorem~\ref{thm:Denjoy-coinvariants-state} that $\im(\tau) \bmod{\Z}$ and $\im(\tau') \bmod{\Z}$ are the $\Z$-submodules of $\Torus$ respectively generated by $A$ and $A'$, which concludes the proof. 
\end{proof}

\begin{Corollary}[2-AI equivalence preserving the class of Denjoy systems]
\label{cor:2-AI-equiv-Denjoy}
The Denjoy systems $(\Denjoy_{\alpha, A}, F_{\alpha, A})$ and $(\Denjoy_{\beta, B}, F_{\beta, B})$ are related by 2-AI equivalencies preserving the class of Denjoy systems if and only if $\alpha=\pm \beta \bmod{\Z}$ and $\Q A = \Q B$.
\end{Corollary}
\begin{proof}
Suppose there is $\epsilon \in \{\pm 1\}$ such that $\alpha = \epsilon \beta$ and $\Q A = \Q B$.
Up to conjugacy, we may replace $(\beta, B)$ by $(\epsilon\beta, \epsilon B)$.
We have $\Q A = \Q (A \cup B) = \Q B$ so the 2-asymptotic factors $\pi_{A} \colon (\Denjoy_{\alpha, A \cup B}, F_{\alpha, A \cup B}) \to (\Denjoy_{\alpha, A}, F_{\alpha, A})$ and $\pi_{B} \colon (\Denjoy_{\alpha, A \cup B}, F_{\alpha, A \cup B}) \to (\Denjoy_{\alpha, B}, F_{\alpha, B})$ are both infinitesimal by Corollary~\ref{cor:2-AI-factor-Denjoy}.

Conversely suppose that the Denjoy systems are 2-AI equivalent through Denjoy systems. By Lemma~\ref{lemma:2AfactorsofDenjoyrotationnumbers}, each 2-AI factor (blowup or blowdown) preserves the rotation number up to a sign, hence $\alpha = \pm \beta \bmod{\Z}$. That $\Q{A} = \Q{B}$ then follows from Theorem~\ref{thm:Denjoy-coinvariants-state} and Corollary~\ref{cor:2-asym-factors-ue-minimal-cantor}. 
% Each 2-AI factor (blowup or blowdown) yields an inclusion of subsets in $\Torus$ and we may thus consider the subset $C\subset \Torus$ consisting of of all blownup points along the chain of 2-AI equivalencies.
%We thus find 2-AI factors $(\Denjoy_{\alpha, A\cup C}, F_{\gamma, A\cup C}) \to (\Denjoy_{\alpha, A\cup C}, F_{\gamma, A\cup C})$ and $(\Denjoy_{\beta, B\cup C}, F_{\beta, B\cup C}) \to (\Denjoy_{\beta, B}, F_{\beta, B})$.
\end{proof}

\begin{Corollary}[2-AI factors of powers of Denjoy systems]
\label{cor:2-AI-factor-Denjoy-power-m}
Consider an irrational $\alpha\in \Torus$ and an $R_\alpha$-invariant countable subset $A \subset \Torus$ containing zero. 
Fix $m \in \N_{> 0}$ and let $A^{\prime}$ be a set containing one point from each $R_{\alpha}$-orbit in $A$ such that $0 \in A^{\prime}$ and let ${A^{\prime}(m) = \{R_{m\alpha}^{k}(x) \colon x \in A^{\prime}, k \in \Z \}}$. There is a 2-AI factor between Denjoy systems $(\Denjoy_{\alpha, A}, F_{\alpha, A}^{m}) \to (\Denjoy_{m\alpha, A_{m}^{\prime}}, F_{m\alpha, A^{\prime}(m)})$. 
Hence the Sturmian system $(X_{m \alpha},\Shift_{m \alpha})$ is 2-AI equivalent to the system $(X_{\alpha}, \Shift_{\alpha}^{m})$. 
\end{Corollary}
\begin{proof}
First observe that $(\Denjoy_{\alpha, A}, F_{\alpha, A}^{m})$ is topologically conjugate to $(\Denjoy_{m\alpha, A}, F_{m\alpha, A})$.
The factor map $\pi \colon (\Denjoy_{m\alpha, A}, F_{m\alpha, A}) \to (\Denjoy_{m\alpha, A^{\prime}(m)}, F_{m\alpha, A^{\prime}(m)})$ obtained by identifying all pairs of points $x^{-},x^{+} \in \Denjoy_{m\alpha,A}$ such that $\pi_{m\alpha,A}(x^{-})=\pi_{m\alpha,A}(x^{+}) \in A \setminus A^{\prime}(m)$ is 2-A, and it is infinitesimal by Corollary~\ref{cor:2-AI-factor-Denjoy} since the subsets $A^{\prime}(m) \subset A$ generate the same $\Q$-submodule of $\Torus$.
Hence there is a 2-AI factor map $(\Denjoy_{\alpha, A}, F_{\alpha, A}^{m}) \to (\Denjoy_{m\alpha, A^{\prime}(m)}, F_{m\alpha, A^{\prime}(m)})$.

The last part is the special case where $A$ is the $R_{\alpha}$-orbit of $0$. % A=\{R_\alpha^n(0)\colon n\in \Z\}$.
\end{proof}

\begin{Remark}[iterated trees and fibered products of 2-A factors]
There exist 2-A factor maps between zero-dimensional systems with non-trivial fibers lying over points belonging to non-trivial asymptotic classes. We record here an example of such. 

Fix an irrational $\alpha \in \R$ and let $(X_{\alpha},\Shift_{\alpha})$ be the associated Sturmian subshift with canonical 2-A factor map $\pi \colon X_{\alpha} \to \Torus$. Let $(Z,U)$ denote the fiber product of $\pi$ with itself, that is $Z=\{(x_1,x_2) \in X_{\alpha} \times X_{\alpha} \colon \pi(x_1) = \pi(x_2)\}\subset X_\alpha\times X_\alpha$ with $U= \Shift \times \Shift$, so that we have a diagram:
\begin{equation}
\xymatrix{
& Z \ar[dr]^{\psi_{2}} \ar[dl]_{\psi_{1}} & \\
X_{\alpha} \ar[dr]^{\pi} & & X_{\alpha} \ar[dl]_{\pi} \\
& \Torus & \\
}
\end{equation}
We claim that $\psi_{1}$ and $\psi_{2}$ are 2-A factor maps. We argue for $\psi_{1}$. Fix $x_1 \in X_\alpha$. By construction, $\psi_{1}^{-1}(x_1)= \{(x_1,x_2)\in X_\alpha\times X_\alpha \colon x_2 \in \pi^{-1}(\pi(x_1))\}$.
Thus $\Card{\psi_{1}^{-1}(x_1)} = \Card{\pi^{-1}(\pi(x_1))}\in \{1,2\}$ is equal to $2$ precisely when $x_1$ belongs to a non-trivial $\Shift$-asymptotic class.
Finally, suppose that $\{(x_{1},x_2), (x_{1}, x_2')\}$ is a non-trivial fiber of $\psi_{1}$. Then since $\pi(x_{2})=\pi(x_{2}')$ we must have $x_{2}\Bumpeq x_{2}'$ for $\Shift$, hence $(x_{1},x_2)\Bumpeq (x_{1}, x_2')$ for $U$.

Note that while $Z$ is zero-dimensional, the system $(Z,U)$ is not minimal since the diagonal $\{(x,x) \in X_{\alpha} \times X_{\alpha}\}$ is proper, closed, and invariant. 
The off diagonal points are isolated, and they form two orbits, which are both asymptotic to the orbits of the nontrivial asymptotic pair in $X_{\alpha}$.
Hence the factor $Z\to \Torus$ has non-trivial fibers of cardinal $4$ consisting of asymptotic points.
% Note that quotienting $Z$ by the involution exchanging factors yields a similar factor map with has non-trivial fibers of cardinal $3$ consisting of asymptotic points.
\end{Remark}

\subsection{2-A factors of interval exchanges and minimal models}
\label{subsec:IES-2A-minimal-model}
% The results in this subsection may be compared with \cite[Section 2.1]{Dahmani_Groups-IET_2019} about minimal models of interval exchange transformations.

Fix our ordered alphabet $\Al=\{1,\dots,d\}$ on $d\in \N_{\ge 1}$ letters, and recall the previous subsections about interval exchange systems $(X,\Shift)\subset (\Al, \Shift)$.
In particular, if a minimal subshift $(X,\Shift)\subset (\Al,\Shift)$ has complexity $C_X(n)=(d-1)n+1$, then the IETs with that encoding all have ido, the same irreducible permutation $\pi\colon \Al\to \Al$, and their lengths $\lambda\colon \Al \to \R_{>0}$ correspond to the invariant measures of $(X, \Shift)$.

Let us first describe the 2-A factors of $(X,\Shift)\subset (\Al,\Shift)$, a minimal IES with complexity $C_X(n)=(d-1)n+1$ and permutation $\pi\colon \Al\to \Al$.

The proofs of the following two propositions are left as exercises. 
They follow from the description in \cite[Section 2.1]{Dahmani_Groups-IET_2019} of a \emph{minimal model} for an IET, combined with the characterisation \cite{Ferenczi-Zamboni_Languages-k-IET_2008} of languages associated to $d$-IETs with ido (relying on the notion of bispecial words).

\begin{Proposition}[cyclic case]
    The system $(X,\Shift)$ factors onto a Sturmian system $(X_\rho, \Shift_\rho)$ if and only if the permutation $\pi$ preserves the cyclic order of $\Al\bmod{d}$.
    %, meaning that $\pi(k)\equiv \pi(k)+r\bmod{d}$ for some $r\in \Z/d\setminus\{0\}$, 
    % $\pi(1),\dots,\pi(d)$ is a cyclic reordering of $(1,\dots,d)$
    If so, then $(X,\Shift)$ is uniquely ergodic and its invariant probability measure $\lambda\in \R_{>0}^\Al$ yields an IET $T$ with discontinuities $\alpha_k=\sum_1^k \lambda_j$ such that $\rho=T(\alpha_d)=\alpha_{\pi(d)}$. Moreover $(X,\Shift)$ is topologically conjugate to the Denjoy system with parameters $\rho$ and $Q\equiv\{\alpha_1,\dots,\alpha_d\}\bmod R_\rho$, so its 2-A factors have already been described. Finally, if two Sturmian systems on $\Al=\{1,2\}$ are topologically conjugate, then they are equal up to the morphism extending the involution $1\leftrightarrow 2$.
\end{Proposition}

Now assume that $\pi$ is not a cyclic permutation preserving the cyclic order.
We say that $\pi\colon \Al \to \Al$ can be \emph{factored} at $k\in \Al\setminus\{d\}$ when $\pi(k+1)=\pi(k)+1\bmod{d}$, that is when $T$ has $\alpha_k$ as a \emph{removable discontinuity}.

\begin{Proposition}[acyclic case]
    Suppose that $\pi$ does not preserve the cyclic order of $\Al\bmod{d}$, namely that $(X,\Shift)$ does not factor onto a Sturmian system. The 2-A factors of $(X,\Shift)$ are obtained by merging disjoint pairs of consecutive letters $\{k,k+1\}$ at which $\pi$ factors.
    The factored systems are minimal IES with ido on a smaller alphabet, and there is a unique \emph{smallest 2-A factor}. Finally consider two minimal IES on $\Al = \{1,\dots, d\}$ which are not factorable: if they are topologically conjugate then they are equal up to the morphism extending global inversion $k\mapsto d+1-k$.
\end{Proposition}

% \begin{proof}
%The proofs are left as exercises. 
% They follow from the description in \cite[Section 2.1]{Dahmani_Groups-IET_2019} of a \emph{minimal model} for IET, combined with characterisation \cite{Ferenczi-Zamboni_Languages-k-IET_2008} of languages associated to $d$-IETs with ido (relying on the notion of bispecial words).
% \end{proof}

To summarize, any minimal IES $(X,\Shift)$ with ido has a \emph{smallest 2-A factor} $(X_{\downarrow}, \Shift_{\downarrow})$ that is an IES with ido, which is uniquely determined as a subshift on an ordered alphabet (up to order reversal of the alphabet).
Let us note that if an IET has encoding $X$, then its minimal model will have encoding $(X_{\downarrow}, \sigma_{\downarrow})$.

\section{Isogenies: virtual 2-AI and flow equivalences}
% Powers of systems: virtual equivalences and isogenies

We now bring together our various equivalence relations with the additional move of passing to powers. Combining these will lead to our notion of isogeny. 

\begin{Definition}[virtual equivalence]
\label{def:virtual-equivalence}
Given an equivalence relation $\mathcal{R}$ on systems, we say that two systems $(X,S)$ and $(Y,T)$ are \emph{virtually $\mathcal{R}$} when there exist $m,n\in \N_{>0}$ such that $(X,S^m)$ and $(Y,T^n)$ are $\mathcal{R}$.
\end{Definition}

Virtual flow equivalence has an interpretation in terms of covers, in the following sense. Given a system $(X,S)$, the product space $X \times \R$ is endowed with the diagonal action of the lattice $\Z\subset \R$ by $n\cdot(x,t)=(S^{n}(x),t-n)$.
For $m\in \Z$, the quotient space $\Sigma_S^{(m)} X := X\times \R \bmod{m\Z}$ by the normal subgroup $m\Z\subset \Z$ is a $\Z/m$-Galois cover of $ \Sigma_S X$ with covering map denoted
\begin{equation}\label{eqn:rhocoveringmap}
\rho^{(m)} \colon \Sigma_{S}^{(m)}X \to \Sigma_{S}X.
\end{equation}
Note that when $m \ge 1$ there is a flow equivalence
\begin{equation*}
\begin{gathered}
\chi^{(m)} \colon \Sigma_{S^{m}}X \to \Sigma^{(m)}_{S}X\\
\chi^{(m)} \colon (x,t) \mapsto (x,m t).
\end{gathered}
\end{equation*}

Hence two systems $(X,S)$ and $(Y,T)$ are virtually flow equivalent if and only if there exists an orientation preserving homeomorphism between the $m,n$-fold covers
\begin{equation*}
    \Sigma^{(m)}_{S}X \stackrel{\cong}{\longrightarrow} \Sigma^{(n)}_{T}Y.
\end{equation*}

\begin{comment}
Suppose that $(X,S)$ is a totally uniquely ergodic zero-dimensional system and let $\mu$ denote the unique $S$-invariant Borel probability measure.
The flow on $\Sigma_{S^{m}}X$ is uniquely ergodic, with the invariant probability measure defined locally by $\nu = \mu \times \textrm{Leb}$. Moreover, $\Sigma_{S}^{(m)}X$ is uniquely ergodic, with invariant probability measure defined locally by $\nu_{m} = \mu \times \frac{1}{m}\textrm{Leb}$.
The flow equivalences $\chi^{(m)} \colon \Sigma_{S^{m}}X \to \Sigma^{(m)}_{S}X$ satisfy $\chi^{(m)}_{*}(\nu) = \nu_{m}$.
Let $\rho^{(m)} \colon \Sigma_{S}^{(m)}X \to \Sigma_{S}X$ denote the $m$-fold covering map defined in~\eqref{eqn:rhocoveringmap}.
The composition $\rho^{(m)} \circ \chi^{(m)} \colon \Sigma_{S^{m}}X \to \Sigma_{S}X$ is a covering map.
\end{comment}

We now come to one of the key definitions for us.
%\begin{Definition}[isogeny]
%We say that two zero-dimensional systems $(X,S)$ and $(Y,T)$ are \emph{isogenous} when they are connected by a sequence of virtual 2-AI equivalences and flow equivalences.
%\end{Definition}
\begin{Definition}[isogeny]
\label{def:isogeny}
We say that two zero-dimensional systems $(X,S)$ and $(Y,T)$ are \emph{isogenous} when there exist systems $(X_{i},S_{i})$ for $1 \le i \le n$ such that $(X_{1},S_{1}) = (X,S)$, $(X_{n},S_{n}) = (Y,T)$, and for every $1 \le i \le n-1$ the systems $(X_{i},S_{i})$ and $(X_{i+1},S_{i+1})$ are either virtually 2-AI equivalent or flow equivalent.
\end{Definition}
Equivalently, two systems are isogenous when they are connected by a sequence of any of the following: virtual topological conjugacies, 2-AI equivalences and flow equivalences.

Certain properties, such as minimality or unique ergodicity are not necessarily preserved by taking powers. Following usual terminology, given a property $\mathtt{P}$ for systems, we say a system $(X,S)$ is totally $\mathtt{P}$ if all powers $(X,S^m)$ for $m\in \N_{>0}$ have property $\mathtt{P}$.

\begin{Remark}[isogenies preserve total ergodicity]
Given two isogenous systems, if one is totally uniquely ergodic then so is the other. 
This follows from the fact that both flow equivalence and 2-AI factors preserve unique ergodicity.
\end{Remark}

\subsection{Isogenies of Cantor systems}

\begin{Lemma}[virtual 2-AI and flow of minimal Cantor systems preserve ratios of states]
\label{lem:Virtual-2-AI-flow_states}
Consider minimal Cantor systems $(X,S)$ and $(Y,T)$ with invariant probability measures $\mu\in \mathcal{M}_{p}(X,S)$ and $\nu\in \mathcal{M}_{p}(Y,T)$. If there exists one of the following:
\begin{enumerate}[noitemsep]
\item a measure preserving virtual conjugacy $h \colon (X,S^{m},\mu) \to (Y,T^{n},\nu)$
\item a measure preserving  2-AI factor $\pi \colon (X,S,\mu) \to (Y,T,\nu)$
\item a flow equivalence $\Psi \colon \Sigma_{S}X \to \Sigma_{T}Y$ and $\lambda>0$ such that $\Psi_{*}(\mu) = \lambda \nu$
\end{enumerate}
then there exists a real $\lambda > 0$ such that the states $\tau_{\mu}$ and $\tau_{\nu}$ satisfy 
$\im(\tau_{\mu}) \otimes \Q = \lambda \cdot \im(\tau_\nu) \otimes \Q$.
\end{Lemma}
\begin{proof}
The first item is straightforward: the image of $\tau_{\mu} \colon \CoInv_{S} \to \R$ is equal to the image of $\tau_{\mu} \colon \CoInv_{S^{m}} \to \R$, and a topological conjugacy $h \colon (X,S^{m}) \to (Y,T^{n})$ which takes $\mu$ to $\nu$ satisfies $\im(\tau_{\mu}) = \im(\tau_{\nu})$. 

The second item follows from Proposition~\ref{prop:2A_infinitesimal-implies-0-measured}.

For the third item, the flow equivalence $\Psi$ induces an isomorphism $\Psi^{*} \colon \CoInv_{T} \to \CoInv_{S}$. Given $[f] \in \CoInv_{S}$, there exists $[g] \in \CoInv_{T}$ such that $\Psi^{*}([g]) = [f]$, and hence
\begin{equation*}
    \tau_{\mu}([f]) = \tau_{\mu}(\Psi^{*}([g])) = \Psi_{*}(\tau_{\mu})([g])=\lambda\tau_{\nu}([g]).
\end{equation*}
It follows that $\im(\tau_{\mu}) \otimes \Q \subset \lambda \cdot \im(\tau_{\nu}) \otimes \Q $. The same argument applied to $\Psi^{-1}$ implies $\im(\tau_{\nu}) \otimes \Q \subset \lambda^{-1} \cdot \im(\tau_{\mu}) \otimes \Q $, so $\lambda \cdot \im(\tau_{\nu}) \otimes \Q \subset \im(\tau_{\mu}) \otimes \Q $.
\end{proof}

\begin{Corollary}[when isogenies preserve ratios of states]
\label{cor:isogenous-Qsates}
Suppose $(X,S)$ and $(Y,T)$ are Cantor systems which are totally minimal and totally uniquely ergodic.
If they are isogenous, then there exists a positive real $\lambda>0$ such that their unique invariant probability measures $\mu$ and $\nu$ satisfy:  
\begin{equation*}
    \im(\tau_\mu)\otimes \Q=\lambda \cdot (\im(\tau_\nu) \otimes \Q)
\end{equation*}
\end{Corollary}
\begin{proof}
Suppose $(X,S)$ and $(Y,T)$ are both totally minimal and totally uniquely ergodic, and let $\mu, \nu$ be the unique invariant measures.
Both 2-A factors and flow equivalences between these minimal Cantor systems induce isomorphisms between their spaces of invariant probability measures, so they take $\mu$ to $\nu$.
Moreover, a virtual conjugacy $h \colon (X,S^{m}) \to (Y,T^{n})$ must take the unique invariant measure $\mu$ under $S^m$ to the unique invariant measure $\nu$ under $T^n$.
Hence the result follows from Lemma~\ref{lem:Virtual-2-AI-flow_states}.
\end{proof}

\subsection{Isogenies of Sturmian systems and Denjoy systems}

\subsubsection*{Isogenies of Sturmian systems}

Recall that Fokkink's Theorem~\ref{thm:flow-equivalence-Sturmian-systems} classifies Sturmian systems up to flow equivalence by the $\PGL_{2}(\Z)$-orbit of their rotation number.

We are now in a position to prove one of our main theorems, classifying Sturmian systems up to isogenies by the $\PGL_2(\Q)$-orbit of their rotation number.

\begin{Theorem}[isogenies of Sturmian systems]
\label{thm:isogenous-Sturmian}
% For $\alpha,\beta\in \R \setminus \Q$. 
The following are equivalent:
\begin{itemize}[noitemsep]
    \item[-] The points $\alpha,\beta\in \R\Proj^1\setminus \Q\Proj^1$ belong to the same $\PGL_{2}(\Q)$-orbit.
    \item[-] The Sturmian subshifts $(X_{\alpha},\Shift_{\alpha})$, $(X_{\beta},\Shift_{\beta})$ are isogenous.
\end{itemize}
\end{Theorem}

\begin{proof}
Suppose that $\alpha,\beta\in \R\setminus \Q$ are $\PGL_{2}(\Q)$-equivalent.
Let $M \in \PGL_{2}(\Q)$ be such that $M \alpha = \beta$ and choose a representative $M\in \Mat_{2}(\Z)$. By Lemma~\ref{lem:Smith-PSL2} we may factorize $M = UDV$ with $U,V \in \SL_{2}(\Z)$ and $D = \operatorname{diag}(m,1)$ for some $m\in \N_{>0}$, so that $U^{-1}\beta = DV\alpha$.
Note that $(X_{\alpha},\Shift_{\alpha})$ is flow equivalent to $(X_{V\alpha},\Shift_{V \alpha})$ and $(X_{\beta},\Shift_{\beta})$ is flow equivalent to $(X_{U^{-1}\beta},\Shift_{U^{-1}\beta})$ by Theorem~\ref{thm:flow-equivalence-Sturmian-systems}. Since $DV\alpha = m \cdot V\alpha$, the systems $(X_{V\alpha},\Shift_{V\alpha})$ and $(X_{DV\alpha},\Shift_{DV\alpha})$ are virtually 2-AI equivalent by Corollary~\ref{cor:2-AI-factor-Denjoy-power-m}, from which one implication follows.

Now suppose that $(X_{\alpha},\Shift_{\alpha})$ and $(X_{\beta},\Shift_{\beta})$ are isogenous.
These systems are totally minimal and totally uniquely ergodic and by Theorem~\ref{thm:Denjoy-coinvariants-state} their states satisfy
$\im(\tau_X)=\Z+\alpha\Z$ and $\im(\tau_Y)=\Z+\beta\Z$ so by Corollary~\ref{cor:isogenous-Qsates} there exists $\lambda>0$ such that we have the equality of oriented lattices $\Q+\alpha\Q=\lambda \Q+\lambda\beta\Q$.
We thus conclude using Lemma~\ref{lem:proportional-rank2-lattices} that $\alpha$ and $\beta$ belong to the same $\PGL_2(\Q)$ orbit.
\end{proof}

\subsubsection*{Isogenies of Denjoy systems}

Recall the classification of Denjoy systems up to flow equivalence from Theorem~\ref{thm:flow-equivalence-Denjoy-systems}.
We propose the following conjecture, classifying Denjoy systems up to isogeny.
This generalizes the case of Sturmian systems.

\begin{Conjecture}[isogenies of Denjoy systems]
\label{conj:isogenies-Denjoy-systems}
% Consider irrationals $\alpha, \beta \in \Torus$, and $R_\alpha, R_\beta$-invariant subsets $A, B\subset \Torus$ containing the origin.
% The Denjoy systems $(\Denjoy_{\alpha, A}, F_{\alpha, A})$ and $(\Denjoy_{\beta, B}, F_{\beta, B})$ are flow equivalent if and only if there exists $M=$ such that:
For $k\in \{0,1\}$, the Denjoy systems $(\Denjoy_{F_k}, F_k)$ with invariants $(\alpha_k, A_k)$ are isogenous if and only if there exists 
$M= (m_{ij}) \in \PGL_2(\Q)$ such that:
\begin{equation*}
    \alpha_1 = \frac{m_{11}\cdot\alpha_0+m_{12}}{m_{21}\cdot\alpha_0+m_{22}} 
    \qquad \mathrm{and} \qquad 
    \Q{A_1} \equiv \frac{1}{m_{21}\cdot\alpha_0+m_{22}} \cdot \Q{A_0}
\end{equation*}
% where $\tilde{A}\subset \R$ is the lift of $A\subset \Torus$ and $\equiv$ denotes equality of subsets of $\R$ up to translation.
where $\Q{A_k}\subset \Torus$ denotes the $\Q$-submodule of $\Torus$ generated by a subset $A_k\subset \Torus$, and $\equiv$ denotes equality of sets in $\Torus$ up to rotations.
\end{Conjecture}

For now, we can prove the following weak version of the conjecture.

\begin{Theorem}[isogenies of Denjoy systems]
\label{thm:isogenous-Denjoy}
For $k \in \{0,1\}$, consider Denjoy systems $(\Denjoy_{F_k}, F_k)$ with invariants $(\alpha_k, A_k)$.
\begin{enumerate}
\item
If there exists 
$M= (m_{ij}) \in \PGL_2(\Q)$ such that:
\begin{equation*}
    \alpha_1 = \frac{m_{11}\cdot\alpha_0+m_{12}}{m_{21}\cdot\alpha_0+m_{22}} 
    \qquad \mathrm{and} \qquad 
   \Q A_1 \equiv \frac{1}{m_{21}\cdot\alpha_0+m_{22}} \cdot \Q A_0
\end{equation*}
then the Denjoy systems $(\Denjoy_{F_k}, F_k)$ are isogenous.
\item
If the Denjoy systems are isogenous through a chain of Denjoy systems, then there exists 
$M= (m_{ij}) \in \PGL_2(\Q)$ such that:
\begin{equation*}
    \alpha_1 = \frac{m_{11}\cdot\alpha_0+m_{12}}{m_{21}\cdot\alpha_0+m_{22}} 
    \qquad \mathrm{and} \qquad 
    \Q A_1 \equiv \frac{1}{m_{21}\cdot\alpha_0+m_{22}} \cdot \Q A_0.
\end{equation*}
\end{enumerate}
\end{Theorem}

\begin{proof}
For part one, consider $M \in \PGL_{2}(\Q)$ as in the assumption, and choose a lift still denoted $M\in \Mat_2(\Z)$ (thus $\det(M)\ne 0$). 
By Lemma~\ref{lem:Smith-PSL2} we may factorize $M = U_{1}DU_{0}$ with $U_{0},U_{1} \in \SL_{2}(\Z)$ and $D = \operatorname{diag}(m,1)$ for some nonzero $m \in \Z$. Thus $U_{1}^{-1}\alpha_{1} = DU_{0}\alpha_{0} = m\cdot U_{0} \alpha_{0}$. 
Let us introduce the following notations:
\begin{equation*}
    U_{1}^{-1} = \begin{psmallmatrix} a & b \\ c & d \end{psmallmatrix} 
    \quad\mathrm{and}\quad 
    U_{0} = \begin{psmallmatrix} e & f \\ g & h \end{psmallmatrix}
\end{equation*}
\begin{equation*}
    \alpha_{1}' = U_{1}^{-1}\alpha_{1}
    \quad \mathrm{and} \quad
    \alpha_{0}' = U_{0}\alpha_{0}
\end{equation*}
\begin{equation*}
    A_{1}' = \tfrac{1}{c \cdot \alpha_{1} + d} \cdot A_{1}
    \quad \mathrm{and}
    \quad
    A_{0}' = \tfrac{1}{g \cdot \alpha_{0} + h} \cdot A_{0}.
\end{equation*}

By Theorem~\ref{thm:flow-equivalence-Denjoy-systems}, the Denjoy system $(\Denjoy_{F_{k}},F_{k})$ is flow equivalent to the Denjoy system with invariants $(\alpha'_{k},A'_{k})$, and we will now show that these two systems are virtually 2-AI equivalent.

By Corollary~\ref{cor:2-AI-factor-Denjoy-power-m}, the Denjoy system $(\Denjoy_{\alpha'_0, A'_0}, F^{m}_{\alpha_0, A'_0})$ is 2-AI equivalent to the Denjoy system with invariants $(m\cdot \alpha'_{0}, A'_{0}) = (\alpha'_{1}, A'_{0})$. Since $U_{1}^{-1}M = DU_{0}$, we have $cm_{11} + dm_{21} = g$ and $cm_{12}+dm_{22} = h$. 
By assumption $\Q A_{1} = \tfrac{1}{m_{21} \cdot \alpha_{0} + m_{22}} \cdot \Q A_{0}$ and $\alpha_{1}(m_{21} \cdot \alpha_{0} + m_{22}) = m_{11} \cdot \alpha_{0} + m_{12}$ so:
\begin{align*}
    \Q A'_{1} &= \tfrac{1}{c \cdot \alpha_{1} + d} \cdot \Q A_{1} = \tfrac{1}{c \cdot \alpha_{1} + d} \cdot \tfrac{1}{m_{21} \cdot \alpha_{0} + m_{22}} \cdot \Q A_{0} \\
    &= \tfrac{1}{(cm_{11}+dm_{21})\alpha_{0} + (c m_{12} + dm_{22})} \cdot \Q A_{0} = \tfrac{1}{g \cdot \alpha_{0} + h} \cdot \Q A_{0} = \Q A'_{0}.
\end{align*}
Thus $\Q A'_{1} = \Q A'_{0}$, so by Corollary~\ref{cor:2-AI-equiv-Denjoy}, the Denjoy systems with invariants $(\alpha'_{1},A'_{0})$ and $(\alpha'_{1},A'_{1})$ are 2-AI equivalent, completing the proof of part one.

For part two, suppose the two Denjoy systems $(\Denjoy_{F_k}, F_k)$ are isogenous through a chain of Denjoy systems.
Recall that any isogeny can be written as a chain of flow equivalences, taking powers or roots, and 2-AI equivalences. 
Any 2-AI factor between Denjoy systems does not change the rotation number, and by Corollary~\ref{cor:2-AI-factor-Denjoy}, it does not change the rational span of the cut points.
Moreover, if $(\Denjoy_{\alpha, A}, F_{\alpha, A})$ is a Denjoy system and $m \in \N$, then $(\Denjoy_{\alpha, A}, F_{\alpha, A}^{m})$ is a Denjoy system with rotation number $m\alpha$ and same cut point set. 
The result then follows from Theorem~\ref{thm:flow-equivalence-Denjoy-systems}.
\end{proof}

\subsection{Isogenies of interval exchange systems}

For an IET $T\colon (\alpha_0,\alpha_d]\to (\alpha_0,\alpha_d]$ with discontinuities $\alpha_k \in (\alpha_0,\alpha_d]$ and length vector $\lambda \in (\R_{>0})^\Al$, its \textsc{SAF}-invariant is the element $SAF(T)\in \R\wedge_\Q \R$ defined by
\begin{equation*}
    SAF(T)=\int_0^1 1\otimes (T(x)-x) = \sum_{k\in \Al} \lambda_k \otimes (T(\alpha_k)-\alpha_k).
\end{equation*}

Arnoux-Fathi showed \cite[Proposition 2.3]{Veech_metric-theory-IET-2-ApproxByPrimitive_1984} that if an interval $I\subset (\alpha_0, \alpha_d]$ intersects every $T$-orbit then the IET induced by $T$ on $I$ has the same \textsc{SAF}-invariant.

The IETs on $(0,1]$ form a group under composition, and Sah showed \cite[Theorem 1.3]{Veech_metric-theory-IET-3-SAF_1984} that its abelianisation is given by the surjective homomorphism $\operatorname{SAF}\colon \operatorname{IET}\to \R\wedge_\Q \R$.

Moreover, it is straightforward from the definition that the \textsc{SAF}-invariant is preserved by 2-A factors, in particular any IET has the same \textsc{SAF}-invariant as its minimal model.

% It follows from these results that if two IET with ido are related by isogenies among IET (preserving the measure), then their \textsc{SAF}-invariants are rational multiples of one another.

\begin{Definition}[rational invariants]
\label{def:rational-invariants-IES}
Consider an IES $(X,\Shift)$ with ido and let $\mu_1,\dots,\mu_e$ be its $e\in \N_{>0}$ ergodic probability measures. Each $\mu_i$ corresponds to an IET $T_{\mu_i}\colon (0,1]\mapsto (0,1]$ with that encoding, to which we associate the $\Q$-vector space $\Q\otimes \im(\tau_{\mu_i})\subset \R$ of dimension at most $d-1$ generated by the $\mu_i$-measures of cylinders in $X$, and the vector $\operatorname{SAF}(T_{\mu_i})\in \R\wedge_\Q \R$.
    
We define the \emph{rational invariants} of $(X,\Shift)$ 
as the product of subspaces $\prod_i \Q\otimes \im(\tau_{\mu_i})\subset \R^e$ up to multiplication by $\lambda \in \R_{>0}$ on $\R^e$ together with the product of vectors $\prod_i \operatorname{SAF}(T_{\mu_i})\in (\R\wedge_\Q \R)^e$ up to multiplication by $m/n\in \Q_{>0}$.
\end{Definition}
\begin{Conjecture}[rational invariants of isogenies]
    \label{conj:isogenies-IES=Rational-Invariants}
    If two IES with ido are isogenous then they have the same rational invariants.
\end{Conjecture}

Let us prove the conjecture in the totally ergodic case for isogenies remaining among IES.

\begin{Proposition}[isogenous totally ergodic IES have same rational invariants]
    \label{prop:isogenies-total-ergodic-IES}
    Consider two IES with ido which are totally uniquely ergodic.
    If they are related by isogenies through IES, then they have the same rational invariants.
\end{Proposition}
\begin{proof}
    Note that an IES with ido has its powers with ido, hence it is totally minimal. The invariance of the $\R$-projectivized rational tracial image subgroup $\im(\tau)\otimes \Q \subset \R$ follows from Corollary~\ref{cor:isogenous-Qsates} (the $\R$-projectivization is needed because of flow equivalences). The invariance of the $\Q$-projectivized vector $\operatorname{SAF}(T_\mu)$ follows from the properties of the $\operatorname{SAF}$-invariant (the $\Q$-projectivization is required because of virtual conjugacies).
\end{proof}

\begin{comment}[most IES]
\label{rem:most-ies-ido-LR}
This result applies most IES in the following sense.
For every irreducible permutation of $\Al=\{1,\dots,d\}$ and for almost every $\lambda \in \R_{>0}^d$, the corresponding IET has ido (by \cite{Keane_IET_1975}) 
% it suffices to have linearly independent lengths over the rationals 
and is totally uniquely ergodic (by \cite[Theorem 1.7]{Veech_metric-theory-IET-1-Spectral_1984}).
\end{comment}

This result applies to the IES with ido that are linearly recurrent (see Subsection~\ref{subsec:subshifts}).

\begin{Corollary}[linearly recurrent IES with ido]
\label{cor:isogenies-IES-ido-LR}
    Consider two minimal IES with ido which are linearly recurrent.
    If they are isogenous through IES, then they have the same rational invariants.
\end{Corollary}

\begin{proof}
Note that if an IES has ido, then its powers are also IES with ido, and hence totally minimal (by Keane).
We saw in Subsection~\ref{subsec:subshifts} (using Boshernitzan's criterion for linear recurrence) that if a subshift is totally minimal and linearly recurrent then it is totally linearly recurrent, hence totally ergodic.
Thus IES with ido that are linearly recurrent are totally ergodic.
\end{proof}

\begin{Remark}[Diophantine properties of IES]
The fact that a real number is irrational if and only if its continued fraction expansion is infinite generalizes to the fact that an IES has ido if and only if its Rauzy induction is infinite.

The fact that an irrational is badly approximable by rationals if and only if its continued fraction expansion has bounded entries generalizes to the fact that an IES with ido satisfies Boshernitzan's criterion if and only if its Rauzy induction yields a primitive proper S-adic representation.

The fact that real numbers with infinite periodic continued fraction expansions correspond to slopes of foliations of the torus which are fixed by elements of the mapping class group generalizes to the fact that IES with ido and infinite periodic Rauzy induction correspond to the attractive foliations of pseudo-Anosov mapping class of translation surfaces.
\end{Remark}

More specifically, we wish to mention the special case of minimal IES which are self-induced.
A system $(X,S)$ is called self-induced when it admits a proper return system $(C, S_C)$ such that $(X,S)$ and $(C,S_C)$ are topologically conjugate.
For a minimal subshift, this is equivalent to saying that it is substitutive in the sense of \cite{Durand-Ormes-Petite_self-induced-systems_2018}. 
For IES, this is equivalent to saying that its Rauzy induction is periodic.
Let us recall the well known relation between self-induced IES and pseudo-Anosov mapping classes of orientable surfaces.

Consider a minimal IES $(X,S)$ on $d$ letters, and assume that it is self-induced. 
Under these conditions it must satisfy Keane's ido condition, and it is linearly recurrent, hence uniquely ergodic.
By unique ergodicity it arises from a unique IET $T\colon (0,1]\to (0,1]$ up to conjugation by $x \mapsto 1-x$ on the complement of the discontinuities.
The Rauzy induction will be infinite (ido condition), and periodic (by the self-induced condition).
Hence the corresponding length data $\lambda\in \R_+^d$ is given by the probability right eigenvector of a primitive matrix $M\in \SL_d(\N)$ for its maximal eigenvalue $\lambda_0 \in \R_+$ (so the entries $\lambda_1,\dots,\lambda_d$ lie in the algebraic number field $\Q(\lambda_0)$).
Moreover, there is a preferred suspension data given by the dual left eigenvector of $M$ for that eigenvalue $\lambda_0$: the zippered-rectangle construction \cite[Section 4]{Yoccoz_IET-translation-surfaces_2010} yields a translation surface $\Bar{\Sigma}_S$ which is marked by a horizontal segment, whose vertical and horizontal directions of $\Bar{\Sigma}_S$ are the attractive and repulsive line-fields of a unique primitive pseudo-Anosov mapping class $A\in \Mod^+(\Bar{\Sigma}_S)$.
After puncturing $\Bar{\Sigma}_S$ at its singularities (namely those of the stable or unstable foliations of $A$), we obtain a punctured surface $\Bar{\Sigma}^*_S$ whose negative Euler characteristic is equal to the number $-\chi(\Bar{\Sigma}^*_S) = d-1$ of internal discontinuities of the IET $T$.

\begin{Corollary}[type preserving isogenies of self-induced IES]
    \label{cor:isogenies-self-induced-IES}
    If two minimal self-induced IES are isogenous through IES, then they have the same rational invariants. 
\end{Corollary}

The case of self-induced IES thus provides a terrain to investigate to what extent the rational invariants of IES determine their isogeny class. This motivates the following.

\begin{Question}[rational invariants of pseudo-Anosov mapping classes]
    \label{quest:rational-invariants-pseudo-Anosov}
    Given an irreducible pseudo-Anosov mapping class of an oriented surface, what are the rational invariants of its attracting measured lamination ?
    Under which conditions do the rational invariants of self-induced IES uniquely determine them up to isogenies (among IES)?
\end{Question}

%\begin{Remark}[quadratic IES]
The IETs whose lengths generate a quadratic extension of $\Q$ were studied in \cite{Boshernitzan-Carroll_quadratic-IET_1997}. The main result is that they coincide with IETs such that repeatedly inducing on sub-intervals bounded by discontinuities yields a finite number of IETs (generalizing Lagrange's theorem for the continued fraction expansion of quadratic numbers).
%%\end{Remark}
In particular they satisfy Keane's idoc, and hence are totally minimal.

\begin{Remark}[purely periodic quadratic IES]
\label{purely-perodic-IES}
Boshernitzan and Carroll also characterize those IETs with quadratic lengths that are purely periodic under Rauzy induction (namely that are self-induced) in terms of pseudo-Anosov classes whose dilatation factor is a quadratic irrational.
Note that \cite[Corollary 4.4]{Boshernitzan-Carroll_quadratic-IET_1997} shows that when the foliation is orientable, these are obtained as branched covers of the torus with a foliation of quadratic irrational slope.
\end{Remark}

\begin{Corollary}[quadratic IES]
    \label{cor:quadratic-IES}
    If two minimal quadratic IES are isogenous through IES's, then they have the same rational invariants.
\end{Corollary}
\begin{proof}
    We already noted that a minimal quadratic IES satisfies Keane's idoc, so it is totally minimal.
    We now show that it is totally uniquely ergodic so that we may apply Proposition~\ref{prop:isogenies-total-ergodic-IES}.
    
    By \cite[Theorem 6.1]{Boshernitzan-Carroll_quadratic-IET_1997} a minimal quadratic IES is uniquely ergodic.
    A power of a minimal IES that is quadratic is a minimal IES which is quadratic (its lengths lie in the same quadratic field), hence it is again uniquely ergodic.
\end{proof}

\begin{Conjecture}[Veech group commensurators]
\label{conj:veech-groups}
Consider two isogenous IESs.
% $(X_i,S_i)$ for $i=1,2$, and assume that they are isogenous.
If one is minimal self-induced, then so is the other.
For such a pair of isogenous minimal self-induced IES, the associated translation surfaces with pseudo-Anosov automorphisms should satisfy the following:
% Now assume this is the case, and consider the associated translation surfaces $\Bar{\Sigma}_{A_i}$ with pseudo-Anosov automorphisms $A_i\in \Mod(\Bar{\Sigma}_{S_i})$.
one is Veech if and only if the other is Veech, in which case the Veech groups have isomorphic commensurators in $\PGL_2(\R)$.
\end{Conjecture}

\subsection{Classifying systems up to isogeny: questions and constructions}
Our attempts at proving Conjectures~\ref{conj:isogenies-Denjoy-systems} and~\ref{conj:isogenies-IES=Rational-Invariants} led to the following questions and constructions, which are interesting in their own right.

Given a family of zero-dimensional systems (such as Denjoy or IES for instance), a natural goal would be to find classification results for the isogeny classes within that family (or at least for isogenies through members of the family).

Ideally, one might hope to find a simple complete set of isogeny-invariants.
Such invariants may involve the structure of the ordered abelian group of coinvariants, and of its morphisms to $\R$ given by ergodic states.

\subsubsection*{Factoring Isogenies}

Recall that isogenies involve compositions of flow equivalencies and virtual 2-AI equivalencies, in any order.
The following factorisation result could be of help in classifying systems up to isogeny (at least in restriction to certain classes).

\begin{Question}[simple isogenies]
\label{conj:ordering-isogenies}
    If two systems $(X,S)$ and $(Y,T)$ are isogenous, must there exist $m,n\in \N_{\ge 1}$ such that the systems $(X,S^m)$ and $(Y,T^n)$ are 2-AI equivalent to systems that are flow equivalent; diagrammatically:
    \begin{equation*}
        (X,S) \xrightarrow{m} (X,S^m)
        \xrightarrow{\text{2-AI}} (X',S') \xleftrightarrow{\text{flow}}
        (Y',T') \xleftarrow[]{\text{2-AI}} (Y,T^n)
        \xleftarrow[]{n}
        (Y,T).
    \end{equation*}
\end{Question}

A positive answer to the question would improve the classification results for Denjoy systems~\ref{thm:isogenous-Denjoy} and IES~\ref{prop:isogenies-total-ergodic-IES} by removing the hypothesis that the chain of equivalences in the isogeny remain within the respective class.

\subsubsection*{Minimal models under 2-A factors}

One is naturally led to search for ``minimal models'' of systems under 2-A or 2-AI factors.
Those may be applied to obtain the smallest systems $(X',S')$ and $(Y',T')$ in Question~\ref{conj:ordering-isogenies}.

\begin{Question}[factoring asymptotic orbits]
    \label{quest:factoring-asymptotic-orbits}
    Given a system $(X,S)$ with a pair of asymptotic points $x\bumpeq x'$, when does there exist a 2-A factor $\pi \colon (X,S)\to (Y,T)$ such that $\pi(x)=\pi(x')$? Or simply such that $\pi(\{S^n(x)\colon n\in \Z\})=\{S^n(x')\colon n\in \Z\}$ ?

    % Note that the system $Y$ may not be Hausdorff.
\end{Question}

The following question aims at a general construction encompassing the results in Subsection~\ref{subsec:IES-2A-minimal-model} which provide minimal models for IES.

\begin{Question}[terminal 2-A factors]
    \label{quest:terminal-2-A-factors}
    Let $(X,S)$ be a system, and consider the category of factors $(X,S)\to (Y,T)$ obtained as composition of 2-A factors (respectively 2-AI factors).
    Does it admit a universal terminal element?
    In restriction to which families of systems does this hold?
    When such minimal models exist, how do they relate to maximal equicontinuous factors?
\end{Question}

\section{Virtual and eventual flow equivalence}

For an equivalence relation on systems, we have already introduced the notion of virtual equivalence, which asks that some powers are equivalent. The stronger notion of eventual equivalence is defined as follows.

\begin{Definition}[eventual equivalence]
\label{def:eventual-equivalence}
Given an equivalence relation $\mathcal{R}$ on systems, we say that two systems $(X,S)$ and $(Y,T)$ are \emph{eventually $\mathcal{R}$} when there exist $k\in \N$ such that $(X,S^m)$ and $(Y,T^m)$ are $\mathcal{R}$ for all $m \ge k$.
\end{Definition}

\begin{Remark}[some easy implications]
\label{rem:implications-virtual-eventual-flow}
Since topological conjugacy implies flow equivalence, it follows that eventual topological conjugacy implies eventual flow equivalence and virtual conjugacy implies virtual flow equivalence.
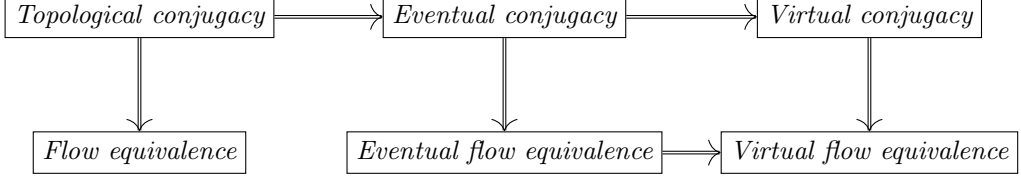
\begin{center}
\begin{tikzpicture}[node distance=4.3cm and 1cm,
box/.style = {draw, align=center}]

    % Define nodes
    \node (A) [box] at (0, 0) {Topological conjugacy};
    \node (B) [box, right of=A, xshift=.5cm] {Eventual conjugacy};
    \node (C) [box, right of=B, xshift=.5cm] {Virtual conjugacy};
    \node (D) [box, below of=A,yshift=2.5cm] {Flow equivalence};
    \node (E) [box, below of=B,yshift=2.5cm] {Eventual flow equivalence};
    \node (F) [box, below of=C,yshift=2.5cm] {Virtual flow equivalence};

    % Draw arrows
    \draw[->, double] (A) -- (B);
    \draw[->, double] (B) -- (C);
    \draw[->, double] (A) -- (D);
    \draw[->, double] (B) -- (E);
    \draw[->, double] (C) -- (F);
    % \draw[->, double] (D) -- (E);
    \draw[->, double] (E) -- (F);
\end{tikzpicture}
\end{center}
\end{Remark}

Eventual topological conjugacy has played an important role, especially in the context of shifts of finite type. Furthermore, for shifts of finite type, the relation between eventual conjugacy and flow equivalence has also been studied in~\cite{Boyle_algebraic-aspects-symbolic-dynamics_2000}, and more recently in~\cite{Boyle_shift-equiv-implies-flow-equiv_2024} where it was shown that for all shifts of finite type, eventual conjugacy implies flow equivalence.

In general, eventual flow equivalence does not imply eventual conjugacy, and flow equivalence does not imply eventual flow equivalence; in fact, both of these implications can fail already for shifts of finite type. We are grateful to Mike Boyle for communicating to us these facts, and showing us the following Theorem (\cite{boylepersonalcommunication}).

\begin{Theorem}[SFT examples] 
% justifying non-implications
\label{thm:boyleexamples}
For nonnegative integral matrix $A\in \Mat_d(\N)$, let $(X_{A},\Shift_{A})$ denote the associated edge shift of finite type. 
\begin{enumerate}
\item
Suppose $C,D \in \Mat_2(\N)$ are primitive with $\det(C) = \det(D) = -1$ but different spectral radii, and set $A=I+C, B=I+D$. Then the shifts of finite type $(X_{A},\Shift_{A})$ and $(X_{B},\Shift_{B})$ are flow equivalent but not eventually flow equivalent.
\item
There exist primitive $E,F\in \Mat_d(\N)$ such that for every $m \ge 1$ the systems $(X_{E},\Shift_{E}^{m})$ and $(X_{F},\Shift_{F}^{m})$ are flow equivalent, but $(X_{E},\Shift_{E})$ and $(X_{F},\Shift_{F})$ are not eventually conjugate.
\end{enumerate}
\end{Theorem}

Explicit examples satisfying item $(2)$ above can be extracted from the proof of their existence, but are cumbersome to write down, and we omit them.

% Since topological conjugacy for Sturmian subshifts is quite rigid, and flow equivalence is more interesting, it is perhaps natural to consider the notion of eventual flow equivalence. We show here that eventual flow equivalence implies topological conjugacy.

\begin{Question}[from eventual flow equivalence to flow equivalence or topological conjugacy]
    In which classes of systems does eventual flow equivalence imply topological conjugacy? In which classes of systems does eventual flow equivalence imply flow equivalence?
\end{Question}

In this Section we show that eventual flow equivalence implies topological conjugacy for Sturmian systems (Theorem~\ref{thm:eventual-flow-Sturmian}) and Denjoy systems (Theorem~\ref{thm:eventual-flow-Denjoy}).
It may be that eventual flow equivalence implies topological conjugacy among large classes of minimal systems, but the proof would require a deeper understanding of topological invariants under flow equivalence and powers, which would have to subsume the difficulties of the self-induced Sturmian case (see proof of Theorem~\ref{thm:eventual-PGL2Z} as well as Example~\ref{ex:p-eventually-PSL2Z-equivalent} and Question~\ref{quest:weak-eventual-PGLZ}).

\subsection{Eventual $\PGL_2(\Z)$-equivalence of real numbers}

The purpose of this subsection is to prove Theorem~\ref{thm:eventual-PGL2Z}.
We first gather arithmetic lemmas (that serve in past, present and future subsections) about rank-$2$ modules and $\PGL_2$-equivalence.

\begin{Lemma}
\label{lem:equal-rank2-lattices}
For $\alpha,\beta \in \R\setminus\Q$, we have $(\Z+\alpha\Z) = (\Z + \beta \Z) \iff \alpha=\pm \beta \bmod{\Z}$.
\end{Lemma}
\begin{proof}
    If $(\Z+\alpha\Z) = (\Z + \beta \Z)$ then there exist $u,v,x,y\in \Z$ such that $\alpha=u+v\beta$ and $\beta=x+y\alpha$. Thus $\alpha=(u+vx)+vy\alpha$ which implies $vy=1$, hence $v=y$ and are equal to $\pm 1$, which proves the implication.
    The converse implication is easy.
\end{proof}

% Recall that an integral domain $\mathcal{O}$ is a ring which has no divisors of $0$: it has a fraction field denoted $\K$.
We will only apply the following Lemma~\ref{lem:proportional-rank2-lattices} with $\K=\R$ or a real quadratic number field, and $\mathcal{O}=\Z$ or an order in such a quadratic number field.
\begin{Lemma}[proportional rank-$2$ lattices]
\label{lem:proportional-rank2-lattices}
Consider a field $\K$ containing a ring $\mathcal{O}$ and $\alpha,\beta \in \K$.
The numbers $\alpha,\beta\in \K$ are $\PGL_{2}(\mathcal{O})$-equivalent if and only if there exists $\lambda \in \K^*$ such that $(\mathcal{O}+\alpha\mathcal{O}) = \lambda(\mathcal{O} + \beta \mathcal{O})$. 
\end{Lemma}

\begin{proof}
% Note that when $\alpha$ or $\beta$ belongs to the fraction field of $\mathcal{O}$, either condition of the equivalence implies that they both belong to that field of fractions and the equivalence is obvious, so we may assume that $(1,\alpha)$ and $(1,\beta)$ both generate modules of rank two over $\mathcal{O}$.

\emph{$\implies$}
Suppose that there exists $M\in \PGL_2(\mathcal{O})$ whose action on $\K\Proj^1$ satisfies $\beta = M\alpha$.
Choose a representative in $M=(m_{ij})\in \GL_{2}(\mathcal{O})$ so that $\beta = \frac{m_{11} \alpha + m_{12}}{m_{21} \alpha + m_{22}}$.
Thus for $x,y\in \mathcal{O}$ we have $x \cdot \beta+y \cdot 1=\left((m_{11}x+m_{21}y)\alpha+(m_{12}x + m_{22} y) \cdot 1\right)\tfrac{1}{m_{21}\alpha + m_{22}}$, and the map $(x,y)\in \mathcal{O}^2\mapsto \left(m_{11}x+m_{21}y,m_{12}x + m_{22} y\right)\in \mathcal{O}^2$ is the linear automorphism corresponding to multiplication by $M$ on the right.
In other words, we have the following equality of $\mathcal{O}$-submodules of rank $\le 2$ of $\K$:
\begin{equation*}
    \mathcal{O}^2 \cdot
    \begin{psmallmatrix} 
    \beta & 0 \\ 
    0 & 1 
    \end{psmallmatrix}
    = \mathcal{O}^2 \cdot
    \begin{psmallmatrix} 
    m_{11} & m_{12} \\ 
    m_{21} & m_{22}
    \end{psmallmatrix}
    \cdot
    \begin{psmallmatrix} 
    \alpha & 0 \\ 
    0 &  1
    \end{psmallmatrix}
    \cdot \tfrac{1}{m_{21}\alpha + m_{22}}
    = \mathcal{O}^2 \cdot
    \begin{psmallmatrix} 
    \alpha & 0 \\ 
    0 &  1
    \end{psmallmatrix}
    \cdot \tfrac{1}{m_{21}\alpha + m_{22}}
\end{equation*}
which proves $(\mathcal{O}+\alpha\mathcal{O})=\lambda(\mathcal{O}+\beta\mathcal{O})$ for $\lambda = 1/(m_{21} \alpha + m_{22})\in \K$.

\emph{$\impliedby$}
Suppose that $\lambda \in \K^*$ satisfies $(\mathcal{O}+\alpha\mathcal{O}) = \lambda(\mathcal{O} + \beta \mathcal{O})$. 
The multiplication by $\lambda$ yields an isomorphism of $\mathcal{O}$-modules $M_\lambda \colon (\mathcal{O}+\beta \mathcal{O}) \to (\mathcal{O} + \alpha \mathcal{O})$, which in the bases $(1,\beta), (1,\alpha)$ yields a matrix $M=(m_{ij})\in \GL_2(\mathcal{O})$ such that $\lambda \cdot \beta = m_{11}\alpha+m_{12}$ and $\lambda \cdot 1 = m_{21}\alpha+m_{22}$, whereby $\beta=M\alpha$ for the action of $M\in \PGL_2(\mathcal{O})$ on $\K\Proj^1$. 
\end{proof}

Let us now recall the Smith normal form.
For a ring $\mathcal{O}$ we denote by $\mathcal{O}^\times$ and $\mathcal{O}^*$ its invertible and non-zero elements respectively.
A subset of $\mathcal{O}$ consists of coprime elements when it generates $\mathcal{O}$ as an $\mathcal{O}$-module.

\begin{Lemma}[Smith normal form]
\label{lem:Smith-PSL2}
    Let $\mathcal{O}$ be a principal ideal domain with  fraction field $\K$. Every $A\in \Mat_2(\mathcal{O})$ with $\det(A)\ne 0$ and globally coprime entries has a factorisation of the form $A=UMV$ with $U,V\in \SL_2(\mathcal{O})$ and $M=\operatorname{diag}(m,1)$ for $m\in \mathcal{O}^*$ whose class $m \bmod{\mathcal{O}^\times}$ is uniquely determined by $A$. Hence the automorphism group $\PGL_2(\K)$ of the projective line $\K\Proj^1$ is generated by the union of its:
    \begin{itemize}[noitemsep, align=left]
        \item[-] subgroup of integral units $\PSL_2(\mathcal{O})$
        \item[-] submonoid of integral multiplications $\mathcal{O}^* \simeq \{x\mapsto mx\colon m\in \mathcal{O}^*\}$.
    \end{itemize}
\end{Lemma}

\begin{proof}
    Every element in $\PGL_2(\K)$ has a representative $A\in \Mat_2(\mathcal{O})$ with $\det(A)\ne 0$ and globally coprime entries, so the second statement will follow from the first.
    
    The Smith normal form (as stated in \cite{Stanley_Smith-Normal-Form-combinatorics_2016}) provides a factorisation $A=PDQ$ for some $P,Q\in \GL_2(\mathcal{O})$ and $D=\operatorname{diag}(1,d)$ where $d \in \mathcal{O}^*$ is unique up to multiplication by units $\mathcal{O}^\times$. % Observe that $\det(A)=d$.
    
    Observe that $J=\begin{psmallmatrix}
        0&1\\1&0
    \end{psmallmatrix}$ conjugates $D=\operatorname{diag}(1,d)$ to $D'=\operatorname{diag}(d,1)$, so $A=P'D'Q'$ with $P'=PJ, Q'=JQ \in \GL_2(\mathcal{O})$.
    We deduce the desired factorisation by setting $U=P'\operatorname{diag}(\det(P'),1)$, $V=\operatorname{diag}(\det(Q'),1) Q'$ and $D=\operatorname{diag}(d/\det(P'Q'),1)$.
\end{proof}

The following two remarks recall some background on quadratic orders and their units that will serve subsequently.

\begin{Remark}[Orders in quadratic fields and Dirichlet unit theorem]
    \label{rem:quadratic-orders}
    Consider a real quadratic field $\K \subset \R$.
    An order is a subring $\mathcal{O}\subset \K$ which is also a $\Z$-module of rank two, and its discriminant $\operatorname{disc}(\mathcal{O})\in \N_{>0}$ is defined as the determinant of a $\Z$-basis.
    The ring of integers $\mathcal{O}_\K$ contains every other order, and its discriminant is denoted $\Delta_\K$. 
    (These so called fundamental quadratic discriminants $\Delta_\K$ can be characterised as the positive integers that are either square-free and $\equiv 1\bmod{4}$ or else $4$ times a square-free number that is $\ne 1\bmod{4}$, see~\cite{Cohn_algebraic-numbers-class-fields_1978, Hatcher_topology-numbers_2022}.)
    Thus, every order $\mathcal{O} \subset \K$ is of the form $\Z+f\mathcal{O}_{\K}$ where $f= [\mathcal{O}_\K \colon \mathcal{O}]\in \N_{>0}$ is called its conductor, so it has discriminant $\Delta=f^2\Delta_\K$.
    Conversely, if $\Delta\in \N$ is such that $\Q(\sqrt{\Delta})=\K$, then there is a unique order $\mathcal{O}_\Delta\subset \mathcal{O}_\K$ of discriminant $\Delta$.
\end{Remark}

\begin{Remark}[Dirichlet unit theorem]
    \label{rem:Dirichlet-unit}
    %Let $\Delta\in \N$ such that $\sqrt{\Delta}\notin \N$.
    Consider a real quadratic field $\K \subset \R$, and denote by $\gamma \in \K \mapsto \Bar{\gamma}\in \K$ the Galois involution so that $\gamma\Bar{\gamma}$ is the norm of $\gamma$.
    By Dirichlet's unit theorem \cite{Cohn_algebraic-numbers-class-fields_1978, Hatcher_topology-numbers_2022}, its ring of integers $\mathcal{O}$ has a group of units $\mathcal{O}_\K^{\times}=\{\gamma\in \mathcal{O}\colon \gamma \Bar{\gamma}=\pm 1\}$.
    
    The positive units admit a unique minimal element $\varepsilon_\K=\min \mathcal{O}_\K^{\times} \cap \R_{>0}$, called the fundamental unit of $\mathcal{O}_\K$. 
    That yields the presentation $\mathcal{O}_\K^{\times} = \{\pm \varepsilon^e \colon e\in \Z\}$.

    For $f\in \N_{>0}$, let $\Delta=f^2\Delta_\K$ and consider the order $\mathcal{O}_\Delta\subset \mathcal{O}_\K$: its units form a subgroup of index $f$ in $\mathcal{O}_\K^\pm$, namely $\mathcal{O}_\Delta^\pm = \{\pm \varepsilon^{fe}\colon e\in \Z\}$.

    Denoting $\varepsilon_\Delta=t+s\tfrac{1}{2} \sqrt{\Delta}$ for $(2t,s)\in \N\times \Z$ yields the fundamental solution of the Pell-Fermat equation $(2t)^2-s^2\Delta=\pm 4$ (where again, the equation with $+4$ always has a solution whereas that with $-4$ may or may not have one).
\end{Remark}

Now we relate the action of $\PGL_2(\Q)$ on $\R\Proj^1$ and real quadratic numbers.

\begin{Lemma}[stabilisers in $\PGL_2(\Q)$]
\label{lem:PGL2Q-stabilizer}
Consider a real quadratic irrational number $\alpha\in \R\Proj^1$, and let $a,b,c\in \Z$ be coprime such that $a\alpha^{2}+b\alpha+c=0$.
Under the projective action of $\PGL_{2}(\Q)$ on $\R\Proj^1$, the stabilizer subgroup of $\alpha$ is parametrized by $[t:s]\in \Q\Proj^1$ according to:
\begin{equation*}
    (t,s)\mapsto A = t \Id + s \begin{psmallmatrix} -b/2 & -c \\ a & +b/2 \end{psmallmatrix}\in \GL_2(\Q)
    \end{equation*}
Note that $\det(A)=t^2-s^2\Delta/4$. In particular, the stabilizer of $\alpha$ under the action of $\SL_2^\pm(\Q)$ corresponds to the rational points $(2t,s)\in \Q\times \Q$ on the Pell-Fermat conic $(2t)^2-s^2\Delta=\pm 4$.
\end{Lemma}
\begin{proof}
Suppose \(A = \begin{psmallmatrix} i & j  \\ k & l\end{psmallmatrix} \in \GL_2(\Q)\) satisfies $A\alpha=\alpha$, namely $k\alpha^{2}+(l-i)\alpha-j = 0$.
If $k=0$ then $i=l$ and $j=0$ hence $A = \Id$ in $\PGL_{2}(\Q)$, and this corresponds to $s=0$.
Otherwise, since the minimal polynomial of $\alpha$ over $\Q$ is projectively unique, there exists $s\in \Q\setminus\{0\}$ such that $(k,(l-i),-j)=s.(a,b,c)$, so we obtain the desired form for $A$.

Conversely any matrix of the given form fixes $\alpha$.
\end{proof}

\begin{Lemma}[stabiliser in $\PGL_2(\Z)$]
    \label{lem:PGL2Z-stabilizer}
    Consider a real quadratic irrational $\alpha\in \R$.
    
    There is a unique triple of coprime integers $a,b,c\in \Z$ such that $\alpha =\tfrac{-b+\sqrt{\Delta}}{2a}$ with $\Delta = b^2-4ac$. There is a unique primitive hyperbolic matrix \(F_\alpha \in \PGL_2(\Z)\) with attractive fixed point $\alpha$, and its entries are expressed in terms of the fundamental$(\pm)$-unit $\varepsilon_\Delta=t+s\tfrac{1}{2}\sqrt{\Delta}\in \mathcal{O}_\Delta^{\times}$ as:
    \begin{equation*} F_\alpha=t\Id+s\begin{pmatrix}
    -b/2 &-c \\ +a & +b/2
    \end{pmatrix}
    \qquad \forall e\in \Z \colon \:
    F_\alpha^e = 
    \tfrac{1}{2} (\varepsilon_\Delta^e +\varepsilon_\Delta^{-e})
    \Id + 
    \tfrac{1}{2\sqrt{\Delta}} (\varepsilon_\Delta^e -\varepsilon_\Delta^{-e})
    \begin{pmatrix}
        -b &-2c \\ +2a & +b
    \end{pmatrix}
    \end{equation*}

    Under the projective action of $\PGL_{2}(\Z)$ on $\R\Proj^1$, the stabiliser of $\alpha \in \R\Proj^1$ coincides with the maximal cyclic subgroup $\{F_\alpha^e \colon e \in \Z \}$ generated by $F_\alpha$.
    
    For two real quadratic irrationals $\alpha, \beta\in \R$, an element $C\in \PGL_2(\Q)$ satisfies $C\alpha = \beta$ if and only if $CF_\alpha C^{-1}=F_\beta$.
\end{Lemma}

\begin{proof}
    These classical results can be derived from the previous lemma and the structure of the discrete subgroup $\PSL_2(\Z)\subset \PGL_2(\R)$.
    We refer to \cite{Lachaud_cont-frac-zeta_1988} and \cite{CLS_phdthesis_2022, CLS_Conj-PSL2K_2022} for details.
\end{proof}

\begin{Remark}[expansion factor is unit]
    \label{rem:expansion-factor-unit}
    Let $M\in \PGL_2(\Z)$ such that $M^2$ is hyperbolic with attractive fixed point $\alpha\in \R$.
    % In that case we have $\tfrac{dM}{d\xi}(\alpha)\in \R$ is a unit of norm $1$.
    We have $M\alpha=\alpha$ and taking the Galois conjugate $M\Bar{\alpha}=\Bar{\alpha}$, therefore we have a diagonalization of $M$:
    \begin{equation*}
        M=\begin{psmallmatrix}
            m_{11} & m_{12} \\
            m_{21} & m_{22}
        \end{psmallmatrix}
        \quad
        % \tfrac{dM(\xi)}{d\xi} = \tfrac{\det(M)}{\lvert m_{21} \xi + m_{22}\rvert^2}
        % \qquad
        M\cdot
        \begin{psmallmatrix}
            \alpha \\ 1
        \end{psmallmatrix}
        = (m_{21} \alpha + m_{22})
        \begin{psmallmatrix}
            \alpha \\ 1
        \end{psmallmatrix}
        \quad
        M\cdot
        \begin{psmallmatrix}
            \Bar{\alpha} \\ 1
        \end{psmallmatrix}
        = (m_{21} \Bar{\alpha} + m_{22})
        \begin{psmallmatrix}
            \Bar{\alpha} \\ 1
        \end{psmallmatrix}
    \end{equation*}
    and the product of eigenvectors $(m_{21} \alpha + m_{22})(m_{21} \Bar{\alpha} + m_{22}) = \det(M)=\pm 1$ reveals that the expansion factor $(m_{21} \alpha + m_{22})\in \Z[\alpha]^\times$ is a quadratic unit. 

    % Note that we may reformulate this in terms of the derivative of $M$ at $\alpha$, namely $\tfrac{dM(\xi)}{d\xi} = \tfrac{\det(M)}{\lvert m_{21} \xi + m_{22}\rvert^2}$ takes a value at $\alpha$ that is a quadratic unit $\tfrac{dM(\alpha)}{d\xi} \in \mathcal{O}_\Delta^\times$.
\end{Remark}

\begin{Remark}[$\PGL_2(\Z)$-equivalence with multiples]
\label{rem:a-equiv-ma-mod-PGL2Z}
For $m\in \N_{>1}$, a real irrational $\alpha$ is $\PGL_2(\Z)$-equivalent to $m\alpha$ if and only if it is the root of a primitive integral quadratic polynomial $ax^2+bx+c\in \Z[x]$ with discriminant $\Delta=b^2-4ac$ such that $m\mid a$ and the Pell-Fermat equation $(2t)^2-s^2\Delta=\pm 4m$ has an integral solution $(t,s)\in \Z\times \Z$ with $2m\mid 2t+sb$.
\end{Remark}
\begin{proof}
Let $m\in \N_{>1}$ and consider a real irrational $\alpha$ which is $\PGL_2(\Z)$ equivalent to $m \alpha$. This implies that $\alpha$ is quadratic, so we may consider the unique triple of coprime integers $a,b,c\in \Z$ such that $\alpha = \frac{-b+\sqrt{\Delta}}{2a}$ where $\Delta=b^2-4ac$.
% (in particular $a\alpha^2+b\alpha+c=0$).
Since the discriminant $\Delta$ of such a quadratic number is invariant under the $\PGL_2(\Z)$-action, we find by writing $m\alpha = \tfrac{-mb+m\sqrt{\Delta}}{2a}=\tfrac{-b'+\sqrt{\Delta}}{2a'}$ that $a=ma'$ and $b'=b$, in particular $m$ must divide $a$.

Writing $M\alpha=m\alpha$ with \(M=\begin{psmallmatrix}
        i&j\\k&l
    \end{psmallmatrix}\in \GL_2(\Z)\), 
we find that $mk\alpha^2+(ml-i)\alpha-j=0$ hence there exists $s\in \Z$ such that $s.(a,b,c)=(mk, ml-i, -j)$.
Taking the discriminant yields an integral solution $(2t,s)$ to the Pell-Fermat equation $(2t)^2-s^2\Delta=4m\det(M)$ with $2t+sb=2ml$. 
% (In other terms $\pm m$ ramifies in $\Q(\sqrt{\Delta})$ (in particular $\pm m$ is a square $\bmod{\Delta}$).

Conversely if the Pell-equation $(2t)^2-s^2\Delta=\pm 4m$ has an integral solution $(t,s)$, we may determine $i,j,k,l\in \Q$ satisfying $il-kj=\pm 1$ from the equations $2t=ml+i$ and $s(a,b,c)=(mk, ml-i, -j)$, as $i=(2t-sb)/2$, $l=(2t+sb)/2m$, $k=sa/m$, $j=-sc$. If $m\mid a$ and $2m\mid (2t+sb)$ then $i,j,k,l\in \Z$.
\end{proof}

\begin{Question}[$\PGL_2(\Z)$-equivalence with multiples]
    \label{quest:a-equiv-ma-mod-PGL2Z}
    Can one provide alternative descriptions for the quadratic numbers in Remark~\ref{rem:a-equiv-ma-mod-PGL2Z} so as to answer the following questions:

    \begin{enumerate}[noitemsep]
        \item Can we construct an infinite family of examples in each quadratic field?
        \item What are the corresponding ideal classes of quadratic orders?
        \item Can we characterise the corresponding geodesics in the modular surface? % in terms of their lifts in congruence covers ?
    \end{enumerate}
\end{Question}

We are now ready to prove the main arithmetic result of this subsection, which we believe is novel and of independent interest.
%, and could have been conjectured by Lagrange or Gauss...

\begin{Theorem}[eventual $\PGL_2(\Z)$-equivalence]
\label{thm:eventual-PGL2Z}
% Two irrationals $\alpha,\beta\in \R\setminus \Q$ are eventually $\PGL_2(\Z)$-equivalent (namely there exists $m\in \N$ such that for all $n \ge m$ the numbers $n\alpha$ and $n\beta$ belong to the same $\PGL_2(\Z)$ orbit) if and only if $\alpha = \pm \beta \bmod{1}$.
Consider irrationals $\alpha,\beta\in \R\setminus \Q$.
Consider a subset of rational primes $\Primes'\subset \Primes$ having full Dirichlet density (say all but finitely many).
If for all $n\in \Primes'$, the numbers $n\alpha,n\beta$ are $\PGL_2(\Z)$-equivalent, then $\alpha=\pm \beta \bmod{\Z}$.
\end{Theorem}

\begin{proof}
Fix $m\in \Primes'$ and consider for each $n \in \Primes'_{\ge m}$ an element 
    \(C_{n} = \begin{psmallmatrix} 
    i_{n} & j_{n} \\ k_{n} & l_{n}
    \end{psmallmatrix}
    \in \PGL_{2}(\Z)\) 
satisfying $C_{n} \cdot (n\alpha) = n\beta$. 
Thus for any $n \in \Primes'_{\ge m}$ we have $B_{n} \alpha  = \beta$ where \(B_{n} = \begin{psmallmatrix} 
1/n & 0 \\ 0 & 1 \end{psmallmatrix} 
C_{n} \begin{psmallmatrix} 
n & 0 \\ 0 & 1 \end{psmallmatrix}\), hence the matrix $A_{n} = B_{m}^{-1}B_{n} \in \PGL_2(\Q)$ with $\det(A_{n})=\pm 1\bmod{(\Q^*)^2}$ stabilizes $\alpha\in \R\Proj^1$. 
% Since there is an infinite set of $n\ge m$ for which $B_{n}$ all have the same determinant in $\{\pm 1\}$, we may restrict to those $n$ and reset $m$ to be the first of them, so that for all $n\ge m$ we have $A_{n}\in \SL_2(\Q)$.
%
Note that since $m\alpha,m\beta$ are $\PGL_2(\Z)$-equivalent, they are either both non-quadratic or both generate the same quadratic field.

\emph{Suppose first that $\alpha,\beta$ have degree $>2$.}
In this case the stabilizer of $\alpha$ is $\{\pm \Id\}$, so we must have $A_{n} = \pm \Id$ so $B_{m} = \pm B_{n}$.
Writing this out in terms of the coefficients of $C_{m}, C_{n}\in \PGL_2(\Z)$ and letting $n\to \infty$ shows that $B_{m}$ is upper triangular.
Thus $C_{m}$ is upper triangular with diagonal entries $\pm 1$ and hence $\alpha=\pm \beta\bmod{\Z}$. 

\emph{Suppose now that $\alpha,\beta$ are quadratic irrational.}
%
% Now we show that $\alpha,\beta$ are $\PGL_2(\Z)$-equivalent. 
Consider the unique primitive $(a,b,c)\in \Z^3$ such that $\alpha=\tfrac{-b+\sqrt{\Delta}}{2a}$ with $\Delta=b^2-4ac$ (not a rational square, hence $ac\ne 0$).
By Lemma~\ref{lem:PGL2Q-stabilizer}, there exist $t_{n} , s_{n} \in \Q$ such that \(A_{n} = t_{n}\Id + s_{n} \begin{psmallmatrix} -b/2 & -c \\ a & b/2 \end{psmallmatrix}\in \PGL_{2}(\Q )\).
If there exist infinitely many $n \in \Primes'_{\ge m}$ such that $A_{n} = \pm \Id$ then we are done as for the non-quadratic case, so we now assume otherwise and consider only those $n \in \Primes'_{\gg m}$ large enough such that $A_{n} \ne \pm \Id$, hence $s_{n}\ne 0$.
The relation $\pm 1 = \det(A_{n}) = t_{n}^2 - s_{n}^2 \Delta/4$ implies that $t_n\ne 0$ (as $\Delta\in \N$ is not a square).

Spelling out $B_{m}^{-1}B_{n} = A_{n}$ in terms of the coefficients of $C_{m}, C_{n} \in \GL_2(\Z)$ yields:
\begin{equation}
    \label{BmBn=An} % \tag{Bn=BmAn}
    \begin{pmatrix} 
    l_{m} i_{n} -\tfrac{1}{m} j_{m} n k_{n}  & l_{m} \tfrac{1}{n} j_{n} - \tfrac{1}{m} j_{m} l_{n} \\ 
    -mk_{m} i_{n} + i_{m} n k_{n} & -mk_{m} \tfrac{1}{n}j_{n} + i_{m} l_{n} 
    \end{pmatrix}
    =
    \begin{pmatrix} 
    t_{n}-s_{n}b/2 & -s_{n}c \\ 
    s_{n}a & t_{n}+s_{n}b/2 \end{pmatrix} 
\end{equation}
Equating the lower-left entries shows that $s_na\in \Z$, so for every prime $p\nmid a$ we have $\val_p(s_n)\ge 0$.
Assume that $n\in \Primes'_{\gg m}$ is coprime to $abcm$ and all non-zero entries of $C_m$.
Equating the difference of diagonal entries shows that $n \mid j_{n}k_{m}$ hence either $n \mid j_n$ or $k_m=0$.
Equating the upper-right entries shows that $n\mid j_{n}l_{m}$ hence either $n\mid j_{n}$ or $l_{m}=0$.
Since $(k_m,l_m)\ne (0,0)$ we must have $n\mid j_n$, implying that $B_{n}\in \PGL_2(\Z)$.
This shows that $\alpha$ and $\beta$ are $\PGL_2(\Z)$-equivalent.

Consequently we have $C_1\in \PGL_2(\Z)$ such that $C_1\alpha=\beta$, and repeating the previous reasoning with $m=1$ shows that if $n\in \Primes'$ is coprime to $abc$ and the non-zero entries of $C_1$, then we have $B_n\in \GL_2(\Z)$ satisfying $B_n=C_1A_n$ where $A_n\in \PGL_2(\Z)$ stabilises $\alpha$. 
\begin{comment}
Spelling out $B_n=C_1A_n$ yields
\begin{equation}
    \label{Bn=C1An} % \tag{Bn=BmAn}
    \begin{pmatrix} 
    i_{n} & j_{n}/n \\ 
    nk_{n} & l_{n} 
    \end{pmatrix}
    =
    \begin{pmatrix} 
    i_{m}(t_{n}-s_{n}b/2)+j_{m}s_{n}a  & 
    j_{m}(t_{n}+s_{n}b/2)-i_{m}s_{n}c \\ 
    k_{m}(t_{n}-s_{n}b/2)+l_{m}s_{n}a & 
    l_{m}(t_{n}+s_{n}b/2)-k_{m}s_{n}c
    \end{pmatrix} 
\end{equation}
\end{comment}

By Lemma~\ref{lem:PGL2Z-stabilizer}, there is $e_n\in \Z$ such that $A_n = F_\alpha^{e_n} = t_{n}\Id + s_{n} \begin{psmallmatrix} -b/2 & -c \\ a & b/2 \end{psmallmatrix}\) where we have, in terms of the fundamental unit $\varepsilon_\Delta\in \mathcal{O}^\times_\Delta$ (as in Remark~\ref{rem:quadratic-orders}), the formulae $t_n= \tfrac{1}{2}(\varepsilon_\Delta^{e_n}+\Bar{\varepsilon}_\Delta^{e_n})$ and $s_n =\tfrac{1}{2\sqrt{\Delta}} (\varepsilon_\Delta^{e_n}-\Bar{\varepsilon}_\Delta^{e_n})$.
Equating the lower-left entries of $B_n = C_1 A_n$ yields:
\begin{equation}
    \label{eq:nkn_Bn=C1.An}
    \rvert nk_n \lvert 
    = \lvert k_1(t_n-s_n b/2)+l_1 s_n a \rvert
    = \lvert \varepsilon_\Delta^{e_n}\Bar{\eta}-\Bar{\varepsilon}_\Delta^{e_n} \eta \rvert
    \quad \mathrm{where} \quad
    2\sqrt{\Delta} \eta 
    = (l_1a-k_1b/2)+k_1\sqrt{\Delta}
\end{equation}
Recall that if $k_1=0$ then $\alpha=\beta \bmod{\Z}$ so we may assume otherwise hence $\eta\ne 0$.
Now assume further that the prime $n\in \Primes'$ is coprime to $\eta$ in $\Q(\sqrt{\Delta})$, and multiply both terms of \eqref{eq:nkn_Bn=C1.An} by $\varepsilon_\Delta^{e_n}/\Bar{\eta}$ to find that $n$ divides $\left| \varepsilon_\Delta^{2e_n}-\eta/\Bar{\eta} \right|$.
Since this holds for a positive density of primes $n\in \Primes'$, we deduce from Schinzel's local-to-global theorem on power residues \cite[Theorem 2]{Schinzel_power-residues-exponential-congruences_1975} that there exists $e_0\in \N$ such that $\epsilon^{2e_0}=\eta/\Bar{\eta}$.

Now observe that $C_1 F_\alpha^{e_0}\in \PGL_2(\Z)$ has its lower left entry equal to $0$, and by definition it sends $\alpha$ to $\beta$, proving that $\alpha = \pm \beta \bmod{\Z}$.
\end{proof}

\begin{Remark}[acknowledgments]
In a previous version of this work, we had left Theorem~\ref{thm:eventual-PGL2Z} as a conjecture, and only showed that when $\alpha,\beta\in \R\setminus \Q$ are eventually $\PGL_2(\Z)$-equivalent:
\begin{itemize}
    \item If $\alpha,\beta$ are not quadratic then $\alpha = \pm \beta \bmod{\Z}$.
    \item If $\alpha,\beta$ are quadratic then for all $n\in \N$ the numbers $n\alpha, n\beta$ are $\PGL_2(\Z)$-equivalent.
    Moreover, we have $\alpha = \pm \beta \bmod{\Q}$.
\end{itemize}
Soon after, Daniel Martin communicated a sketch to complete the proof of the quadratic case using Kummer theory, which we simplified and refined using Schinzel's local-to-global principle.

In addition, François Maucourant communicated to us the following geometric argument showing that if quadratic $\alpha,\beta\in \R\setminus \Q$ are eventually $\PGL_2(\Z)$-equivalent, then $\alpha=\pm \beta \bmod{\Q}$.
% The group $\PSL_2(\R)$ acts by projective transformations of the upper-half plane $\HP=\{z\in \C\colon \Im(z)>0\}$, as the group of orientation preserving isometries of the hyperbolic plane.
% A real quadratic number $\gamma$ determines a geodesic between $\gamma$ and its Galois conjugate $\Bar{\gamma}$, whose height is $h(z) = \tfrac{1}{2}\lvert \gamma-\Bar{\gamma}\rvert$.
Define the height of a real quadratic number $\gamma$ with Galois conjugate $\Bar{\gamma}$ as $h(z) = \lvert \gamma-\Bar{\gamma}\rvert$, and say that $\gamma$ maximizes the height in its $\PGL_2(\Z)$-class when for any $\gamma'$ in its $\PGL_2(\Z)$-orbit, we have $h(\gamma')\le h(\gamma)$.
Notice that if $\gamma$ is a quadratic integer, then it maximises its height in its $\PGL_2(\Z)$-class, since for every 
$C=\begin{psmallmatrix}
    i&j\\k&l
\end{psmallmatrix}\in \PGL_2(\Z)$ we have
\begin{equation*}
    \lvert C\gamma - C\Bar{\gamma} \rvert 
    = \left| \tfrac{\det(C) \cdot (\gamma - \Bar{\gamma})}{(j\gamma+l)(j\Bar{\gamma}+l)} \right|
    = \tfrac{\lvert \gamma - \Bar{\gamma}\rvert}{\lvert (j\gamma+l)(j\Bar{\gamma}+l)\rvert}
\end{equation*}
and the denominator is the norm of the quadratic integer $\lvert j\gamma+l\rvert$, that is a non-zero integer.
There exists (a common denominator) $n\in\N$ such that $n\alpha, n\beta$ are both integers in $\Q(\sqrt{\Delta})$, so they both maximize the height in their common $\PGL_2(\Z)$-class, hence $n\lvert\alpha-\Bar{\alpha}\rvert = n\lvert\beta-\Bar{\beta}\rvert$.
We deduce that $\alpha = \pm \beta \bmod{\Q}$.
\end{Remark}

Let us mention an example which hints at some of the subtlety in Theorem~\ref{thm:eventual-PGL2Z}.

\begin{Example}[almost counter-example]
    \label{ex:p-eventually-PSL2Z-equivalent}
    Denoting $\phi = \tfrac{1}{2}(1+\sqrt{5})$, 
    set $\alpha = \tfrac{1}{2}\phi=\Ecf{0;1,(4,4)^\N}$ and $\beta = \tfrac{1}{2}(\phi-1)=\Ecf{0;3,(4,4)^\N}$. We have that $\alpha \ne \pm \beta \bmod{\Z}$ but they are $\PGL_2(\Z)$ equivalent.
    Moreover, for many $f\in \N$, we have that $f\alpha$ and $f\beta$ are $\PGL_2(\Z)$-equivalent: this holds when $f$ is even (since $f\alpha-f\beta=f/2$), and experimentally for many odd values of $f$ (in particular all those for which the order of discriminant $5f^2$ is principal). 
    However by Theorem~\ref{thm:eventual-PGL2Z}, there must be some $f\in \N$ such that $f\alpha, f\beta$ are not $\PGL_2(\Z)$-equivalent: the smallest is $f=9$, for which we have $9\alpha = \Ecf{7; 3, (1, 1, 3, 1, 9, 3,)^\N}$ and $9\beta = \Ecf{2; 1, (3, 1, 1, 3, 9, 1)^\N}$.
    
    A similar phenomenon can be observed with $\alpha = \tfrac{1}{5}\sqrt{2}$ and $\beta = \tfrac{1}{5}(\sqrt{2}+2)$.
\end{Example}

The previous example shows that the assumption in Theorem~\ref{thm:eventual-PGL2Z} on the full Dirichlet density of primes would be difficult to generalize.
The following question makes this more precise.

\begin{Question}[$\N'$-virtual $\PGL_2(\Z)$-equivalence]
    \label{quest:weak-eventual-PGLZ}
    For a subset of $\N' \subset \N$, say that $\alpha,\beta$ are $\N'$-virtually $\PGL_2(\Z)$-equivalent when for all $n\in \N'$ the numbers $n\alpha,n\beta$ are $\PGL_2(\Z)$-equivalent.

    For which subsets $\N'$ does $\N'$-virtual $\PGL_2(\Z)$-equivalence imply $\PGL_2(\Z)$-equivalence, or further equality up to sign and addition by an integer?
    
    Theorem~\ref{thm:eventual-PGL2Z} replies both positively when $\N'$ is a subset of primes of full Dirichlet density while Example~\ref{ex:p-eventually-PSL2Z-equivalent} shows that both are false when $\N'$ consists of the positive powers of a single prime (and we believe a similar construction extends this counter examples to $\N'$ the positive powers of any finite subset of primes.
    What if $\N'$ is a subset of primes of positive Dirichlet density?
\end{Question}

\subsection{Eventual equivalencies between Sturmian systems}

\label{subsec:eventual-flow-Sturmian}

Recall the diagram for virtual and eventual equivalencies (conjugacy and flow) of dynamical systems from Remark~\ref{rem:implications-virtual-eventual-flow}.
The goal of this section is to complete it for Sturmian systems into the following diagram; in which all non-indicated implications are false.
The main result is Theorem~\ref{thm:eventual-flow-Sturmian} saying that eventual flow equivalence implies topological conjugacy.

\begin{center}
\begin{tikzpicture}[node distance=4.5cm and 1cm,
box/.style = {draw, minimum height=12mm, inner xsep=1mm, align=center},
%sy+/.style = {yshift= 2mm},
%sy-/.style = {yshift=-2mm},
every edge quotes/.style = {align=center}]
    % Define nodes
    \node (A) [box] at (0, 0) {Topological conjugacy \\ $\pm \alpha \bmod{\Z}$ \\ Lemma~\ref{lem:equal-rank2-lattices}};
    \node (B) [box, right of=A, xshift=.5cm] {Eventual conjugacy \\ $\pm \alpha \bmod{\Z}$ \\ Proposition~\ref{prop:eventual-conj-to-conj-Sturmian}};
    \node (C) [box, right of=B, xshift=.5cm] {Virtual conjugacy \\ $\pm \alpha \bmod{\Q}$ \\ Proposition~\ref{prop:virtual-conjugacy-Sturmian}};
    \node (D) [box, below of=A, yshift=2.5cm] {Flow equivalence \\ $\alpha \bmod{\PGL_2(\Z)}$ \\ Theorem~\ref{thm:flow-equivalence-Sturmian-systems}};
    \node (E) [box, below of=B,yshift=2.5cm] {Eventual flow equivalence \\ $\pm \alpha \bmod{\Z}$ \\ Theorem~\ref{thm:eventual-flow-Sturmian}};
    \node (F) [box, below of=C, yshift=2.5cm] {Virtual flow equivalence \\ $n\alpha\bmod{\PGL_2(\Z)}$ \\ Remark \ref{rem:virtual-flow-equivalence-Sturmian}};

    % Draw arrows
    \draw[->, double] (A) -- (B);
    \draw[->, double] (B) -- (C);
    \draw[->, double] (A) -- (D);
    \draw[->, double] (B) -- (E);
    \draw[->, double] (C) -- (F);
    % \draw[->, double] (D) -- (E);
    \draw[->, double] (E) -- (F);
    \draw[->, double] (B) -- (A);
    \draw[->, double] (E) -- (B);
    \draw[->, double] (E) -- (A);
    \draw[->, double] (A) -- (E);
\end{tikzpicture}
\end{center}

Let us begin by characterizing virtual conjugacy of Sturmian systems.

\begin{Proposition}[virtual conjugacy of Sturmian systems]
    \label{prop:virtual-conjugacy-Sturmian}
    Two Sturmian systems $(X_{\alpha},\Shift_{\alpha})$ and $(X_{\beta},\Shift_{\beta})$ are virtually conjugate if and only if $\alpha=\pm \beta\bmod{\Q}$.
\end{Proposition}
\begin{proof}
    The proof is easy and left as an exercise.
\end{proof}

\begin{Remark}[virtual flow equivalence of Sturmian systems]
    \label{rem:virtual-flow-equivalence-Sturmian}
    Two Sturmian systems $(X_{\alpha},\Shift_{\alpha})$ and $(X_{\beta},\Shift_{\beta})$ are virtually flow equivalent if and only if $n\alpha=\pm n\beta\bmod{\PGL_2(\Z)}$, in other terms when there is $C\in \bigcup_{n\in \N_{>0}}\operatorname{diag}(n,1)\PGL_2(\Z)\operatorname{diag}(n,1)^{-1}$ such that $C\cdot \alpha = \beta$.
\end{Remark}

For Sturmian systems, it is not hard to show that eventual conjugacy implies conjugacy.
\begin{Proposition}[eventually conjugate implies conjugate]
\label{prop:eventual-conj-to-conj-Sturmian}
Sturmian subshifts are eventually conjugate if and only if they are conjugate.
\end{Proposition}
\begin{proof}
For the nontrivial direction, suppose $(X_{\alpha},\Shift_{\alpha}^{n})$ and $(X_{\beta},\Shift_{\beta}^{n})$ are topologically conjugate for all $n \ge m$. Let $p \ge m$ be prime. Then we have $p \alpha = \epsilon_{0} p \beta + k_0$ and $(p+1)\alpha = \epsilon_{1} (p+1) \beta + k_1$ for some $k_i \in \Z$ and $\epsilon_{i} \in \{-1,+1\}$. Thus $\alpha = \epsilon_{0} \beta + \frac{k_0}{p}$ and $\alpha = \epsilon_{1} \beta + \frac{k_1}{p+1}$, so $(\epsilon_{0}-\epsilon_{1}) \beta = \frac{k_1}{p+1}-\frac{k_0}{p}$. 
We must have $\epsilon_{0}=\epsilon_{1}$ since $\beta$ is irrational, hence $\frac{k_1}{p+1}=\frac{k_0}{p}$.
Since $p$ is prime, it must divide $k_0$, so $\alpha = \epsilon_{0} \beta \bmod{\Z}$.
\end{proof}

We now consider the relation of eventual flow equivalence for Sturmian systems. For nonquadratic irrationals, we will show that it implies topological conjugacy. This also gives some motivation for our introduction of isogenies, yielding a richer equivalence relation which connects to the arithmetic of continued fractions.

\begin{Remark}[if Sturmian have flow equivalent powers then powers coincide]
\label{rem:powers-Denjoy-flow}
Note that if Sturmian systems $(X_{\alpha},\Shift_{\alpha})$ and $(X_{\beta},\Shift_{\beta})$ are such that $(X_{\alpha},\Shift_{\alpha}^m)$ and $(X_{\beta},\Shift_{\beta}^n)$ are flow equivalent for some $m,n\ge 1$, then $m=n$. 
Indeed, the Denjoy system $(X_{\alpha}^m,\Shift_{\alpha}^m)$ has exactly $m$ distinct asymptotic orbits, hence $\Sigma_{\Shift_{\alpha}^{m}}X$ has exactly $m$ pairs of asymptotic leaves. Since a flow equivalence must carry asymptotic leaves to asymptotic leaves, this implies that $m=n$.
\end{Remark}

\begin{Lemma}[eventual flow equivalence of Sturmian systems]
\label{lem:eventual-flow-Sturmian}
Consider for $k\in\{0,1\}$ two Sturmian systems $(X_{\alpha_k},\Shift_{\alpha_k})$  and let $m\in \N_{\ge 1}$. The systems $(X_{\alpha_0},\Shift_{\alpha_0}^{m})$ and $(X_{\alpha_1},\Shift_{\alpha_1}^{m})$ are flow equivalent if and only if there exists $C_{m} \in \PGL_{2}(\Z)$ such that $C_{m} \cdot (m\alpha_0) = m\alpha_1$. 
In particular $(X_{\alpha_0},\Shift_{\alpha_0})$ and $(X_{\alpha_1},\Shift_{\alpha_1})$ are eventually flow equivalent if and only if for all but finitely many $m\in \N_{\ge 1}$, there exists $C_{m} \in \PGL_{2}(\Z)$ such that $C_{m} \cdot (m\alpha) = m\beta$.
\end{Lemma}
\begin{proof}
We cannot apply Theorem~\ref{thm:flow-equivalence-Sturmian-systems} since a proper power of Sturmian system is not Sturmian. 
However for $k=0,1$, both $(X_{\alpha_k},\Shift_{\alpha_k}^{m})$ are Denjoy systems with invariants $(m\alpha_k,A_k)$, 
where $\tilde{A}_{k} = \Z+\alpha_k\Z \subset \R$ are the lifts of $A_{k}$ to $\R$. Choosing a representative \(\begin{psmallmatrix}
    a_{m} & b_{m} \\ c_{m} & d_{m}
\end{psmallmatrix} \in \GL_{2}(\Z)\) of $C_{m} \in \PGL_{2}(\Z)$, we have by Lemma~\ref{lem:proportional-rank2-lattices} that $\frac{1}{c_{m}\alpha_0+d_{m}}\tilde{A}_{0} = \frac{1}{c_{m}\alpha+d_{m}}\Z +\alpha_0 \Z = \Z + \alpha_1 \Z = \tilde{A}_{1}$, so the result follows from Theorem~\ref{thm:flow-equivalence-Denjoy-systems}.
\end{proof}

\begin{Theorem}[eventual flow implies conjugacy]
\label{thm:eventual-flow-Sturmian}
For Sturmian systems, eventual flow equivalence implies topological conjugacy (hence flow equivalence and eventual conjugacy).
\end{Theorem}

\begin{proof}
    Suppose that two Sturmian systems $(X_{\alpha},\Shift_{\alpha})$ and $(X_{\beta},\Shift_{\beta})$ are eventually flow equivalent.
    By Lemma~\ref{lem:eventual-flow-Sturmian}, their rotation numbers $\alpha, \beta$ must be eventually $\PGL_2(\Z)$-equivalent: there exists $m$ such that for all $n \ge m$ there exists \(C_{n} \in \GL_{2}(\Z)\) 
    such that $C_{n} \cdot (n\alpha) = n\beta$. 
    By Theorem~\ref{thm:eventual-PGL2Z}, this implies that $\alpha=\pm\beta \bmod{\Z}$.
\end{proof}

\begin{Corollary}[flow equivalent but not eventually flow equivalent]
For Sturmian systems, flow equivalence does not imply eventual flow equivalence.
\end{Corollary}

\begin{comment}
\begin{Proposition}
    For Sturmian systems, virtual flow equivalence does not imply flow equivalence.
\end{Proposition}
\end{comment}

The following examples exhibit Sturmian systems which are virtually flow equivalent but not flow equivalent.

\begin{Example}[virtual flow equivalent but not flow equivalent]
\label{ex:virtual-flow-not-flow}
Let $\alpha$ be an irrational which is not quadratic, and let \(C_1 = \begin{psmallmatrix} 2 & 1 \\ 1 & 1 \end{psmallmatrix}\). 

The Sturmian systems associated to $\alpha$ and $\beta = C_1 \alpha$ are flow equivalent by Theorem~\ref{thm:flow-equivalence-Sturmian-systems}. 
However for $n \ge 2$, the Sturmian systems with rotation numbers $\frac{1}{n}\alpha$ and $\frac{1}{n}\beta$ are not flow equivalent. 
Indeed, if they were there would exist \(C_n \in \GL_{2}(\Z)\) such that $C_n \frac{1}{n}\alpha = \frac{1}{n}\beta$, hence $B_{n} \alpha = \beta$ where \(B_{n} = \begin{psmallmatrix} n & 0 \\ 0 & 1 \end{psmallmatrix} C_n \begin{psmallmatrix} 1/n & 0 \\ 0 & 1 \end{psmallmatrix} \in \GL_2(\Q)\). 
Thus $C_1^{-1} B_{n} \in \GL_2(\Q)$ stabilizes $\alpha$, so $B_{n}=\pm C_1$, hence the upper right coefficient of $C_n$ is $\pm 1/n$, leading to a contradiction once $n \ge 2$.
\end{Example}

A slightly more involved example shows that even for quadratic irrationals, two Sturmian systems can be virtually flow equivalent but not flow equivalent.

\begin{Example}[virtual flow equivalent but not flow equivalent for quadratic]
\label{ex:virtual-flow-not-flow-quadratic}
Let $\alpha = \sqrt{2}$, \(C_1 = \begin{psmallmatrix} 2 & 1 \\ 1 & 1 \end{psmallmatrix} \in \GL_{2}(\Z)\) and set $\beta = C_1 \alpha = 3-\sqrt{2}$.

The Sturmian systems associated to $\alpha$ and $\beta$ are flow equivalent by Theorem~\ref{thm:flow-equivalence-Sturmian-systems}. 
We claim that the Sturmian systems associated to $\frac{1}{2}\alpha$ and $\frac{1}{2}\beta$ are not flow equivalent.
For if they were, there would exist \(C_2 = \begin{psmallmatrix} i & j \\ k & l \end{psmallmatrix} \in \GL_{2}(\Z)\) such that \(B_{2} = \begin{psmallmatrix} i & 2j \\ k/2 & l \end{psmallmatrix}\) satisfies $B_2\alpha = \beta$.
Hence \(C_1^{-1}B_{2}=\begin{psmallmatrix}
    i-k/2 & 2j-l \\ -i+k & -2j+2l
\end{psmallmatrix} \in \GL_{2}(\Q)\) stabilizes $\alpha$ which by Lemma~\ref{lem:PGL2Q-stabilizer} means it has the form \(C_1^{-1} B_{2} = \begin{psmallmatrix} r & 2s \\ s & r \end{psmallmatrix}\) for some $r, s \in \Q$.
Equating the ratio of the off-diagonal coefficients leads to $2\mid l$.
Equating the difference of diagonal coefficients leads to $2\mid k$.
This contradicts $\gcd(k,l)=1$.
\end{Example}

\begin{Remark}[self-isogenies of Sturmian systems]
\label{rem:flow-equiv-to-power-Sturmian}
For an integer $m>1$, the Sturmian systems with parameters $\alpha$ and $m\alpha$ are flow equivalent if and only if $\alpha$ is the root of a primitive integral quadratic polynomial $ax^2+bx+c\in \Z[x]$ with discriminant $\Delta=b^2-4ac$ such that $m\mid a$ and the Pell-Fermat equation $(2t)^2-s^2\Delta=\pm 4m$ has a solution $(t,s)\in \Z\times \Z$ with $2m\mid 2t+sb$.
This follows from Remark~\ref{rem:a-equiv-ma-mod-PGL2Z}. 
Question~\ref{quest:a-equiv-ma-mod-PGL2Z} asks about alternative descriptions for such numbers $\alpha$, and if one may construct infinite families of examples in each quadratic field.
\end{Remark}

\subsection{Eventual flow equivalence of Denjoy systems}
\label{subsec:eventual-flow-Denjoy}

Recall the diagram for virtual and eventual equivalencies (conjugacy and flow) of dynamical systems from Remark~\ref{rem:implications-virtual-eventual-flow}.
The goal of this section is to complete it for Denjoy systems into the following diagram; in which all non-indicated implications are false.
The main result is Theorem~\ref{thm:eventual-flow-Denjoy} saying that eventual flow equivalence implies topological conjugacy.

\begin{center}
\begin{tikzpicture}[node distance=4.5cm and 1cm,
box/.style = {draw, minimum height=12mm, inner xsep=1mm, align=center},
%sy+/.style = {yshift= 2mm},
%sy-/.style = {yshift=-2mm},
every edge quotes/.style = {align=center}]
    % Define nodes
    \node (A) [box] at (0, 0) {Topological conjugacy \\ $\pm (\alpha, A) \bmod{\Z}$ \\ Lemma~\ref{lem:equal-rank2-lattices}};
    \node (B) [box, right of=A, xshift=.5cm] {Eventual conjugacy \\ $\pm (\alpha,A) \bmod{\Z}$ \\ Theorem \ref{thm:eventual-flow-Denjoy} };
    \node (C) [box, right of=B, xshift=.5cm] {Virtual conjugacy \\ $(\pm\alpha \bmod{\Q}, A)$ \\ Proposition \ref{prop:virtual-conjugacy-Denjoy}} ;
    \node (D) [box, below of=A, yshift=2.5cm] {Flow equivalence \\ $(\alpha,A) \bmod{\PGL_2(\Z)}$ \\ Theorem \ref{thm:flow-equivalence-Denjoy-systems} };
    \node (E) [box, below of=B,yshift=2.5cm] {Eventual flow equivalence \\ $\pm (\alpha, A) \bmod{\Z}$ \\ Theorem \ref{thm:eventual-flow-Denjoy}};
    \node (F) [box, below of=C, yshift=2.5cm] {Virtual flow equivalence \\ $(n\alpha, A)\bmod{\PGL_2(\Z)}$ \\ Remark~\ref{rem:virtual-flow-equivalence-Denjoy} };

    % Draw arrows
    \draw[->, double] (A) -- (B);
    \draw[->, double] (B) -- (C);
    \draw[->, double] (A) -- (D);
    \draw[->, double] (B) -- (E);
    \draw[->, double] (C) -- (F);
    % \draw[->, double] (D) -- (E);
    \draw[->, double] (E) -- (F);
    \draw[->, double] (B) -- (A);
    \draw[->, double] (E) -- (B);
    \draw[->, double] (E) -- (A);
    \draw[->, double] (A) -- (E);
\end{tikzpicture}
\end{center} 

We begin by characterizing virtual conjugacy of Denjoy systems.

\begin{Proposition}[virtual conjugacy of Denjoy systems]
    \label{prop:virtual-conjugacy-Denjoy}
Consider for $k \in \{0,1\}$ Denjoy systems $(\Denjoy_{F_k}, F_k)$ with (rotation number, cut-point set) invariants $(\alpha_k, A_k)$ (with $0\in A_k$). 
They are virtually conjugate if and only if $\alpha_0=\pm \alpha_1\bmod{\Q}$ and $A_0 \equiv A_1$.
\end{Proposition}

\begin{proof}
    This follows from the classification of Denjoy systems up to topological conjugacy (§\ref{subsec:denjoy-systems}), noting that for $n\in \N_{>0}$ the systems $(\Denjoy_{\alpha_k, A_k}, F_{\alpha_k, A_k}^{n})$ are Denjoy with invariants $(n\alpha_k,A_k)$.
\end{proof}

The following contrasts the Remark \ref{rem:powers-Denjoy-flow} for Sturmian systems.

\begin{Remark}[if Denjoy have conjugate powers then powers may not coincide]
For any $m,n\in \N_{\ge 1}$, one may find pairs of Denjoy systems $(\Denjoy_{F_0}, F_0), (\Denjoy_{F_1}, F_1)$ with rotation numbers satisfying $m\alpha_0=n\alpha_1$
such that $(\Denjoy_{F_0}, F_0^m)$ and $(\Denjoy_{F_1}, F_1^n)$ are topologically conjugate (hence flow equivalent).
We leave the construction as an exercise.
\end{Remark}

\begin{Remark}[virtual flow equivalence of Denjoy systems]
    \label{rem:virtual-flow-equivalence-Denjoy}
    For $k \in \{0,1\}$, the Denjoy systems $(\Denjoy_{F_k}, F_k)$ with invariants $(\alpha_k, A_k)$ % (satisfying $0\in A_k$). 
    are virtually flow equivalent if and only if there exists $M\in \bigcup_{n\in \N_{>0}}\operatorname{diag}(n,1)\PGL_2(\Z)\operatorname{diag}(n,1)^{-1}$ such that \eqref{eq:Denjoy-flow-equivalence} holds. 
\end{Remark}

\begin{Theorem}[eventual flow implies conjugacy]
    \label{thm:eventual-flow-Denjoy}
    For Denjoy systems, eventual flow equivalence implies topological conjugacy (hence flow equivalence and eventual conjugacy).
\end{Theorem}
\begin{proof}
    Consider for $k \in \{0,1\}$ Denjoy systems $(\Denjoy_{F_k}, F_k)$ with (rotation number, cut-point set) invariants $(\alpha_k, A_k)$, normalised so that $0\in A_k$ and $\alpha_k\in (0,1)$.
    % Denote $\tilde{A}_k\subset \R$ the lifts of $A_k\subset \R\bmod{\Z}$: they are $\Z$-modules invariant by translation of $\alpha_k$.

    % Recall that by Corollary~\ref{cor:2-AI-factor-Denjoy-power-m}: for all $n \in \N'$, letting $A_k^{\prime}$ a set containing one point from each $R_{\alpha_k}$-orbit in $A_k$ such that $0 \in A_k^{\prime}$ and ${A_{k}^{\prime}(n) = \{R_{n\alpha}^{e}(x) \colon x \in A_k^{\prime}, e \in \Z \}}$, there is a 2-AI factor between Denjoy systems $(\Denjoy_{\alpha_k, A_k}, F_{\alpha_k, A_k}^{n}) \to (\Denjoy_{n\alpha, A_{k}(n)^{\prime}}, F_{n\alpha_k, A_{k}^{\prime}(n)})$.

    For each $n\in \N_{>0}$ the systems $(\Denjoy_{\alpha_k, A_k}, F_{\alpha_k, A_k}^{n})$ are Denjoy with invariants $(n\alpha_k,A_k)$, so by Theorem~\ref{thm:flow-equivalence-Denjoy-systems} they are flow equivalent if and only if there is $M^{(n)}= (m^{(n)}_{ij}) \in \GL_2(\Z)$ such that:
    \begin{equation}
        \label{eq:Denjoy-n-flow}
        n\alpha_1 = \frac{m^{(n)}_{11}\cdot n\alpha_0+m^{(n)}_{12}}{m^{(n)}_{21}\cdot n\alpha_0+m^{(n)}_{22}} 
        \qquad \mathrm{and} \qquad 
        A_1 \equiv \frac{1}{m^{(n)}_{21}\cdot n\alpha_0+m^{(n)}_{22}} \cdot A_0
    \end{equation}
    where $\equiv$ denotes equality of subsets of $\Torus$ up to rotation.

    By assumption the systems $(\Denjoy_{\alpha_k, A_k}, F_{\alpha_k, A_k})$ are eventually flow equivalent so there is a cofinite subset $\N'\subset \N_{>0}$ such that for all $n\in \N'$ we have \eqref{eq:Denjoy-n-flow}. 
    In particular, $\alpha_0,\alpha_1$ are eventually $\PGL_2(\Z)$-equivalent hence by Theorem~\ref{thm:eventual-PGL2Z} we have $\alpha_0=\pm \alpha_1\bmod{\Z}$.
    Up to exchanging the sign of $(\alpha_1, A_1)$, we now assume $\alpha_0=\alpha_1=\alpha\in (0,1)$.
    If $\alpha$ is not quadratic, then the stabiliser of $n\alpha$ in $\PGL_2(\Z)$ is trivial so for all (hence for some) $n\in \N'$, we have $M^{(n)}=\pm \Id$ which proves $A_1\equiv \pm A_0$ as desired.
    Now suppose that $\alpha$ generates the quadratic order $\mathcal{O}_{\Delta}$ of discriminant $\Delta$ with fundamental unit $\varepsilon_\Delta\in \Z[\alpha]$. 
    Note that for all $n\in \N$, the quadratic order $\Z[n\alpha]$ has index $n$ in $\Z[\alpha]$ so it is $\Z[n\alpha]=\mathcal{O}_{n^2\Delta}$ hence by Remark~\ref{rem:Dirichlet-unit} it has fundamental unit $\varepsilon_\Delta^n$.
    
    For every $n\in \N'$, we know by Remark~\ref{rem:expansion-factor-unit} that the number $\lambda^{(n)} = m_{21}^{(n)} n\alpha_0+m_{22}^{(n)}\in \Z[n\alpha]$ is a unit in $\Z[n\alpha]=\mathcal{O}_{n^2\Delta}$  hence by Remark~\ref{rem:Dirichlet-unit} that there is $e_n\in \N$ such that $\pm \lambda^{(n)}=\varepsilon_\Delta^{n.e_n}$.
    If for some $n\in \N'$ we have $\lambda^{(n)}=\pm \Id$ then we are done, so assume otherwise namely $e_n> 0$.
    
    Fix $m,n\in \N'$ we have $\varepsilon_\Delta^{me_m-ne_n}A_0\equiv A_0$.
    Let $G\subset \Z[\alpha]^\times$ be the subgroup generated by all $\pm \varepsilon_\Delta^{me_m-ne_n}$, so that $G A_0\equiv A_0$ and $G A_1\equiv A_1$.
    Denoting $g=\gcd\{me_m-ne_n \colon m,n\in \N'\}>0$ we have $G=\{\pm (\epsilon_\Delta^g)^e \colon e\in \Z\}$.
    Now observe that $\N'\cap g\N$ is non-empty, and for $n\in \N'\cap g\N$ we have $\pm \lambda^{(n)}=\varepsilon_\Delta^{n e_n} \in G$, therefore $A_0\equiv \lambda^{(n)} A_1 \equiv A_1$.
\end{proof}

\begin{comment}
\section*{Declarations}
\paragraph{Ethical approval:}
This declaration is not applicable.
\paragraph{Funding:}
This declaration is not applicable.
\paragraph{Availability of data and materials:}
This declaration is not applicable.
\end{comment}

%\nocite{Boyle-Carlsen-Eilers_flow-equivalence-subshifts_2017}
\bibliographystyle{alpha}
\bibliography{biblio.bib}

@book{Cohn_algebraic-numbers-class-fields_1978,
    AUTHOR = {Cohn, Harvey},
    TITLE = {A classical invitation to algebraic numbers and class fields},
    SERIES = {Universitext},
    NOTE = {With two appendices by Olga Taussky: ``Artin's 1932 G\"{o}ttingen lectures on class field theory'' and ``Connections between algebraic number theory and integral matrices''},
    PUBLISHER = {Springer-Verlag, New York-Heidelberg},
    YEAR = {1978},
    PAGES = {xiii+328},
    ISBN = {0-387-90345-3},
    MRCLASS = {12-01 (12Axx)},
    MRNUMBER = {506156},
    MRREVIEWER = {Ezra Brown},
}

@book{Hatcher_topology-numbers_2022,
    author = {Hatcher, Allen},
    title = {Topology of numbers},
    year = {2022},
    publisher={CUP}, 
    note = {\url{https://pi.math.cornell.edu/~hatcher/TN/TNpage.html}}
}

@Article{Lachaud_cont-frac-zeta_1988,
    author={Lachaud, Gilles},
    title={Continued fractions, binary quadratic forms, quadratic fields, and zeta functions},
    journal={Algebra and topology (Taejon)},
    year={1988},
    volume={Korea Inst. Tech.},
    pages={1--56}
}

@phdthesis{CLS_phdthesis_2022,
  TITLE = {{Arithmetic and Topology of Modular knots}},
  AUTHOR = {Simon, Christopher-Lloyd},
  URL = {https://tel.archives-ouvertes.fr/tel-03755147},
  SCHOOL = {{Universit{\'e} de Lille}},
  YEAR = {2022},
  MONTH = Jun,
  TYPE = {Th{\`e}se},
  HAL_ID = {tel-03755147},
  HAL_VERSION = {v1}, 
  NOTE = {\href{https://tel.archives-ouvertes.fr/tel-03755147/file/AriTopoModuKnots_14-07-2022.pdf}{PDF on HAL}}
}

@Misc{CLS_Conj-PSL2K_2022,
    author = {Simon, Christopher-Lloyd},
    title = {Conjugacy classes in $\operatorname{PSL}_2(\mathbb{K})$},
    note = {Submitted for publication, \href{https://arxiv.org/abs/2210.02481}{arxiv}},
    year = {2022}
}

@article{Aka_Continued-fractions-quadratics_2020,
    AUTHOR = {Aka, Menny},
     TITLE = {Continued fractions of arithmetic sequences of quadratics},
   JOURNAL = {Expo. Math.},
  FJOURNAL = {Expositiones Mathematicae},
    VOLUME = {38},
      YEAR = {2020},
    NUMBER = {3},
     PAGES = {397--406},
      ISSN = {0723-0869,1878-0792},
   MRCLASS = {11K55 (11A55 11J70)},
  MRNUMBER = {4149180},
MRREVIEWER = {Thomas\ Ward},
       DOI = {10.1016/j.exmath.2019.05.004},
       URL = {https://doi.org/10.1016/j.exmath.2019.05.004},
}

@article{Aka-Shapira_Continued-fractions-quadratics_2018,
    AUTHOR = {Aka, Menny and Shapira, Uri},
     TITLE = {On the evolution of continued fractions in a fixed quadratic
              field},
   JOURNAL = {J. Anal. Math.},
  FJOURNAL = {Journal d'Analyse Math\'{e}matique},
    VOLUME = {134},
      YEAR = {2018},
    NUMBER = {1},
     PAGES = {335--397},
      ISSN = {0021-7670,1565-8538},
   MRCLASS = {37A45 (11A55 30B70)},
  MRNUMBER = {3771486},
MRREVIEWER = {Derong\ Kong},
       DOI = {10.1007/s11854-018-0012-4},
       URL = {https://doi.org/10.1007/s11854-018-0012-4},
}

@article {Shallit_Reals-bounded-partial-quotients:survey_1992,
    AUTHOR = {Shallit, Jeffrey},
     TITLE = {Real numbers with bounded partial quotients: a survey},
   JOURNAL = {Enseign. Math. (2)},
  FJOURNAL = {L'Enseignement Math\'{e}matique. Revue Internationale. 2e
              S\'{e}rie},
    VOLUME = {38},
      YEAR = {1992},
    NUMBER = {1-2},
     PAGES = {151--187},
      ISSN = {0013-8584},
   MRCLASS = {11A55 (11J70)},
  MRNUMBER = {1175525},
MRREVIEWER = {Douglas\ Hensley},
}

@article{GPS_asymptotic-index_2001,
    AUTHOR = {Giordano, Thierry and Putnam, Ian F. and Skau, Christian F.},
     TITLE = {{$K$}-theory and asymptotic index for certain almost
              one-to-one factors},
   JOURNAL = {Math. Scand.},
  FJOURNAL = {Mathematica Scandinavica},
    VOLUME = {89},
      YEAR = {2001},
    NUMBER = {2},
     PAGES = {297--319},
      ISSN = {0025-5521,1903-1807},
   MRCLASS = {46L80 (19K99 37B05 46L35 54F45 54H20)},
  MRNUMBER = {1868179},
MRREVIEWER = {Vicumpriya\ S.\ Perera},
       DOI = {10.7146/math.scand.a-14343},
       URL = {https://doi.org/10.7146/math.scand.a-14343},
}

@article{Boyle-Chuysurichay_MCG-subshift-finite-type_2018,
    AUTHOR = {Boyle, Mike and Chuysurichay, Sompong},
     TITLE = {The mapping class group of a shift of finite type},
   JOURNAL = {J. Mod. Dyn.},
  FJOURNAL = {Journal of Modern Dynamics},
    VOLUME = {13},
      YEAR = {2018},
     PAGES = {115--145},
      ISSN = {1930-5311},
   MRCLASS = {37B10 (20F10 20F38)},
  MRNUMBER = {3918260},
MRREVIEWER = {Alfredo M. G. Costa},
       DOI = {10.3934/jmd.2018014},
       URL = {https://doi.org/10.3934/jmd.2018014},
}

@article{Boyle-Carlsen-Eiler_Flow-isotopy_2017,
    AUTHOR = {Boyle, Mike and Carlsen, Toke Meier and Eilers, S{\o}ren},
     TITLE = {Flow equivalence and isotopy for subshifts},
   JOURNAL = {Dyn. Syst.},
  FJOURNAL = {Dynamical Systems. An International Journal},
    VOLUME = {32},
      YEAR = {2017},
    NUMBER = {3},
     PAGES = {305--325},
      ISSN = {1468-9367},
   MRCLASS = {37B10 (37B50)},
  MRNUMBER = {3669803},
MRREVIEWER = {Kengo Matsumoto},
       DOI = {10.1080/14689367.2016.1207753},
       URL = {https://doi-org.proxy-um.researchport.umd.edu/10.1080/14689367.2016.1207753},
        note = {Corrigendum: ``{F}low equivalence and isotopy for subshifts''},
}

@article{Parry-Sullivan_topo-invariant-flows_1975,
	Author = {Parry, Bill and Sullivan, Dennis},
	Fjournal = {Topology. An International Journal of Mathematics},
	Issn = {0040-9383},
	Journal = {Topology},
	Mrclass = {54H20 (58F99)},
	Mrnumber = {MR0405385 (53 \#9179)},
	Mrreviewer = {Gustav A. Hedlund},
	Number = {4},
	Pages = {297--299},
	Title = {A topological invariant of flows on {$1$}-dimensional spaces},
	Volume = {14},
	Year = {1975},
}

@article{Putnam_Cstar-algebras-IET_1992,
    AUTHOR = {Putnam, Ian F.},
     TITLE = {{$C^*$}-algebras arising from interval exchange
              transformations},
   JOURNAL = {J. Operator Theory},
  FJOURNAL = {Journal of Operator Theory},
    VOLUME = {27},
      YEAR = {1992},
    NUMBER = {2},
     PAGES = {231--250},
      ISSN = {0379-4024},
   MRCLASS = {46L55 (28D05 46L80)},
  MRNUMBER = {1249645},
MRREVIEWER = {Dana P. Williams},
}

@article{Putnam_Cstar-algebras-minimal-homeo-Cantor_1989,
    AUTHOR = {Putnam, Ian F.},
     TITLE = {The {$C^*$}-algebras associated with minimal homeomorphisms of
              the {C}antor set},
   JOURNAL = {Pacific J. Math.},
  FJOURNAL = {Pacific Journal of Mathematics},
    VOLUME = {136},
      YEAR = {1989},
    NUMBER = {2},
     PAGES = {329--353},
      ISSN = {0030-8730},
   MRCLASS = {46L80 (19K99 46L55 58F11)},
  MRNUMBER = {978619},
MRREVIEWER = {Klaus Thomsen},
       URL = {http://projecteuclid.org/euclid.pjm/1102650733},
}

@article{PSS_C-star-Denjoy_1986,
    AUTHOR = {Putnam, Ian and Schmidt, Klaus and Skau, Christian},
     TITLE = {{$C^\ast$}-algebras associated with {D}enjoy homeomorphisms of
              the circle},
   JOURNAL = {J. Operator Theory},
  FJOURNAL = {Journal of Operator Theory},
    VOLUME = {16},
      YEAR = {1986},
    NUMBER = {1},
     PAGES = {99--126},
      ISSN = {0379-4024},
   MRCLASS = {46L55 (46L80 58F11 58F22)},
  MRNUMBER = {847334},
MRREVIEWER = {Thierry\ Fack},
}

@book{Fogg_Subsitutions-in-dyn-arit-comb_2002,
    AUTHOR = {Fogg, N. Pytheas},
     TITLE = {Substitutions in dynamics, arithmetics and combinatorics},
    SERIES = {Lecture Notes in Mathematics},
    VOLUME = {1794},
      NOTE = {Edited by V. Berth\'{e}, S. Ferenczi, C. Mauduit and A. Siegel},
 PUBLISHER = {Springer-Verlag, Berlin},
      YEAR = {2002},
     PAGES = {xviii+402},
      ISBN = {3-540-44141-7},
   MRCLASS = {37-06 (05A99 11B85 11J70 11R06 37B10 68R15)},
  MRNUMBER = {1970385},
MRREVIEWER = {Teturo Kamae},
       DOI = {10.1007/b13861},
       URL = {https://doi.org/10.1007/b13861},
}

@article{Durand-Ormes-Petite_self-induced-systems_2018,
    AUTHOR = {Durand, Fabien and Ormes, Nicholas and Petite, Samuel},
     TITLE = {Self-induced systems},
   JOURNAL = {J. Anal. Math.},
  FJOURNAL = {Journal d'Analyse Math\'{e}matique},
    VOLUME = {135},
      YEAR = {2018},
    NUMBER = {2},
     PAGES = {725--756},
      ISSN = {0021-7670},
   MRCLASS = {37B20 (37B05 37B10 37B40 54H20)},
  MRNUMBER = {3829615},
MRREVIEWER = {Antonio Linero Bas},
       DOI = {10.1007/s11854-018-0051-x},
       URL = {https://doi.org/10.1007/s11854-018-0051-x},
}

@article{Schmieding-Yang_Map-subshift_2021,
    AUTHOR = {Schmieding, Scott and Yang, Kitty},
     TITLE = {The mapping class group of a minimal subshift},
   JOURNAL = {Colloq. Math.},
  FJOURNAL = {Colloquium Mathematicum},
    VOLUME = {163},
      YEAR = {2021},
    NUMBER = {2},
     PAGES = {233--265},
      ISSN = {0010-1354},
   MRCLASS = {37B10 (20F38 37A20)},
  MRNUMBER = {4173180},
MRREVIEWER = {Anna Giordano Bruno},
       DOI = {10.4064/cm7933-2-2020},
       URL = {https://doi.org/10.4064/cm7933-2-2020},
}

@article{Putnam_Morita-equivalence-Denjoy_1988,
    AUTHOR = {Putnam, Ian F.},
     TITLE = {Strong {M}orita equivalence for the {D}enjoy {$C^*$}-algebras},
   JOURNAL = {Canad. Math. Bull.},
  FJOURNAL = {Canadian Mathematical Bulletin. Bulletin Canadien de
              Math\'{e}matiques},
    VOLUME = {31},
      YEAR = {1988},
    NUMBER = {4},
     PAGES = {439--447},
      ISSN = {0008-4395},
   MRCLASS = {46L55 (46L05 46L80 54H20)},
  MRNUMBER = {971571},
MRREVIEWER = {Judith A. Packer},
       DOI = {10.4153/CMB-1988-064-7},
       URL = {https://doi.org/10.4153/CMB-1988-064-7},
}

@book{Fokkink_structure-trajectories_1991,
    AUTHOR = {Fokkink, Robbert Johan},
     TITLE = {The structure of trajectories},
      NOTE = {Thesis (Ph.D.)--Technische Universiteit Delft (The
              Netherlands)},
 PUBLISHER = {ProQuest LLC, Ann Arbor, MI},
      YEAR = {1991},
     PAGES = {112},
   MRCLASS = {Thesis},
  MRNUMBER = {2714473},
       URL =
              {http://gateway.proquest.com/openurl?url_ver=Z39.88-2004&rft_val_fmt=info:ofi/fmt:kev:mtx:dissertation&res_dat=xri:pqdiss&rft_dat=xri:pqdiss:C211667},
}

@article{Barge-Williams_Denjoy-classification_2000,
    AUTHOR = {Barge, Marcy and Williams, R. F.},
     TITLE = {Classification of {D}enjoy continua},
   JOURNAL = {Topology Appl.},
  FJOURNAL = {Topology and its Applications},
    VOLUME = {106},
      YEAR = {2000},
    NUMBER = {1},
     PAGES = {77--89},
      ISSN = {0166-8641},
   MRCLASS = {37B45 (54F15 54H20)},
  MRNUMBER = {1769334},
MRREVIEWER = {Wayne Lewis},
       DOI = {10.1016/S0166-8641(99)00070-X},
       URL = {https://doi.org/10.1016/S0166-8641(99)00070-X},
}

@article{Boyle-Handelman_orbit-flow-oredered-cohomology_1996,
    AUTHOR = {Boyle, Mike and Handelman, David},
     TITLE = {Orbit equivalence, flow equivalence and ordered cohomology},
   JOURNAL = {Israel J. Math.},
  FJOURNAL = {Israel Journal of Mathematics},
    VOLUME = {95},
      YEAR = {1996},
     PAGES = {169--210},
      ISSN = {0021-2172},
   MRCLASS = {46L55 (46L85 54H20 55M10)},
  MRNUMBER = {1418293},
MRREVIEWER = {A. I. Danilenko},
       DOI = {10.1007/BF02761039},
       URL = {https://doi.org/10.1007/BF02761039},
}

@article{ItzaOrtiz_Denjoy-flow_2015,
    AUTHOR = {Itz\'{a} Ortiz, Benjam\'{\i}n A.},
     TITLE = {Classification of generalized {D}enjoy continua},
   JOURNAL = {Houston J. Math.},
  FJOURNAL = {Houston Journal of Mathematics},
    VOLUME = {41},
      YEAR = {2015},
    NUMBER = {4},
     PAGES = {1295--1311},
      ISSN = {0362-1588},
   MRCLASS = {37B45 (46L55)},
  MRNUMBER = {3455359},
MRREVIEWER = {Jan P. Boro\'{n}ski},
}

@article{Dolce-Francesco_IET-admissibility-induction_2017,
    AUTHOR = {Dolce, Francesco and Perrin, Dominique},
     TITLE = {Interval exchanges, admissibility and branching {R}auzy
              induction},
   JOURNAL = {RAIRO Theor. Inform. Appl.},
  FJOURNAL = {RAIRO Theoretical Informatics and Applications. Informatique
              Th\'eorique et Applications},
    VOLUME = {51},
      YEAR = {2017},
    NUMBER = {3},
     PAGES = {141--166},
      ISSN = {0988-3754,1290-385X},
   MRCLASS = {37E05 (37B10 68R15)},
  MRNUMBER = {3743415},
MRREVIEWER = {Abdumajid\ Safarovich\ Begmatov},
       DOI = {10.1051/ita/2017004},
       URL = {https://doi.org/10.1051/ita/2017004},
}

@incollection{Yoccoz_IET-translation-surfaces_2010,
    AUTHOR = {Yoccoz, Jean-Christophe},
     TITLE = {Interval exchange maps and translation surfaces},
 BOOKTITLE = {Homogeneous flows, moduli spaces and arithmetic},
    SERIES = {Clay Math. Proc.},
    VOLUME = {10},
     PAGES = {1--69},
 PUBLISHER = {Amer. Math. Soc., Providence, RI},
      YEAR = {2010},
      ISBN = {978-0-8218-4742-8},
   MRCLASS = {37-02 (28D05 37A25 37C40 37D25 37D40 37D50)},
  MRNUMBER = {2648692},
MRREVIEWER = {Thomas\ A.\ Schmidt},
}

@article{Boshernitzan-Carroll_quadratic-IET_1997,
    AUTHOR = {Boshernitzan, M. D. and Carroll, C. R.},
     TITLE = {An extension of {L}agrange's theorem to interval exchange
              transformations over quadratic fields},
   JOURNAL = {J. Anal. Math.},
  FJOURNAL = {Journal d'Analyse Math\'ematique},
    VOLUME = {72},
      YEAR = {1997},
     PAGES = {21--44},
      ISSN = {0021-7670,1565-8538},
   MRCLASS = {58F11 (11R45 58F03)},
  MRNUMBER = {1482988},
MRREVIEWER = {P.\ Mendes},
       DOI = {10.1007/BF02843152},
       URL = {https://doi.org/10.1007/BF02843152},
}

@article{Masur_IET-MF_1982,
    AUTHOR = {Masur, Howard},
     TITLE = {Interval exchange transformations and measured foliations},
   JOURNAL = {Ann. of Math. (2)},
  FJOURNAL = {Annals of Mathematics. Second Series},
    VOLUME = {115},
      YEAR = {1982},
    NUMBER = {1},
     PAGES = {169--200},
      ISSN = {0003-486X},
   MRCLASS = {28D05 (10K10 58F11 58F18)},
  MRNUMBER = {644018},
MRREVIEWER = {D.\ Newton},
       DOI = {10.2307/1971341},
       URL = {https://doi.org/10.2307/1971341},
}

@article{Veech_Gauss-measure-induction-IET-ae-UE_1982,
    AUTHOR = {Veech, William A.},
     TITLE = {Gauss measures for transformations on the space of interval
              exchange maps},
   JOURNAL = {Ann. of Math. (2)},
  FJOURNAL = {Annals of Mathematics. Second Series},
    VOLUME = {115},
      YEAR = {1982},
    NUMBER = {1},
     PAGES = {201--242},
      ISSN = {0003-486X},
   MRCLASS = {28D05 (58F11)},
  MRNUMBER = {644019},
MRREVIEWER = {Michael\ Keane},
       DOI = {10.2307/1971391},
       URL = {https://doi.org/10.2307/1971391},
}

@article{Ferenczi-Zamboni_Languages-k-IET_2008,
    AUTHOR = {Ferenczi, S\'ebastien and Zamboni, Luca Q.},
     TITLE = {Languages of {$k$}-interval exchange transformations},
   JOURNAL = {Bull. Lond. Math. Soc.},
  FJOURNAL = {Bulletin of the London Mathematical Society},
    VOLUME = {40},
      YEAR = {2008},
    NUMBER = {4},
     PAGES = {705--714},
      ISSN = {0024-6093,1469-2120},
   MRCLASS = {37B10 (68R15)},
  MRNUMBER = {2441143},
MRREVIEWER = {Val\'erie\ Berth\'e},
       DOI = {10.1112/blms/bdn051},
       URL = {https://doi.org/10.1112/blms/bdn051},
}

@article{Barge-Kallendonk-Schmieding_factors-tiling_2012,
    AUTHOR = {Barge, Marcy and Kellendonk, Johannes and Schmieding, Scott},
     TITLE = {Maximal equicontinuous factors and cohomology for tiling
              spaces},
   JOURNAL = {Fund. Math.},
  FJOURNAL = {Fundamenta Mathematicae},
    VOLUME = {218},
      YEAR = {2012},
    NUMBER = {3},
     PAGES = {243--268},
      ISSN = {0016-2736},
   MRCLASS = {37Bxx (54H11)},
  MRNUMBER = {2982777},
       DOI = {10.4064/fm218-3-3},
       URL = {https://doi.org/10.4064/fm218-3-3},
}

@article{Masui_Denjoy_2009,
    AUTHOR = {Masui, Kenichi},
     TITLE = {Denjoy systems and substitutions},
   JOURNAL = {Tokyo J. Math.},
  FJOURNAL = {Tokyo Journal of Mathematics},
    VOLUME = {32},
      YEAR = {2009},
    NUMBER = {1},
     PAGES = {33--53},
      ISSN = {0387-3870},
   MRCLASS = {37B10 (05A05)},
  MRNUMBER = {2541153},
MRREVIEWER = {Clemens Fuchs},
       DOI = {10.3836/tjm/1249648408},
       URL = {https://doi.org/10.3836/tjm/1249648408},
}

@article{McMullen_diophantine-quadratic-field_2009,
    AUTHOR = {McMullen, Curtis T.},
     TITLE = {Uniformly {D}iophantine numbers in a fixed real quadratic
              field},
   JOURNAL = {Compos. Math.},
  FJOURNAL = {Compositio Mathematica},
    VOLUME = {145},
      YEAR = {2009},
    NUMBER = {4},
     PAGES = {827--844},
      ISSN = {0010-437X,1570-5846},
   MRCLASS = {11J17 (11J70 37A45)},
  MRNUMBER = {2521246},
MRREVIEWER = {Jeffrey\ O.\ Shallit},
       DOI = {10.1112/S0010437X09004102},
       URL = {https://doi.org/10.1112/S0010437X09004102},
}

@article{Dahmani_Groups-IET_2019,
    AUTHOR = {Dahmani, Fran\c cois},
     TITLE = {Groups of interval exchange transformations},
   JOURNAL = {Winter Braids Lect. Notes},
  FJOURNAL = {Winter Braids Lecture Notes},
    VOLUME = {6},
      YEAR = {2019},
     PAGES = {Exp. No. I, 22},
      ISSN = {2426-0312},
   MRCLASS = {37E05 (20B27)},
  MRNUMBER = {4374308},
}

@article{Keane_IET_1975,
    AUTHOR = {Keane, Michael},
     TITLE = {Interval exchange transformations},
   JOURNAL = {Math. Z.},
  FJOURNAL = {Mathematische Zeitschrift},
    VOLUME = {141},
      YEAR = {1975},
     PAGES = {25--31},
      ISSN = {0025-5874,1432-1823},
   MRCLASS = {28A65},
  MRNUMBER = {357739},
MRREVIEWER = {M.\ Lo\`eve},
       DOI = {10.1007/BF01236981},
       URL = {https://doi.org/10.1007/BF01236981},
}

@book {AransonBelitskyZhuzhoma1996,
    AUTHOR = {Aranson, S. Kh. and Belitsky, G. R. and Zhuzhoma, E. V.},
     TITLE = {Introduction to the qualitative theory of dynamical systems on
              surfaces},
    SERIES = {Translations of Mathematical Monographs},
    VOLUME = {153},
      NOTE = {Translated from the Russian manuscript by H. H. McFaden},
 PUBLISHER = {American Mathematical Society, Providence, RI},
      YEAR = {1996},
     PAGES = {xiv+325},
      ISBN = {0-8218-0369-7},
   MRCLASS = {58F25 (34C28 34C35)},
  MRNUMBER = {1400885},
MRREVIEWER = {Michael Hurley},
       DOI = {10.1090/mmono/153},
       URL = {https://doi.org/10.1090/mmono/153},
}

@article {HedlundMorse1940,
    AUTHOR = {Morse, Marston and Hedlund, Gustav A.},
     TITLE = {Symbolic dynamics {II}. {S}turmian trajectories},
   JOURNAL = {Amer. J. Math.},
  FJOURNAL = {American Journal of Mathematics},
    VOLUME = {62},
      YEAR = {1940},
     PAGES = {1--42},
      ISSN = {0002-9327},
   MRCLASS = {70.1X},
  MRNUMBER = {745},
MRREVIEWER = {C. B. Tompkins},
       DOI = {10.2307/2371431},
       URL = {https://doi.org/10.2307/2371431},
}

@article {Glasner-Weiss-1995,
    AUTHOR = {Glasner, Eli and Weiss, Benjamin},
     TITLE = {Weak orbit equivalence of {C}antor minimal systems},
   JOURNAL = {Internat. J. Math.},
  FJOURNAL = {International Journal of Mathematics},
    VOLUME = {6},
      YEAR = {1995},
    NUMBER = {4},
     PAGES = {559--579},
      ISSN = {0129-167X},
   MRCLASS = {46L55 (54H20)},
  MRNUMBER = {1339645},
MRREVIEWER = {A. I. Danilenko},
       DOI = {10.1142/S0129167X95000213},
       URL = {https://doi.org/10.1142/S0129167X95000213},
}

@book {DenkerGrillenbergerSigmundBook,
    AUTHOR = {Denker, Manfred and Grillenberger, Christian and Sigmund,
              Karl},
     TITLE = {Ergodic theory on compact spaces},
    SERIES = {Lecture Notes in Mathematics, Vol. 527},
 PUBLISHER = {Springer-Verlag, Berlin-New York},
      YEAR = {1976},
     PAGES = {iv+360},
   MRCLASS = {28A65 (54H20)},
  MRNUMBER = {457675},
MRREVIEWER = {W. Parry},
}

@article{Veech_metric-theory-IET-1-Spectral_1984,
    AUTHOR = {Veech, William A.},
     TITLE = {The metric theory of interval exchange transformations. {I}.
              {G}eneric spectral properties},
   JOURNAL = {Amer. J. Math.},
  FJOURNAL = {American Journal of Mathematics},
    VOLUME = {106},
      YEAR = {1984},
    NUMBER = {6},
     PAGES = {1331--1359},
      ISSN = {0002-9327,1080-6377},
   MRCLASS = {28D05 (28D20 58F11)},
  MRNUMBER = {765582},
MRREVIEWER = {Andr\'es\ del\ Junco},
       DOI = {10.2307/2374396},
       URL = {https://doi.org/10.2307/2374396},
}

@article{Veech_metric-theory-IET-2-ApproxByPrimitive_1984,
    AUTHOR = {Veech, William A.},
     TITLE = {The metric theory of interval exchange transformations. {II}.
              {A}pproximation by primitive interval exchanges},
   JOURNAL = {Amer. J. Math.},
  FJOURNAL = {American Journal of Mathematics},
    VOLUME = {106},
      YEAR = {1984},
    NUMBER = {6},
     PAGES = {1361--1387},
      ISSN = {0002-9327,1080-6377},
   MRCLASS = {28D05 (58D20 58F11)},
  MRNUMBER = {765583},
MRREVIEWER = {Andr\'es\ del\ Junco},
       DOI = {10.2307/2374397},
       URL = {https://doi.org/10.2307/2374397},
}

@article{Veech_metric-theory-IET-3-SAF_1984,
    AUTHOR = {Veech, William A.},
     TITLE = {The metric theory of interval exchange transformations. {III}.
              {T}he {S}ah-{A}rnoux-{F}athi invariant},
   JOURNAL = {Amer. J. Math.},
  FJOURNAL = {American Journal of Mathematics},
    VOLUME = {106},
      YEAR = {1984},
    NUMBER = {6},
     PAGES = {1389--1422},
      ISSN = {0002-9327,1080-6377},
   MRCLASS = {28D05 (58D20 58F11)},
  MRNUMBER = {765584},
MRREVIEWER = {Andr\'es\ del\ Junco},
       DOI = {10.2307/2374398},
       URL = {https://doi.org/10.2307/2374398},
}

@article {Boyle_shift-equiv-implies-flow-equiv_2024,
    AUTHOR = {Boyle, Mike},
     TITLE = {Shift equivalence implies flow equivalence for shifts of
              finite type},
   JOURNAL = {Illinois J. Math.},
  FJOURNAL = {Illinois Journal of Mathematics},
    VOLUME = {70},
      YEAR = {2026},
    NUMBER = {1},
     PAGES = {33--51},
      ISSN = {0019-2082,1945-6581},
   MRCLASS = {37B10 (15B36)},
  MRNUMBER = {5076497},
       DOI = {10.1215/00192082-12306614},
       URL = {https://doi.org/10.1215/00192082-12306614},
}

@incollection{Boyle_algebraic-aspects-symbolic-dynamics_2000,
    AUTHOR = {Boyle, Mike},
     TITLE = {Algebraic aspects of symbolic dynamics},
 BOOKTITLE = {Topics in symbolic dynamics and applications ({T}emuco, 1997)},
    SERIES = {London Math. Soc. Lecture Note Ser.},
    VOLUME = {279},
     PAGES = {57--88},
 PUBLISHER = {Cambridge Univ. Press, Cambridge},
      YEAR = {2000},
   MRCLASS = {37B10},
  MRNUMBER = {1776756},
}

@ARTICLE{Schmieding-Simon_hexp_2024,
       author = {{Schmieding}, Scott and {Simon}, Christopher-Lloyd},
        title = "{Geometry and Transcendence of the Hexponential}",
      journal = {arXiv e-prints},
         year = 2024,
        month = feb,
          eid = {arXiv:2402.17628},
        pages = {arXiv:2402.17628},
          doi = {10.48550/arXiv.2402.17628},
archivePrefix = {arXiv},
       eprint = {2402.17628},
 primaryClass = {math.NT},
       adsurl = {https://ui.adsabs.harvard.edu/abs/2024arXiv240217628S}
}

@article {Markley_homeo-circle_1970,
    AUTHOR = {Markley, Nelson G.},
     TITLE = {Homeomorphisms of the circle without periodic points},
   JOURNAL = {Proc. London Math. Soc. (3)},
  FJOURNAL = {Proceedings of the London Mathematical Society. Third Series},
    VOLUME = {20},
      YEAR = {1970},
     PAGES = {688--698},
      ISSN = {0024-6115},
   MRCLASS = {54.80},
  MRNUMBER = {268873},
MRREVIEWER = {E. Hemmingsen},
       DOI = {10.1112/plms/s3-20.4.688},
       URL = {https://doi.org/10.1112/plms/s3-20.4.688},
}

@article{Poincare_courbes-equa-diff_1885,
author = {H. Poincaré},
journal = {Journal de Mathématiques Pures et Appliquées},
language = {fre},
pages = {167-244},
title = {Sur les courbes définies par les équations différentielles (III)},
url = {http://eudml.org/doc/235596},
volume = {1},
year = {1885},
}

@article{Stanley_Smith-Normal-Form-combinatorics_2016,
    AUTHOR = {Stanley, Richard P.},
     TITLE = {Smith normal form in combinatorics},
   JOURNAL = {J. Combin. Theory Ser. A},
  FJOURNAL = {Journal of Combinatorial Theory. Series A},
    VOLUME = {144},
      YEAR = {2016},
     PAGES = {476--495},
      ISSN = {0097-3165,1096-0899},
   MRCLASS = {15-02 (05E18 15A21)},
  MRNUMBER = {3534076},
       DOI = {10.1016/j.jcta.2016.06.013},
       URL = {https://doi.org/10.1016/j.jcta.2016.06.013},
}

@article{Durand_characterization-substitutive-sequences-return-words_1998,
    AUTHOR = {Durand, Fabien},
     TITLE = {A characterization of substitutive sequences using return
              words},
   JOURNAL = {Discrete Math.},
  FJOURNAL = {Discrete Mathematics},
    VOLUME = {179},
      YEAR = {1998},
    NUMBER = {1-3},
     PAGES = {89--101},
      ISSN = {0012-365X,1872-681X},
   MRCLASS = {68R15 (68Q68)},
  MRNUMBER = {1489074},
MRREVIEWER = {Val\'erie\ Berth\'e},
       DOI = {10.1016/S0012-365X(97)00029-0},
       URL = {https://doi.org/10.1016/S0012-365X(97)00029-0},
}

@incollection{Berthe-Delecroix_beyond-substitutive-S-adic_2014,
    AUTHOR = {Berth\'e, Val\'erie and Delecroix, Vincent},
     TITLE = {Beyond substitutive dynamical systems: {$S$}-adic expansions},
 BOOKTITLE = {Numeration and substitution 2012},
    SERIES = {RIMS K\^oky\^uroku Bessatsu},
    VOLUME = {B46},
     PAGES = {81--123},
 PUBLISHER = {Res. Inst. Math. Sci. (RIMS), Kyoto},
      YEAR = {2014},
   MRCLASS = {37B10 (37B25)},
  MRNUMBER = {3330561},
MRREVIEWER = {Paul\ Surer},
}

@article{Durand-Host-Skau_Substitutional-systems-Bratteli-dimension-groups_1999,
    AUTHOR = {Durand, F. and Host, B. and Skau, C.},
     TITLE = {Substitutional dynamical systems, {B}ratteli diagrams and
              dimension groups},
   JOURNAL = {Ergodic Theory Dynam. Systems},
  FJOURNAL = {Ergodic Theory and Dynamical Systems},
    VOLUME = {19},
      YEAR = {1999},
    NUMBER = {4},
     PAGES = {953--993},
      ISSN = {0143-3857,1469-4417},
   MRCLASS = {46L55 (37B99 46L80 54H20)},
  MRNUMBER = {1709427},
MRREVIEWER = {A.\ I.\ Danilenko},
       DOI = {10.1017/S0143385799133947},
       URL = {https://doi.org/10.1017/S0143385799133947},
}

@misc{Delecroix_translation-surfaces-Salta_2015,
    title = {Interval Exchange transformations},
    author = {Delecroix, Vincent},
    year = {2015},
    note = {Salta Lecture notes},
}

@incollection{Boyle-Schmieding_Stable-Algebra_2024,
    AUTHOR = {Boyle, Mike and Schmieding, Scott},
     TITLE = {Symbolic dynamics and the stable algebra of matrices},
 BOOKTITLE = {Groups and graphs, designs and dynamics},
    SERIES = {London Math. Soc. Lecture Note Ser.},
    VOLUME = {491},
     PAGES = {266--422},
 PUBLISHER = {Cambridge Univ. Press, Cambridge},
      YEAR = {2024},
   MRCLASS = {37B10 (15B36 15B48 19Bxx 19D10)},
  MRNUMBER = {4813459},
}

@article{Schinzel_power-residues-exponential-congruences_1975,
    AUTHOR = {Schinzel, A.},
     TITLE = {On power residues and exponential congruences},
   JOURNAL = {Acta Arith.},
  FJOURNAL = {Polska Akademia Nauk. Instytut Matematyczny. Acta Arithmetica},
    VOLUME = {27},
      YEAR = {1975},
     PAGES = {397--420},
      ISSN = {0065-1036},
   MRCLASS = {12A20},
  MRNUMBER = {379432},
MRREVIEWER = {D.\ J.\ Lewis},
       DOI = {10.4064/aa-27-1-397-420},
       URL = {https://doi.org/10.4064/aa-27-1-397-420},
}

@article{Boyle-Carlsen-Eilers_flow-equivalence-subshifts_2017,
    AUTHOR = {Boyle, M. and Carlsen, T. M. and Eilers, S.},
     TITLE = {Corrigendum: ``{F}low equivalence and isotopy for subshifts''
              [{MR}3669803]},
   JOURNAL = {Dyn. Syst.},
  FJOURNAL = {Dynamical Systems. An International Journal},
    VOLUME = {32},
      YEAR = {2017},
    NUMBER = {3},
     PAGES = {ii},
      ISSN = {1468-9367,1468-9375},
   MRCLASS = {37B50},
  MRNUMBER = {3669811},
       DOI = {10.1080/14689367.2017.1322774},
       URL = {https://doi.org/10.1080/14689367.2017.1322774},
}

@misc{boylepersonalcommunication,
  author = "Boyle, Mike",
  date = "2025-02-11",
  howpublished = "personal communication"
}

\end{document}